\documentclass[a4paper,Haag duality]{mathscan}
\newtheorem{thm}{Theorem}[section] 
\newtheorem{pro}[thm]{Proposition}

\theoremstyle{definition}           
\newtheorem{rem}[thm]{Remark}       %%%%% the counter [thm] is optional 
\newcommand{\NI}{\noindent}

\newcommand{\bea}{\begin{eqnarray}}
\newcommand{\eea}{\end{eqnarray}}
\newcommand{\dsp}{\displaystyle}
\def \b #1 {\bf #1}
\newcommand{\IC}{\mathbb{C}}
\newcommand{\IR}{\mathbb{R}}

\newcommand{\IE}{I\!\!E}

\newcommand{\IT}{I\!\!T}
\newcommand{\IZ}{\mathbb{Z}}
\newcommand{\cal}{\mathcal}
\newcommand{\clk}{{\cal K}}

\newcommand{\cla}{{\cal A}}
\newcommand{\clz}{{\cal Z}}
\newcommand{\cli}{{\cal I}}
\newcommand{\cls}{{\cal S}}

\newcommand{\clf}{{\cal F}}

\newcommand{\clh}{{\cal H}}
\newcommand{\clp}{{\cal P}}
\newcommand{\clo}{{\cal O}}
\newcommand{\clb}{{\cal B}}

\newcommand{\cle}{{\cal E}}
\newcommand{\clj}{{\cal J}}
\newcommand{\cln}{{\cal N}}

\newcommand{\clm}{{\cal M}}

\newcommand{\raro}{\rightarrow}

\newcommand{\vsp}{\vskip 1em}

\def \qed {\hfill \vrule height6pt width 6pt depth 0pt}
\newcommand{\be}{\begin{equation}}
\newcommand{\ee}{\end{equation}}
\newcommand{\ben}{\begin{eqnarray*}}
\newcommand{\een}{\end{eqnarray*}}
\begin{document}

\title{Translation invariant pure state on $\otimes_{\IZ}\!M_d(\IC)$ and Haag duality }

\author{ Anilesh Mohari }
\thanks{...}

\address{ The Institute of Mathematical Sciences, CIT Campus, Taramani, Chennai-600113 }

\email{anilesh@imsc.res.in}

\keywords{Uniformly hyperfinite factors. Cuntz algebra, Popescu dilation, Kolmogorov's property, 
Arveson's spectrum, Haag duality }

\subjclass{46L}

\thanks{ This paper has grown over the years starting with initial work in the 
middle of 2005. The author gratefully acknowledge discussion with Ola Bratteli and Palle E. T. Jorgensen 
for inspiring participation in sharing the intricacy of the present problem. Finally the author is indebted to 
Taku Matsui for valuable comments on an earlier draft of the present problem where the author made an attempt to prove 
Haag duality. }

\begin{abstract}

We prove Haag duality property of any translation invariant pure state on $\clb = \otimes_{\IZ}\!M_d(\IC), \;d \ge 2$, where $\!M_d(\IC)$ is 
the set of $d \times d$ dimensional matrices over the field of complex numbers. We also prove a necessary and sufficient condition for a 
translation invariant factor state to be pure on $\clb$. 

\end{abstract}

\maketitle 

\section{ Introduction }

\vsp 
A state $\omega$ on a $C^*$-algebra $\clb$ is called a factor if the center of the von-Neumann algebra $\pi_{\omega}(\clb)''$ is trivial, 
where $(\clh_{\omega},\pi_{\omega},\Omega)$ is the GNS space associated with $\omega$ on $\clb$ [BR vol-I]. A state $\omega$ on $\clb$ is called 
pure if $\pi_{\omega}(\clb)''=\clb(\clh_{\omega})$, the algebra of all bounded operators on $\clh_{\omega}$. Here we fix our convention that 
Hilbert spaces that are considered here equipped always with inner product $<.,.>$ which is linear in the second variable and conjugate linear in the first 
variable. In this paper our primary objective is to study states on $C^*$-algebra that naturally arise in quantum spin chain models on a 
lattice. 

\vsp 
Let $\clb=\otimes_{\IZ^k} \!M_d(\IC)$ be the uniformly hyper-finite $C^*$-algebra over the lattice $\IZ^k$ of dimension $k \ge 1$, where $\!M_d(\IC)$ denote the algebra of $d \times d$-matrices 
over the field of complex number $\IC$. A state $\omega$ on $\clb$ is called translation invariant if $\omega(x)=\omega(\theta_{\bar{m}}(x))$ where $\bar{m}=(m_1,m_2,.,m_k)$ and $\theta_{\bar{m}}$ 
is the translation induced by $\bar{\IZ} \raro  \bar{\IZ} + \bar{m}$ for all $\bar{z} \in \IZ^k$. It is well known since late 60's [Pow] that a translation invariant state $\omega$ 
on $\clb$ is a factor state if and only if 

\be 
\mbox{sup}_{x \in \clb_{\Lambda_n^c}, ||x|| \le 1}|\omega(xy)) - \omega(x)\omega(y)| \raro 0
\ee 
for all $y \in \clb$ as $n \raro \infty$, where $\Lambda_{n}$ is the local algebra with support in the finite set $\{\bar{m}: -n \le m_k \le n \}$. 

Such a criterion is used extensively to show that KMS states of a translation invariant Hamiltonian on the lattice form a simplex and 
its extreme points are translation invariant factor states. For more details and an account until 1980 we refer to [BR vol-II] and also [Sim vol-I] for a later edition. Such an elegant asymptotic 
criterion is missing for a translation invariant pure state. Here one of our objectives is to do so.  

\vsp 
For the sake of simplicity, we will consider the simplest situation namely one lattice dimensional quantum mechanical spin systems. Now onwards we consider the lattice 
to be one dimensional. We briefly set the standard notations and known relations in the following text. The quantum spin chain that
we consider here is described by a $\mbox{UHF}$ $C^*$-algebra denoted by $\clb=\otimes_{\IZ}M_d(\IC)$. Here $\clb$
is the $C^*$ -completion of the infinite tensor product of the algebra $M_d(\IC)$ of $d \times d$ complex matrices [Sa],
each component of the tensor product element is indexed by an integer $j$. Let $Q$ be a matrix in $M_d(\IC)$. We denote the element
$Q^{(j)}=...\otimes 1 \otimes 1 ... 1 \otimes Q \otimes 1 \otimes ... 1\otimes ,,. $, where $Q$
appears in the $j$-th component. Given a subset $\Lambda$ of $\IZ$, $\clb_{\Lambda}$ is defined as the $C^*$-sub-algebra
of $\clb$ generated by all $Q^{(j)}$ with $Q \in {\bf M}_d(\IC)$, $j \in \Lambda$. We also set
$$\clb_{loc}= \bigcup_{\Lambda:|\Lambda| < \infty } \clb_{\Lambda}$$
where $|\Lambda|$ is the cardinality of $\Lambda$. Let $\omega$ be a state on $\clb$. The restriction of $\omega$
to $\clb_{\Lambda}$ is denoted by $\omega_{\Lambda}$. We also set $\omega_{R}=\omega_{[1,\infty)}$ and $\omega_{L}=
\omega_{(-\infty,0]}$. The translation $\theta_k$ is an automorphism of $\clb$ defined by $\theta_k(Q^{(j)})=Q^{(j+k)}$.
Thus $\theta_1$ and $\theta_{-1}$ are unital $*$-endomorphisms on $\clb_R$ and $\clb_L$ respectively. We say $\omega$ is
translation invariant if $\omega \circ \theta_k = \omega$ on $\clb$ ( $\omega \circ \theta_1 = \omega$ on $\clb$ ).
In such a case $(\clb_R,\theta_1,\omega_{R})$ and $(\clb_L,\theta_{-1},\omega_{L})$ are two unital $*$-endomorphisms with
invariant states.

\vsp 
We will consider a Hamiltonian in one dimensional lattice of the following form
\be 
H= \sum_{ k \in \IZ} \theta^k(h_0)
\ee 
for $h^*_0=h_0 \in \clb_{loc}$ 
where the formal sum gives an auto-morphism $\alpha=(\alpha_t:t \in \IR)$ via the thermodynamic limit 
of $\alpha^{\Lambda}_t(x)=e^{itH_{\Lambda}}xe^{-itH_{\Lambda}}$ for a net of finite subsets of the lattice 
$\Lambda \uparrow \IZ$ whose surface energies are uniformly bounded, where $H_{\Lambda}=\sum_{k \in \Lambda} \theta^k(h_0)$ [Ru,BR2]. 
Such a thermodynamic limit automorphism $\alpha$ is uniquely determined by $H$. In such a case, i.e. translation invariant Hamiltonian $H$ having  
finite range interaction,  KMS state at a given inverse temperature exists and is unique [Ara1],[Ara2], [Ki] and inherits translation and other symmetry 
of the Hamiltonian. Thus low temperature limit points of unique KMS states give ground states for the Hamiltonian $H$ inheriting translation 
and other symmetry of Hamiltonian. It is a well known fact that ground states of a translation invariant Hamiltonian form a face in the convex set of states 
on $\clb$ and its extreme points are pure. In general ground states need to be unique and there are other non translation invariant ground states for a translation 
invariant Hamiltonian [Ma4]. Ising model admits non translation invariant ground states known as N\'{e}el state [BR2]. However, ground 
states that appear as low temperature limit of KMS states of a translation invariant Hamiltonian inherit translation and other symmetry (that we would consider in a follow up 
paper in more details) of the Hamiltonian. In particular if ground state for a translation invariant Hamiltonian model of type (2) is unique, then the ground state is a translation invariant 
pure state.  
    
\vsp 
Unlike classical spin chain problem, any translation invariant state $\omega$ on $\clb$ gives rise to a quantum Markov state in the sense of Luigi Accardi [Ac] and more specifically 
finite or infinitely correlated translation invariant state ( [FNW1], see also [BJ], [BJKW],[Mo3] for their natural generalization ). Here we briefly recall now explaining these two related 
concepts and explain the basic setup of the present problem and some difficulties that crop up. A detailed account is given in section 2 and then section 3 holds key results.   

\vsp
First we recall that the Cuntz algebra $\clo_d ( d \in \{2,3,.., \} )$ [Cun] is the universal unital $C^*$-algebra generated by the elements $\{s_1,s_2,...,
s_d \}$ subject to the following Cuntz relations:
\be 
s^*_is_j=\delta^i_j1,\;\;\sum_{1 \le i \le d } s_is^*_i=1
\ee

\vsp
There is a canonical action of the group $U(d)$ of unitary $d \times d$ matrices on $\clo_d$ given by
$$\beta_g(s_i)=\sum_{1 \le j \le d}\overline{g^j_i}s_j$$
for $g=(g^i_j) \in U(d)$. In particular the gauge action is defined by
$$\beta_z(s_i)=zs_i,\;\;z \in \IT =S^1= \{z \in \IC: |z|=1 \}.$$
If UHF$_d$ is the fixed point sub-algebra under the gauge action, then UHF$_d$ is the closure of the
linear span of all Wick ordered monomials of the form
$$s_{i_1}...s_{i_k}s^*_{j_k}...s^*_{j_1}$$
which is also isomorphic to the UHF$_d$ algebra
$$\clb_R =\otimes^{\infty}_1M_d(\IC)$$
so that the isomorphism carries the Wick ordered monomial above into the matrix element
$$e^{i_1}_{j_1}(1)\otimes e^{i_2}_{j_2}(2) \otimes....\otimes e^{i_k}_{j_k}(k) \otimes 1 \otimes 1 ....$$
and the restriction of $\beta_g$ to $UHF_d$ is then carried to action
$$Ad(g)\otimes Ad(g) \otimes Ad(g) \otimes ....$$

\vsp
We also define the canonical endomorphism $\lambda$ on $\clo_d$ by
\be 
\lambda(x)=\sum_{1 \le i \le d}s_ixs^*_i
\ee
The isomorphism carries $\lambda$ restricted to UHF$_d$ into the one-sided shift
$$y_1 \otimes y_2 \otimes ... \raro 1 \otimes y_1 \otimes y_2 ....$$
on $\otimes^{\infty}_1 M_d$. Note that $\lambda \beta_g = \beta_g \lambda $ on UHF$_d$.

\vsp
Let $d \in \{2,3,..,,..\}$ and $\IZ_d$ be a set of $d$ elements.  $\cli$ be the set of finite sequences
$I=(i_1,i_2,...,i_m)$ where
$i_k \in \IZ_d$ and $m \ge 1$. We also include empty set $\emptyset \in \cli$ and set $s_{\emptyset }=1=s^*_{\emptyset}$,
$s_{I}=s_{i_1}......s_{i_m} \in \clo_d $ and $s^*_{I}=s^*_{i_m}...s^*_{i_1} \in \clo_d$. 

\vsp 
We fix a translation invariant state $\omega$ on $\clb$ and denote by $\omega_R$ the restriction of $\omega$ to $\clb_R$. Using weak$^*$ compactness of 
the convex set of states on a $C^*$-algebra, a standard averaging method ensures that the set 
$$K_{\omega}=\{ \psi \in \cls(\clo_d): \psi \lambda= \psi,\; \psi_| \mbox{UHF}_d = \omega_R \}$$ 
is a non-empty compact subset of $\cls(\clo_d)$, where $\cls(\clo_d)$ is the weak$^*$ compact convex set of states on $\clo_d$. 
Further extremal elements in $K_{\omega}$ is a factor state if and only if $\omega_R$ is a factor state and any two
such extremal elements $\psi,\psi'$ are related by $\psi'=\psi \beta_z$ for some $z \in S^1=\{z \in \IC:\;|z|=1 \}$ by Lemma 7.4. in [BJKW] 
where $\beta_z(s_i)=zs_i$ is the automorphism on $\clo_d$ determined uniquely by universal property of Cuntz algebra. 

\vsp 
Irrespective of the factor property of $\omega$, we may choose an element $\psi$ of $K_{\omega}$ and consider the 
GNS space $(\clh,\pi,\Omega)$ associated with state $\psi$ on $\clo_d$.  We set 
$P \in \pi(\clo_d)''$ to be the support projection of $\psi$ i.e. $P=[\pi(\clo_d)'\Omega]$. 
Invariance property of the state $\psi=\psi \lambda$ will ensure that $P\Lambda(I-P)P=0$ where 
$$\Lambda(X)=\sum_i S_iXS_i^*$$ 
is the canonical endomorphism on $\pi_{\psi}(\clo_d)''$ with $S_i=\pi_{\psi}(s_i)$. This  
verifies that 
\be 
S_i^*P=PS^*_iP,\;1 \le i \le d
\ee 
We define a family of contractions $\{v_i:1 \le i \le d\}$ in $\clm$ by 
$v_i=PS_iP,1 \le i \le d $ where we set von-Neumann algebra 
$\clm=P\pi_{\phi}(\clo_d)''P$ acting on Hilbert subspace $\clk$ where $\clk$ is range of $P$. Thus we get $\clm = \{v_i: 1 \le i \le d \}''$ 
and a unital completely positive map $\tau(x)=P\Lambda(PxP)P=\sum_iv_ixv_i^*$ for all $x \in \clm$. Furthermore a crucial 
point to be noted that the support projection of $\psi$ in $\pi(\clo_d)''$ being equal to $P$, by our construction we have 
\be 
\{x \in \clb(\clk): \sum_i v_ixv_i^* = x \}= \clm'
\ee 

\vsp 
Conversely let $\clm$ be a von-Neumann algebra acting on a Hilbert space $\clk$. A family of contractions $\{v_i:1 \le i \le d \}$ in $\clm$ is called 
Popescu's elements if $\sum_iv_iv_i^*=1$. Given a Popescu's elements $\clp=\{\clk,\clm, v_i,\;1 \le i \le d, \sum_iv_iv^*_i=1 \}$, 
the map 
$$s_Is_J^* \raro v_Iv_J^*$$ 
is unital completely positive from $\clo_d$ to $\clm$ and thus Stinespring minimal dilation gives a representation $\pi:\clo_d \raro \clb(\clh)$, a Hilbert space $\clh$ with 
a projection $P$ with range equal to $\clk$ such that 
$$P\pi(s_i)^*P=\pi(s_i)^*P=v_i^*$$
and $\{\pi(s_I)\clk: |I| < \infty \}$ is total in $\clh$. For a faithful normal state $\phi$ on $\clm$
we define state $\psi$ on $\clo_d$ by 
\be 
\psi(s_Is_J^*)=\phi(v_Iv_J^*)
\ee 
The crucial point that we arrive at Proposition 2.4 that $P$ is the support projection for $\pi(\clo_d)''$ if and only if (6) holds. 

\vsp 
We verify also with $v_i^*=PS^*_iP$ that 
\be 
\omega_R(|e_{i_1}><e_{j_1}| \otimes |e_{i_2}><e_{j_2}| \otimes ...\otimes |e_{i_n}><e_{j_n}|)= \phi(v_Iv_J^*)
\ee       
where $v_I=v_{i_n}...v_{i_2}v_{i_1}$ and $v^*_J=v_{j_1}^*v^*_{j_2}...v^*_{j_n}$. The relation (7) can now be recast as 
a quantum Markov state as follows: Let $\clk$ be the Hilbert subspace $P$ of $\clh$ and $\clm$ be a von-Neumann sub-algebra of 
$\clb(\clk)$. Let $V^*:\clk \raro \clk \otimes \IC^d$ be an isometry and in an orthonormal basis $(e_i)$ for $\IC^d$, we have 
$V^*=(v_1^*,v_2^*,..,v_d^*)$ with $v_i \in \clm$. We define
$$\IE: \!M_d(\clm) \raro \clm$$ 
by 
$$\IE(X)=VXV^*=\sum_{1 \le i,j \le d} v_ix^i_jv_j^* $$ 
where $X=((x^i_j))$. Let $\phi$ be a state on $\clm$ such that 
$$\phi(\IE(B \otimes I_d))=\phi(B),\;\forall B \in \clm $$ For each $A \in M_d(\IC)$, define  
$\IE_{A}:\clm \raro \clm$ 
by $$B \raro \IE(B \otimes A)$$ Then 
\be 
\omega(A_1 \otimes A_1 \otimes A_3 ...\otimes A_m)=\phi(\IE_{A_1} \circ \IE_{A_2} \circ ..\IE_{A_m}(I_{\clk}))
\ee
defines a $\lambda$-invariant state on $\clb_R$ and the inductive limit state of $\clb_R \raro ^{\lambda} \clb_R$ [Sa] 
gives a translation invariant state on $\clb$. $\IE$ naturally gives a Markov map i.e. a unital completely positive map on $\clm$ defined by  
$\tau(x)=\IE( x \otimes I_d)=\sum_i v_ixv_i^*$. We have $\phi \tau =\phi$ on $\clm$. 

\vsp 
The state $\phi(x)=<\Omega,x \Omega>$ on $\clm$ being faithful and invariant of $\tau:\clm \raro \clm$ we find a unique unital completely positive map 
$\tilde{\tau}:\clm' \raro \clm'$ satisfying the duality  relation 
\be 
<y\Omega,\tau(x)\Omega>= <\tilde{\tau}(y)\Omega,x\Omega>
\ee
for all $x \in \clm$ and $y \in \clm'$, where $\clm'$ is the commutant of $\clm$ in $\clb(\clh)$. For a proof we refer to section 8 in the monograph [OP] or [Mo1]. 

\vsp
$\phi$ being also a faithful state, $\Omega \in \clk$ is a cyclic and separating vector for $\clm$ and 
the closure of the close-able operator $S_0:x\Omega \raro x^*\Omega, S$ possesses a polar decomposition
$S=\clj \Delta^{1/2}$, where $\clj$ is an anti-unitary and $\Delta$ is a non-negative self-adjoint operator 
on $\clk$. Tomita's [BR] theorem says that 
$\Delta^{it} \clm \Delta^{-it}=\clm,\;t \in \IR$ and $\clj \clm \clj=\clm'$, where $\clm'$ is the
commutant of $\clm$. We define the modular automorphism group
$\sigma=(\sigma_t,\;t \in \IT )$ on $\clm$
by
$$\sigma_t(x)=\Delta^{it}x\Delta^{-it}$$ which satisfies the modular relation
$$\phi(x\sigma_{-{i \over 2}}(y))=\phi(\sigma_{{i \over 2}}(y)x)$$
for any two analytic elements $x,y$ for the automorphism. A more useful 
form for modular  relation here
$$\phi(\sigma_{-{i \over 2}}(x^*)^* \sigma_{-{i \over 2}}(y^*))=\phi(y^*x)$$ 
which shows that $\clj x\Omega= \sigma_{-{i \over 2}}(x^*)\Omega$. 
$\clj$ and $\sigma=(\sigma_t,\;t \in \IR)$ are called Tomita's conjugation operator and 
modular automorphisms associated with $\phi$. Since $\tau(x)=v_kxv_k^*$ 
is an inner map i.e. each $v_k \in \clm$, we have an explicit formula for $\tilde{\tau}$ as follows. 
  
\vsp 
We set $\tilde{v}_k = \overline{ \clj \sigma_{i \over 2}(v^*_k) \clj } \in \clm'$. That $\tilde{v}_k$ is indeed well 
defined as an element in $\clm'$ given in section 8 in [BJKW]. By KMS or modular relation [BR vol-I] we verify that 
$$\sum_k \tilde{v}_k \tilde{v}_k^*=1$$ 
and
\be 
\tilde{\tau}(y)=\sum_k \tilde{v}_ky\tilde{v}^*_k
\ee
and 
\be
\phi(v_Iv^*_J)= \phi(\tilde{v}_{\tilde{I}}\tilde{v}^*_{\tilde{J}})
\ee
where $\tilde{I}=(i_n,..,i_2,i_1)$ if $I=(i_1,i_2,...,i_n)$. Moreover $\tilde{v}^*_I\Omega = 
\clj \sigma_{i \over 2}(v_{\tilde{I}})^*\clj\Omega= \clj \Delta^{1 \over 2}v_{\tilde{I}}\Omega
=v^*_{\tilde{I}}\Omega$. We also set $\tilde{\clm}$ to be the von-Neumann algebra generated by 
$\{\tilde{v}_k: 1 \le k \le d \}$. Thus $\tilde{\clm} \subseteq \clm'$. The major problem that we 
will address in the text when do we have the following equality: 
\be 
\{x \in \clb(\clk): \sum_k \tilde{v}_kx\tilde{v}_k^*= x \} =\clm
\ee 
Equality in (13) will ensure that $P:\tilde{\clh} \raro \clk$ is also the support projection of $\tilde{\pi}(\clo_d)''$ where 
$\tilde{\pi}$ is the Popescu's prescription of Stinespring representation $\tilde{\pi}: \clo_d \raro \clb(\tilde{\clh})$ 
associated with the completely positive map $s_Is_J^* \raro \tilde{v}_I\tilde{v}_J^*,\;|I|,|J| < \infty $ so that
$P\tilde{\pi}(s^*_i)P=\tilde{\pi}(s_i^*)P=\tilde{v}^*_i$. Details has been worked out in Proposition 2.4.

\vsp 
Thus so far we have taken an arbitrary element $\psi \in K_{\omega}$ and worked with its support projection to arrive at a representation of $\omega$ given in (7) or (8) by Popescu's elements 
$\clp=\{(\clk,v_i \in \clm,1 \le i \le d\;\Omega): \sum_kv_iv_i^*=I \}$. However by Lemma 7.4 in [BJKW] for a factor state $\omega$, if we choose an extreme point $\psi \in K_{\omega}$, 
two such extreme points $\psi$ and $\psi'$ in $K_{\omega}$ are related by $\psi'=\phi \beta_z$ for some $z \in S^1$, $\clp$ is uniquely determined modulo a unitary conjugation. In other words we find a 
one-one correspondence between 
\be
\omega \Leftrightarrow \omega_R \Leftrightarrow K^e_{\omega} \Leftrightarrow \clp_e 
\ee  
modulo unitary conjugation where $K^e_{\omega}$ denotes the set of extreme points in $K_{\omega}$ and $\clp_e$ the set of Popescu's elements associated with extreme points 
of $K_{\omega}$ on support projection of the state as described above. Further in such a case $\clm=\{v_k:1 \le k \le d\}''$ 
is a factor and $(\clm,\tau,\phi)$ is an ergodic {\it quantum dynamical system } [La,Ev]. A unital completely positive map $\tau$ on a von-Neumann algebra $\clm$ 
with an invariant normal state $\phi$ [BJKW,Mo1] is called {\it ergodic} if 
\be 
{1 \over N}\sum_{0 \le k \le N-1}\tau^k(x) \raro  \phi(x)I
\ee 
as $N \raro \infty$ in weak$^*$ topology for all $x \in \clm$. Thus any symmetry of $\omega$ will act on Popescu elements $\clp_e$ via this correspondence. It would be worthwhile 
to have a result generalizing this correspondence in a more general situation that Ruy Exel developed [Ex].

\vsp 
For a translation invariant factor state $\omega$ on $\clb$, we say it admits {\it Haag duality } property if 
\be 
\pi_{\omega}(\clb_R)'= \pi_{\omega}(\clb_L)''
\ee
It is clear that such a factor state is pure. A pure mathematical question that arises here whether converse is true? i.e. 
Do we always have Haag duality property for a translation invariant pure state of $\clb$?

\vsp 
In case $\pi_{\omega}(\clb_R)''$ is a type-I factor state, Haag duality property follows easily. In fact we can find Hilbert spaces $\clh^-_{\pi},\;\clh^+_{\pi}$ such that
$\clh_{\pi}$ is unitary equivalent to $\clh^-_{\pi} \otimes \clh^+_{\pi}$ so that $\pi_{\omega}(\clb_R)'' \equiv \clb(\clh^-_{\pi}) \otimes I_{\clh^+_{\pi}}$
and $\pi_{\omega}(\clb_L)'' \equiv I_{\clh^-_{\pi}} \otimes \clb(\clh^+_{\pi})$. A simple proof goes as follows: $\pi_{\omega}(\clb_R)''$ being a type-I factor, its commutant 
is also a type-I factor. Thus we have an inclusion of type-I sub-factor $\pi_{\omega}(\clb_L)'' \subseteq \pi_{\omega}(\clb_R)'$. Without loss of generality we assume now that
$\pi_{\omega}(\clb_R)'=\clb(\clh^-_{\pi})$ and $\pi_{\omega}(\clb_L)'' \subseteq \clb(\clh^-_{\pi})$. Now we again 
use type-I factor property of $\pi_{\omega}(\clb_L)''$ to write $\clh^-_{\pi} \equiv \clh^-_{\pi}(1) \otimes \clh^-_{\pi}(2)$ 
for two Hilbert spaces $\clh^-_{\pi}(k),\; k=1,2$ so that $\pi_{\omega}(\clb_L)'' \equiv \clb(\clh^-_{\pi}(1)) \otimes I_{\clh^-_{\pi}(2)}$. Since $\omega$ is pure we have 
$\pi_{\omega}(\clb_R)' \bigcap \pi_{\omega}(\clb_L)'=\IC$ and so $\clh^-_{\pi}(2)=\IC$. Thus we arrive at 
our conclusion. The real trouble lies in the fact that the factors $\pi_{\omega}(\clb_R)''$ and $\pi_{\omega}(\clb_L)''$ could be of 
type-III and for a type-III factor we may have non-trivial inclusion with trivial relative commutant. Of course, such a splitting 
relation is not true since tensor product of two type-III factors will give a type-III factor. In [Mo5] we will explain in detail
how Haag duality property finds profound importance in studying {\it reflection symmetry} of a pure translation invariant state and its split property. 

\vsp 
A notion of duality appeared first in the framework of local field theory in Minkowski's space-time formulated by Rudolf Haag [Hag]. We also refer [DHR] for a detailed historical account and 
its subsequent adaptation in conformal field theory. Method that we develop here to prove Haag duality may find some relevance in giving a proof for Haag duality 
property in local field theory as our proof seems to use the underlying group symmetry of the state and simplicity of the $C^*$-algebra $\clb$.       

\vsp 
Before we go further into the results proven in this paper, besides Haag duality property (16), we give a brief history of the present topic and related results. Functional relation (9) is called 
{\it quantum Marvov state } by Luigi Accardi [Ac] as a generalization of classical Markov state or more generally of a classical Gibbs state. He has shown that this functional property holds 
good for unique KMS state for a class of Hamiltonian in a one-dimensional infinite quantum spin lattice with a finite range interaction studied previously by [Ara1]. In [FNW1] M. Fannes, Bruno 
Nachtergaele and R.F. Werner investigated mathematical structures of {\it valence bound states} introduced in [AKLT] and found its close relation with quantum Markov states which they have unified 
in a framework that we have discussed above where $\clm$ is a matrix algebra. When $\clm$ is a matrix algebra in relations (9), it is called {\it finitely correlated state } and the general mathematical 
structures of such a translation invariant states on $\clb$ were further investigated in details in [FNW2] [FNW3], [Ma1],[Ma2],[Ma3] and there afterward in [BJ] and [BJKW] for more general translation invariant states on $\clb$. 
For a brief account on the historical notes about its relevance to more deeper problems in statistical mechanics, we refer interested readers to Bruno Nachtergaele's expository paper [Br]. If $\clk=\IC$, then 
Popescu elements are just some complex numbers i.e. $v_k=\lambda_k$, then $\omega$ is a pure tensor product state. We call such a state a {\it Bernoulli state}. Thus Bernoulli state once restricted to the 
the diagonal algebra $\{S_IS_I^*:|I| < \infty \}$ will give a classical {\it Bernoulli state}. On the other hand Gelu Popescu develops a dilation theory, analogous to that of B. Sz.-Nagy and C. Foiaş 
for a single contraction, for an infinite sequence $v_i$ of non-commuting operators satisfying the condition $\sum v_iv_i^* \le I$. Here we closely follow the presentation of Popescu's dilation as given 
in [BJKW] and Theorem 2.1 in section 2 is a finer version of Popescu's theorem [Po] and {\it commutant lifting theorem } in the present form is a new feature that we explore to an extent in this text, 
particularly while giving proof of Theorem 3.6.     

\vsp 
A finitely correlated pure state $\omega$ gives a type-I factor state once restricted to $\clb_R$. On the other hand Araki and Matsui [AMa] found that the unique ground state of $XY$ model is not 
finitely correlated and in fact once restricted to $\clb_R$ gives a type-III$_1$ factor state. In a recent paper [Mo3] it had been shown that for any translation invariant pure state $\omega$,  $\omega_R$ 
is either a type-I or a type-III factor state. This feature makes classification of translation dynamics an interesting problem which we now describe briefly. Given two translation invariant states $\omega_1$ 
and $\omega_2$ on $\clb$, when can we say their translation dynamics $(\clb,\theta,\omega_1)$ and $(\clb,\theta,\omega_2)$ are isomorphic? i.e. When can we say that there exists an automorphism $\alpha$ on $\clb$ 
such that $\omega_2 \alpha =\omega_1$ and $\theta \alpha = \alpha \theta$ ?

\vsp 
Let $\Theta_k:\pi_k(\clb)'' \raro \pi_k(\clb)''$ be the associated automorphisms where $(\clh,\pi_k,\Omega_k)$ are GNS spaces for $(\clb,\omega_k)
,\;k=1,2$. $(\clb,\theta,\omega_1)$ and $(\clb,\theta,\omega_2)$ are said to be {\it weakly isomorphic} if there exists a unitary 
operator $U:\clh_1 \raro \clh_2$ so that $U\Omega_1=\Omega_2$ and $U\Theta_1(X)U^*=\Theta_2(UXU^*)$ for all $X \in \pi_1(\clb)''$. Kolmogorov-Sinai dynamical entropy [Pa] is an 
invariance for translation dynamics on classical spin chain on lattice $\IZ$ and a celebrated result due to D. Ornstein [Or1] also says that Kolmogorov-Sinai dynamical entropy is
a complete invariance for classical Markov states. However it is not a complete invariance [Or2] for translation invariants states on classical spin chain on $\IZ$. Thus the main obstruction comes 
from the fact that a classical translation invariant state need not be a classical stationary Markov state. However such an obstruction is absent in the full quantum situation as we have shown above that 
any translation invariant state is a stationary quantum Markov state. In [Mo6] we have achieved a partial success in proving that two such states gives weakly isomorphic dynamics if both satisfy Kolmogorov's property 
[AM,Mo2] described below. In that proof the simplicity property of $C^*$-algebra $\clb_R$ played a vital role. This result gives a rare hope for a complete classification for translation dynamics hopefully with some
new symmetries and spacial correlation properties of the states under consideration. 

\vsp 
We now briefly recall the Kolmogorov property [Mo2]. Given a translation invariant state $\omega$ on $\clb$, we set increasing sequence of projections $e_n=[\pi(\theta^n(\clb_R))'\Omega],\;n \in \IZ$ in the GNS space 
$(\clh,\pi,\Omega)$ associated with state $\omega$ on $\clb$. It is simple to check that 
\be 
S_me_nS_m^*=e_{n+m}
\ee 
where $S_m$ is the unitary operator implementing $\theta^m$. The family of operators $(S_n,e_n-|\Omega><\Omega|)$ gives rise to a system of imprimitivity if and only if 
$e_{n} \downarrow |\Omega><\Omega|$ as $n \downarrow -\infty$ [Mac]. We say $\omega$ admits Kolmogorov property if $e_n \downarrow |\Omega><\Omega|$ as $n \downarrow -\infty$. In particular Kolmogorov property
implies purity of $\omega$. But the converse statement is not true in general [Mo6,Appendix]. This makes classification of translation dynamics an interesting mathematical problem. 
As a next step of our goal, we would be aiming to classify translation dynamics with pure states. Such a problem demands a comprehensive understanding about translation invariant pure states on $\clb$. 
In [Mo4] we have shown that a translation invariant pure state $\omega$ on $\clb$ can give only type-I or type-III
factor states once we restrict to $\clb_R$ (Theorem 3.4 in [Mo4] ). One natural question that arises now for two such translation dynamics with pure states. How does restrictions of those states 
to $\clb_R$ determine whether their dynamics are isomorphic or weakly isomorphic? If both give type-I factor states, answer is affirmative for weak isomorphism as type-I property gives Kolmogorov property (Theorem 3.4 in [Mo4]). 
One related important question that also arises here how Kolmogorov property which is little stronger then purity can ensure existence of {\it free energy density } for a translation invariant state? For the definition of free energy 
state and its existence for finitely correlated state, we refer to [HMOP]. We will not address this classification problem here by studying known invariance. Rather we will confine our interest to investigate 
translation invariant pure states with additional symmetry by studying associated quantum Markov state $\psi$ and Markov map $(\clm,\tau,\phi)$. 

\vsp 
Now we explain the basic ingredients in the proof of Haag duality property (16). To that end for the time being we fix a translation invariant factor state $\omega$ on $\clb$ and an extreme point $\psi \in K_{\omega}$. 
We consider the GNS space $(\clh,\pi,\Omega)$ associated with the state $\psi$ on $\clo_d$ and associated Popescu's elements $\clp=(\clk,\clm,v_k,\;1 \le k \le d,\;\Omega)$ arises on support projection 
$P=[\pi(\clo_d)'\Omega]$. Now consider the dual Popescu's elements $\tilde{\clp}=(\clk,\tilde{\clm},\tilde{v}_k;1 \le k \le d,\;\Omega)$ and the completely positive map from $\tilde{\clo}_d$ to $\clb(\clk)$ defined by 
$$\tilde{s}_I\tilde{s}_J^* \raro \tilde{v}_I\tilde{v}^*_J,\;|I|,|J| < \infty$$ 
Let $\pi:\tilde{\clo_d} \raro \clb(\tilde{\clh})$ be the minimal Stinespring representation so that 
$P\pi(\tilde{s}_I\tilde{s}_J^*)P=\tilde{v}_I\tilde{v}_J^*$ for all $|I|,|J| < \infty $. In particular 
we have 
\be 
P\pi(\tilde{s}^*_i)P=\pi(\tilde{s}^*_i)P=\tilde{v}^*_i
\ee
We also consider the state $\tilde{\psi}$ on $\tilde{\clo}_d$ given by
\be 
\tilde{\psi}(\tilde{s}_I\tilde{s}_J^*)=\phi(\tilde{v}_I\tilde{v}^*_J)
\ee
By relation (12) we check that $\tilde{\psi}|\mbox{UHF}_d=\omega|\clb_L$ (see section 3). 
Such relations are perfectly symmetric while moving from $\psi$ to $\tilde{\psi}$ 
except the fact that though $P$ is the support projection of $\psi$ in $\pi(\clo_d)''$, 
it is not guaranteed that $P$ equals to $[\pi(\tilde{\clo}_d)'\Omega]$. We give an explicit example to support this claim in the note 
that follows the proof of Proposition 3.1. 

\vsp 
We consider the GNS space $(H_0,\pi,\Omega)$ associated with $\omega$ on $\clb$. Let $e_0=[\pi_{\omega}(\clb_R)'\Omega]$ and $\tilde{e}_0=[\pi_{\omega}(\clb_L)'\Omega]$ be the 
support projections of $\omega$ in $\pi(\clb_R)''$ and $\pi(\clb_L)''$ respectively. We set projection $q_0=e_0 \tilde{e}_0$ and take $\clk_0$ to be the subspace of $\clh$ 
determined by the projection $q_0$. Also set von-Neumann algebras $\clm^1_0=q_0\pi_{\omega}(\clb_R)''q_0$ and $\tilde{\clm}^1_0=q_0\pi_{\omega}(\clb_L)''q_0$. 
So by our construction we have $\tilde{\clm}^1_0 \subseteq (\clm^1_0)'$. Further the vector state $\phi(x)=<\Omega,x \Omega>$ on $\clb(\clk_0)$ is faithful and normal on 
$\clm_0^1$ and $\tilde{\clm}_0^1$. Further $\omega$ being a factor state, we will also have 
factor property of $\clm_0^1$ and $\tilde{\clm}_0^1$ ( Theorem 2.4 ). What is less obvious is when 
we can expect cyclic property i.e. $[\clm^1_0\Omega]=[\tilde{\clm}^1_0\Omega]=I_{\clk_0}$, identity of $\clk_0$. 
The following theorem answers all non trivial questions that we have raised so far.

\vsp 
\begin{thm} Let $\omega$ be a translation invariant factor state on $\clb$ and $\psi$ be an extremal element in $K_{\omega}$. Then the following statements are equivalent:
 
\vsp 
\NI (a) $[\clm^1_0\Omega]=I_{\clk_0},\;[\tilde{\clm}^1_0\Omega]=I_{\clk_0}$;

\vsp 
\NI (b) $(\clm^1_0)'=\tilde{\clm}^1_0$; 

\vsp         
\NI (c) $\pi_{\omega}(\clb_R)'=\pi_{\omega}(\clb_L)''$;

\vsp 
\NI (d) $[\pi_{\omega}(\clb_R)'\Omega]=[\pi_{\omega}(\clb_L)''\Omega]$;

\vsp 
\NI (e) $\{x \in \clb(\clk): \sum_k \tilde{v}_kx\tilde{v}_k^*= x \} =\clm$; 

\vsp 
\NI (f) $\omega$ is pure.

\vsp 
In such a case $\tilde{\clm}^1_0 = \{x \in \tilde{\clm}: \beta_z(x)=x: z \in H \}$ where 
$H=\{ z \in S^1:\psi = \psi \beta_z \}$.  

\end{thm}
\vsp 
Before we elaborate further on equivalence of above statements we briefly recall results on translation invariant pure state on $\clb =\otimes_{\IZ}M_d(\IC)$ that finds its relevance while proving Haag duality property. There is a one to one affine map between translation invariant states on $\clb$ and translation invariant states on $\clb_R =\otimes_{\IZ_+}M_d(\IC)$ by $\omega \raro \omega_R=\omega_{|\clb_R}$. The inverse map is the inductive limit 
state of $(\clb_R,\psi_R) \raro^{\lambda_n} (\clb_R,\psi_R)$ where $(\lambda_n:n \ge 0)$ is the canonical semi-group of right shifts on $\clb_R$. Pure states on a $\mbox{UHF}$ algebra are studied in the general framework of [Pow]. Such a situation has been investigated also in detail at various degrees of generality in [BJP] and [BJKW] motivated by the development a $C^*$ algebraic method in the study of iterative function systems and its associated wavelet theory. One interesting result in [BJP] says that any translation invariant pure state on $\clb_R$ is also a product state and the canonical endomorphism associated with two such states are unitary equivalent. However such a statement is not true for two translation invariant pure states on $\clb$ as their restriction to $\clb_R$ need not be isomorphic. Theorem 3.4 in [Mo4] says that $\omega_R$ is either a type-I or a type-III factor state on $\clb_R$. Both type of factors are known to exist in the literature of quantum statistical mechanics [BR vol-II,Si,Ma1]. Thus the classification problem of translation dynamics on $\clb$ with invariant pure states on $\clb$ up to unitary isomorphism is a delicate one. In this context one interesting problem that remain open is whether mean entropy [Ru] is an invariance for translation dynamics.   

\vsp 
Since a $\theta$ invariant state $\omega$ on $\clb$ is completely determined by its restriction $\omega_R$ to $\clb_R$, in principle it is possible to describe various properties of $\omega$ including purity by studying their restriction $\omega_R$. 
Theorem 3.2 in [Mo3] gives a precise answer: $\omega$ is pure if and only if there exists a sequence of positive contractive elements $x_n \in \clm$ such that 
$$x_n \raro I,\;\;x_{m+n}\tau_n(x) \raro \phi(x)I$$ 
in strong operator topology for all $x \in \clm$ and $m \ge 1$. As an application of this result, we prove that (e) implies (f). 
This statement can be taken as the correct version of Theorem 7.1 of [BJKW]. Theorem 7.1 in [BJKW] has aimed towards a sufficient condition on Popescu elements $\clp$ for purity of the translation invariant state. However the statement and its proof are faulty as certain argument used in the proof is not time reversal symmetric and a factor state with Popescu elements on support projection satisfies the conditions of the statement of Theorem 7.1. One natural remedy to add additional hypothesis that (e) holds. In particular Lemma 7.6 in [BJKW] needs that additional assumption related to the support projection of the dual state $\tilde{\psi} \in K_{\tilde{\omega}}$. Besides this additional structure proof of Lemma 7.8 in [BJKW] is also not complete unless we find a proof for $\tilde{\clm}=\clm'$ ( we retained same notations here in the text) for such a factor state $\omega$. Such a problem could have been solved if there were any method which shows directly that Takesaki's conditional [Ta] expectation exists from $\clm'$ onto $\tilde{\clm}$. For von-Neumann algebras $\cln \subseteq \clm$, a normal unital completely positive map $E_c:\clm \raro \cln$ is called 
{\it normal conditional expectation } if 
\be 
E_c(yxz)=yE_c(x)z
\ee
for all $y,z \in \cln$ and $x \in \clm$. 
A theorem of Takesaki's [Ta, also see AcC] says that conditional expectation $E_c$ preserving a faithful normal state $\phi$ on $\clm$ exists if and only if the modular group $\sigma=(\sigma_t:t \in \IR)$ 
associated with $\phi$, which preserves $\clm$, also preserves $\cln$, i.e. $\sigma_t(x) \in \cln$ for all $t \in \IR$ and $x \in \cln$. Thus main body of the proof for Theorem 7.1 even with the 
additional natural hypothesis $p_0=q_0$, in [BJKW] is incomplete. 

\vsp 
Besides the proof of (e) implies (f), other implication follows along with while we will prove the hardest part of this theorem namely (f) implies (c) i.e. purity implies Haag duality property (16). For the proof 
we have explored the set of representation of $\clb$ quasi-equivalent to $\pi_{\omega}$ and equip it with a strict partial ordering depending on our situation to prove Haag duality. Mackey's system of imprimitivity [Mac] 
plays a crucial role even though a pure state not necessarily give rise to a Mackey's system of imprimitivity generated by the support projection $e_0$ with respect to shift. Though we have worked here with amalgamated representation 
of $\tilde{\clo}_d \otimes \clo_d$ in [BJKW], it seems that just for Haag duality one can avoid doing so. It seems that the underlining group $\IZ$ can easily be replaced by $\IZ^k$ for some $k \ge 2$ and 
wedge duality for a pointed cone can be proved by following the same ideas. We defer this line of analysis leaving it for work as its relation with problems in quantum spin 
chain in higher dimensional lattice needs some additional structure. We also defer application of Haag duality property in studying symmetry and correlation of a translation invariant pure state to another
paper [Mo5].  

\vsp
The paper is organized as follows. In section 2 we study Popescu's dilation associated with a translation invariant state on Cuntz algebra $\clo_d$ and review `commutant lifting theorem' investigated in [BJKW]. 
The proof presented here remove the murky part of the proof of Theorem 5.1 in [BJKW]. In section 3 we explore both the notion of Kolmogorov's shift and its intimate relation with Mackey's imprimitivity system 
to explore a duality argument introduced in [BJKW]. We find a useful necessary and sufficient condition (Theorem 1.1 (a) ) in terms of support projection of Cuntz's state for a translation invariant factor state 
$\omega$ on $\clb$ to be pure. The criterion on support projection is crucial for our main mathematical result Theorem 3.6. 

\begin{rem} 
The paper ``On Haag Duality for Pure States of Quantum Spin Chain'' 
by authors: M. Keyl, Taku Matsui, D. Schlingemann, R. F. Werner, Rev. Math. Phys. 20:707-724,2008 has an incomplete proof for Haag duality property as Lemma 4.3 in that paper has a faulty argument.   
\end{rem}

\section{ States on $\clo_d$ and the commutant lifting theorem }

\vsp
In this section we essentially recall results from [BJKW] and organize it with additional remarks and arguments as 
it needed to understand the present problem investigated in section 3. In the following we recall a commutant lifting 
theorem ( Theorem 5.1 in [BJKW] ), crucial for our purpose.

\begin{thm}
Let $v_1,v_2,...,v_d$ be a family of bounded operators on a Hilbert space
$\clk$ so that $\sum_{1 \le k \le d} v_kv_k^*=I$. Then there exists a unique up to
isomorphism Hilbert space $\clh$, a projection $P$ on $\clk$ and a family of isometries
$\{S_k:,\;1 \le k \le d \}$ satisfying Cuntz's relation so that
\be
PS^{*}_{k}P=S_k^*P=v^*_k
\ee
for all $1 \le k \le d$ and $\clk$ is cyclic for the representation i.e. the vectors
$\{ S_I\clk: |I| < \infty \}$ are total in $\clh$.

Moreover the following holds:

\NI (a) $\Lambda^n(P) \uparrow I$ as $n \uparrow \infty$ where $\Lambda(X)=\sum_k S_kXS_k^*$;

\NI (b) For any $D \in \clb_{\tau}(\clk)=\{x \in \clb(\clk): \tau(x)=\sum_{1 \le k \le d} v_kxv_k^*=x \}$, 
$\Lambda^n(D) \raro X'$ weakly as $n \raro \infty$ for some $X'$ in the commutant $\{S_k,S^*_k: 1 \le k \le d \}'$ so that $PX'P=D$. Moreover
the self adjoint elements in the commutant $\{S_k,S^*_k: 1 \le k \le d \}'$ is isometrically
order isomorphic with the self adjoint elements in $\clb_{\tau}(\clk)$ via the
surjective map $X' \raro PX'P$, where $\clb_{\tau}(\clk)=\{x \in \clb(\clk): \sum_{1 \le k \le d }
v_kxv^*_k=x \}.$

\NI (c) $\{v_k,v^*_k,\;1 \le k \le d \}' \subseteq \clb_{\tau}(\clk)$ and equality holds if and only if
$P \in \{S_k,S_k,\;1 \le k \le d \}''$.

\vsp 
If $(w_i)$ be another such an Popescu elements on a Hilbert space $\clk'$ such that there exists an operator 
$u:\clk \raro \clk'$ so that $\sum_k w_kuv_k^*=u$ then there exists an operator $U:\clh_v \raro \clh_w$ so 
that $\pi'(x)U=U\pi(x)$ where $(\clh_w,\pi',S'_i)$ are Popescu dilation of $(w_i)$ and $\pi'$ is the associated
minimal representation of $\clo_d$. In particular $U$ is isometry, unitary if $u$ is so respectively. If $u$ is unitary 
and $\clk=\clk'$ then we can as well take $\clh_v=\clh_w$.  

\end{thm}

\vsp
\begin{proof}
Following Popescu [Po] we define a completely positive map
$R: \clo_d \raro \clb(\clk)$ by
\be
R(s_Is^*_J)=v_Iv^*_J
\ee
for all $|I|,|J| < \infty$. The representation $S_1,..,S_d$ of $\clo_d$ on $\clh$ thus may be taken to be the Stinespring
dilation of $R$ [BR, vol-2] and uniqueness up to unitary equivalence follows from uniqueness of the Stinespring
representation. That $\clk$ is cyclic for the representation follows from the minimality property of the Stinespring
dilation. 

\vsp 
For (a) let $Q$ be the limiting projection of $\Lambda^n(P)$ as $n \uparrow \infty$. Then we have $\Lambda(Q)=Q$ i.e. $Q\Lambda(I-Q)Q=0$ 
and so $(I-Q)S_k^*Q=0$. Interchanging the role of $Q$ with $I-Q$, we get $QS_k^*(I-Q)=0$. This shows $QS_k^*=S_k^*Q$ for all $1 le k \le d$. 
Since $Q$ is a projection, taking adjoint in the relation, we get $Q \in \{S_k,S^*_k\}'$. That $Q \ge P$ is obvious since $\Lambda^n(P) \ge P$ 
for all $n \ge 1$. In particular $QS_If=S_If$ for all $f \in \clk$ and $|I| < \infty$. Hence $Q=I$ by the
cyclicity of $\clk$. 

\vsp 
For (b) essentially we deffer from the argument used in Theorem 5.1 in [BJKW]. We fix any $D
\in \clb_{\tau}(\clk)$ and note that $P\Lambda_k(D)P=\tau_k(D)=D$ for any $k \ge 1$. Thus for any integers $n > m$
we have
$$\Lambda_m(P)\Lambda_n(D)\Lambda_m(P)=\Lambda_m(P\Lambda_{n-m}(D)P)=\Lambda_m(D)$$
Hence for any fixed $m \ge 1$ limit $<f,\Lambda_n(D)g>$ as $n \raro \infty$ exists for all
$f,g \in \Lambda_m(P)$. Since the family of operators $\Lambda_n(D)$ is uniformly bounded and $\Lambda_m(P) \uparrow I$
as $m \raro \infty$, a standard density argument guarantees that the weak operator limit of $\Lambda_n(D)$ exists as
$n \raro \infty$. Let $X'$ be the limit. So $\Lambda(X')=X'$, by Cuntz's relation, $X' \in \{S_k,S^*_k:1 \le k \le k \}'$.
Since $P\Lambda_n(D)P=D$ for all $n \ge 1$, we also conclude that $PX'P=D$ by taking limit $n \raro \infty$. Conversely,
it is obvious that $P\{ S_k,S^*_k:\; k \ge 1 \}'P \subseteq \clb_{\tau}(\clk)$. Hence
we can identify $P\{S_k,S^*_k:\; k \ge 1 \}'P$ with $\clb_{\tau}(\clk)$.

\vsp
Further it is obvious that $X'$ is self-adjoint if and only if $D=PX'P$ is self-adjoint ( since $\Lambda^n(P) \uparrow I$ as $n \uparrow \infty$, 
by taking limit $\Lambda^n(PXP)=\Lambda^n(PX^*P)$ shows that $X^*=X$ if $D^*=D$ ). Now fix any self-adjoint element $D \in \clb_{\tau}(\clk)$. Since 
the identity operator on $\clk$ is an element in $\clb_{\tau}(\clk)$ for any $\alpha \ge 0$ for which $- \alpha P  \le D \le \alpha P$, we have $\alpha 
\Lambda_n(P) \le \Lambda_n(D) \le \alpha \Lambda_n(P)$ for all $n \ge 1$. By taking limit $n \raro \infty$ we
conclude that $- \alpha I \le X' \le \alpha I$, where $PX'P=D$. Since operator norm of a self-adjoint element $A$
in a Hilbert space is given by
$$||A||= \mbox{inf}_{\alpha \ge 0}\{\alpha: - \alpha I  \le A \le \alpha I \}$$
we conclude that $||X'|| \le ||D||$. That $||D||=||PX'P|| \le ||X'||$ is obvious, $P$ being a
projection. Thus the map is isometrically order isomorphic taking self-adjoint elements of the commutant
to self-adjoint elements of $\clb_{\tau}(\clk)$.

\vsp
We are left to prove (c). Inclusion is trivial. For the last part, we assume first that $P \in \pi(\clo_d)''$. For any element $D$
in $\clb_{\tau}(\clk)$ there exists an element $X'$ in $\{S_k,S^*_k,\;1 \le k \le d\}'$ so that $PX'P=D$. Since $P$ commutes with $X'$, 
we verify that 
$$Dv^*_k= PX'PS^*_kP=PX'S^*_kP=PS^*_kX'P=PS^*_kPX'P=v^*_kD$$ 
We also have $D^* \in \clb_{\tau}(\clk)$ and thus $D^*v^*_k=v^*_kD^*$. Hence $D \in \{v_k,v^*_k: 1 \le k \le d \}'$. Since $P\pi_{\hat{\omega}}
(\clo_d)'P = \clb(\clk)_{\tau}$, we conclude that $\clb(\clk)_{\tau} \subseteq \clm'$. Thus equality holds whenever
$P \in \{S_k,S_k^*,\;1 \le k \le d \}''$. 

\vsp 
For the converse note that by commutant lifting property self-adjoint elements
of the commutant $\{S_k,S^*_k,1 \le k \le d \}'$ is order isometric with the algebra $\clm'$ via the map $X' \raro PX'P$.
Hence $P \in \{S_k,S^*_k,1 \le k \le d \}''$ by Proposition 4.2 in [BJKW]. 

\vsp 
For the proof of intertwining relation and their property we refer to main body of the proof of 
Theorem 5.1 in [BJKW]. 
\end{proof}

\vsp
A family $(v_k,1 \le k \le d)$ of contractive operators on a Hilbert space $\clk$ is called Popescu's elements if $\sum_kv_kv_k^*=I$ 
and the dilation $(\clh,P,\clk,S_k,1 \le k \le d)$ in Theorem 2.1 is called Popescu's dilation to Cuntz elements. In the 
following proposition we deal with a family of minimal Popescu elements for a state on $\clo_d$. 

\begin{pro}
There exists a canonical one-one correspondence between the following objects:

\vsp
\NI (a) States $\psi$ on $\clo_d$

\vsp
\NI (b) Function $C: \cli \times \cli \raro \IC$ with the following properties:

\NI (i) $C(\emptyset, \emptyset)=1$;

\NI (ii) for any function $\lambda:\cli \raro \IC$ with finite support we have
         $$\sum_{I,J \in \cli} \overline{\lambda(I)}C(I,J)\lambda(J) \ge 0$$

\NI (iii) $\sum_{i \in \IZ_d} C(Ii,Ji) = C(I,J)$ for all $I,J \in \cli$.

\vsp
\NI (c) Unitary equivalence class of objects $(\clk,\Omega,v_1,..,v_d)$ where

\NI (i) $\clk$ is a Hilbert space and $\Omega$ is a unit vector in $\clk$;

\NI (ii) $v_1,..,v_d \in \clb(\clk)$ so that $\sum_{i \in \IZ_d} v_iv^*_i=1$;

\NI (iii) the linear span of the vectors of the form $v^*_I\Omega$, where $I \in \cli$, is dense in $\clk$.

\vsp
The correspondence is given by a unique completely positive map $R: \clo_d \raro \clb(\clk)$, where 

\NI (i) $R(s_Is^*_J)=v_Iv^*_J;$

\NI (ii) $\psi(x)=<\Omega,R(x)\Omega>;$

\NI (iii) $\psi(s_Is^*_J)=C(I,J)=<v^*_I \Omega,v^*_J\Omega>;$

\NI (iv) For any fixed $g \in U_d$, the completely positive map $R_g:\clo_d \raro \clb(\clk)$ defined by $R_g= R
\circ \beta_g$ give rises to a Popescu system given by $(\clk,\Omega,\beta_g(v_i),..,\beta_g(v_d))$ where
$\beta_g(v_i)=\sum_{1 \le j \le d} \overline{g^i_j} v_j.$

\end{pro}

\vsp
\begin{proof} 
For a proof we simply refer to Proposition 2.1 in [BJKW]. 
\end{proof} 

\vsp
The following is a simple consequence of Theorem 2.1 valid for a $\lambda$-invariant state $\psi$ on $\clo_d$. 
This proposition will have very little application in the main body of this paper but this gives a clear picture explaining the
delicacy of the present problems. 

\begin{pro}

Let $\psi$ be a state on $\clo_d$ and $(\clh,\pi,\Omega)$ be the GNS space associated with 
$(\clo_d,\psi)$. We set $S_i=\pi(s_i)$ and normal state $\psi_{\Omega}$ on $\pi(\clo_d)''$ defined by 
$$\psi_{\Omega}(X) = <\Omega, X \Omega>$$
Let $P$ be the projection on the closed subspace $\clk$ generated by the vectors $\{S^*_I\Omega: |I| < \infty \}$ and 
\be
v_k=PS_kP
\ee
for $1 \le k \le d$. Then following holds:

\vsp
\NI (a) $\{v^*_I\Omega: |I| < \infty \}$ is total in $\clk$.

\NI (b) $\sum_{1 \le k \le d} v_kv_k^*=I$;

\NI (c) $S_k^*P=PS_k^*P$ for all $1 \le k \le d$; 

\NI (d) For any $I=(i_1,i_2,...,i_k),J=(j_1,j_2,...,j_l)$ with $|I|,|J| < \infty$ we have 
\be 
\psi(s_Is^*_J)= <\Omega,v_Iv^*_J\Omega>
\ee 
and the vectors $\{ S_If: f \in \clk,\;|I| < \infty \}$ are total in
the GNS Hilbert space associated with $(\clo_d,\psi)$. Further such a family $(\clk,\;v_k,\;1 \le k \le d,\;\omega)$ satisfying
(a) to (d) is determined uniquely up to isomorphism. 

\vsp
Conversely given a Popescu system $(\clk,v_k,\;1 \le k \le d,\Omega)$ satisfying (a) and (b) there exists a unique   state 
$\psi$ on $\clo_d$ so that (c) and (d) are satisfied.

Furthermore the following statements are valid:
\vsp
\NI (e)  If the normal state $\phi(x)=<\Omega,x \Omega>$ on the von-Neumann algebra
$\clm=\{ v_i, v^*_i \}''$ is invariant for the Markov map $\tau(x)=\sum_{1 \le k
\le d} v_ixv^*_i,\;x \in \clm$ then $\psi$ is $\lambda$ invariant and $\phi$ is faithful 
on $\clm$.

\vsp
\NI (f) If $P \in \pi(\clo_d)''$ then following are equivalent:

\NI (i) $\psi$ is an ergodic state for $(\clo_d,\lambda)$; 

\NI (ii) $(\clm,\tau,\phi)$ is ergodic;

\NI (ii) $\clm$ is a factor. 
\end{pro} 

\vsp
\begin{proof} 
We fix a state $\psi$ and consider the GNS space  $(\clh,\pi,\Omega)$
associated with $(\clo_d,\psi)$ and set $S_i=\pi(s_i)$. It is obvious 
that $S^*_kP \subseteq P$ for all $1 \le k \le d $, thus $P$ is the minimal subspace 
containing $\Omega$ and invariant by all $\{ S^*_k;\;1 \le k \le d \}$ i.e.
\be
PS^*_kP=S^*_kP
\ee
Thus $v^*_k=PS^*_kP=S_k^*P$ and so $\sum_k v_kv^*_k = \sum_k PS_kS^*_kP= P$ which is 
identity operator in $\clk$. This completes the proof of (a) (b) and (c).

\vsp
For (d) we note that
$$\psi(s_Is^*_J)$$
$$=<\Omega, S_IS^*_J\Omega>$$
$$=<\Omega,PS_IS^*_JP\Omega>$$
$$=<\Omega,v_Iv^*_J\Omega>.$$
Since $\clh$ is spanned by the vectors $\{S_IS^*_J\Omega:|I|,|J| < \infty \}$ and
$\clk$ is spanned by the vectors $\{ S^*_J\Omega=v^*_J\Omega:|I| < \infty \}$, 
$\clk$ is cyclic for $S_I$ i.e. the vectors $\{ S_I\clk: |I| < \infty \}$ spans 
$\clh$. Uniqueness up to isomorphism follows as usual by total property of vectors 
$v^*_I\Omega$ in $\clk$. 

\vsp
Conversely for a Popescu's elements $(\clk,v_i,\Omega)$ satisfying (a) and (b), we consider
the family $(\clh,S_k,\;1 \le k \le d,P)$ of Cuntz's elements defined as in
Theorem 2.1. We claim that $\Omega$ is a cyclic vector for the representation $\pi(s_i) \raro
S_i$. Note that by our construction vectors $\{S_If,f \in \clk: |I| < \infty \}$ are total
in $\clh$ and $v^*_J\Omega=S^*_J\Omega$ for all $|J| < \infty$. Thus by our hypothesis that
vectors $\{v^*_J\Omega:|I| < \infty \}$ are total in $\clk$, we verify that
vectors $\{S_IS^*_J\Omega: |I|,|J| < \infty \}$ are total in $\clh$. Hence $\Omega$ is a cyclic for
the representation $s_i \raro S_i$ of $\clo_d$.

\vsp
We are left to prove (e) and (f). It simple to note by (d) that $\psi \lambda=\psi$ i.e.
$$\sum_i <\Omega,S_iS_IS^*_JS^*_i\Omega>= \sum_i <\Omega,v_iv_Iv^*_Jv^*_i\Omega> $$
$$=<\Omega,v_Iv_J^*\Omega> = <\Omega,S_IS^*_J\Omega>$$ for all $|I|,|J| < \infty$ where 
in the second equality we have used our hypothesis that the vector state $\phi$ on $\clm$ 
is $\tau$-invariant. In such case we aim now to show that $\phi$ is faithful on $\clm$. To that end 
let $p'$ be the support projection in $\clm$ for $\tau$ invariant state $\phi$. Thus $\phi(1-p') = 0$ 
i.e. $p'\Omega=\Omega$ and by invariance we also have $\phi(p'\tau(1-p')p')=\phi(1-p')=0$. Since $p'\tau(1-p')p' \ge 0 $ and an element in 
$\clm$, by minimality of support projection, we 
conclude that $p'\tau(1-p')p'=0$. Hence $p'\Omega=\Omega$ and $p'v^*_kp'=v^*_kp'$ for all $1 \le k \le d$. 
Thus $p'v^*_I\Omega=v^*_I\Omega$ for all 
$|I| < \infty$. As $\clk$ is the closed linear span of the vectors $\{ v_I^*\Omega: |I| < \infty \}$, we conclude 
that $p'=p$. In other words $\phi$ is faithful on $\clm$. This completes the proof for (e).

\vsp
We are left to show (f) where we assume that $P \in \pi(\clo_d)''$. $\Omega$ being a cyclic vector for $\pi(\clo_d)''$, the weak$^*$ limit of the increasing 
projection $\Lambda^k(P)$ is $I$. Thus by Theorem 3.6 in [Mo1] we have $(\pi(\clo_d)'',\Lambda,\psi_{\Omega})$ is ergodic if and only if the reduced dynamics 
$(\clm,\tau,\phi)$ is ergodic. For the last part of the statement, we need to show for a projection $e \in \clm$, $\tau(e)=e$ if and only if $e \in \clm \bigcap \clm'$.
It is an easy consequence since $\tau(e)=e$ says that $e\tau(I-e)e=0$ and so $(1-e)v_k^*e=0$. Changing the role of $e$ by $I-e$, we also get $ev_k^*(I-e)=0$ for
all $1 \le k \le d$. Thus we get that $e$ commutes with each $v_k$. $P$ being in $\pi(\clo_d)''$, $v_k \in \clm$. So $e \in \clm'$. Thus by a standard theorem [La,Ev,Fr,BJKW] 
ergodic property is equivalent to factor property of $\clm$. 
\end{proof} 

\vsp
Before we move to our next result, we comment here that in general for a $\lambda$ invariant state on $\clo_d$, 
the normal state $\phi$ on $\clm=\{v_k,v_k^*:1 \le k \le d \}''$ need not be invariant for $\tau$. To that 
end we consider ( [BR] vol-II page 110 ) the unique KMS state $\psi=\psi_{\beta}$ for the automorphism $\alpha_t(s_i)=e^{it}s_i$ on $\clo_d$. $\psi$ is $\lambda$ 
invariant and $\psi_{|\mbox{UHF}_d}$ is the unique faithful trace. $\psi$ being a KMS state for an automorphism, the normal state induced by the cyclic vector on 
$\pi_{\psi}(\clo_d)''$ is also separating for $\pi(\clo_d)''$. As $\psi \beta_z = \psi$ for all $z \in S^1$ we have 
$<\Omega,\pi(s_I)\Omega>=<\Omega,\beta_z(s_I)\Omega> = z^{|I|}<\Omega,\pi(s_I)\Omega>$ for all $z \in S^1$ and 
so $<\Omega,\pi(s_I)\Omega>=0$ for all $|I| \ge 1$. In particular $<\Omega,v_I^*\Omega>=0$ where $(v_i)$ are defined 
as in Proposition 2.3 and thus $<v_i\Omega,v^*_I\Omega>=<\Omega,v^*_iv_I\Omega>=0$ for all $1 \le i \le d$. Hence $v_i\Omega=0$. By Proposition 2.3 (e), $\Omega$ is 
separating for $\clm$ and so we get $v_i=0$ for all $1 \le i \le d$ and this contradicts that $\sum_i v_iv_i^*=1$. Thus we conclude by Proposition 2.3 (e) that 
$\phi$ is not $\tau$ invariant on $\clm$. This example also indicates that the support projection of a $\lambda$ invariant state $\psi$ in $\pi(\clo_d)''$ need 
not be equal to the minimal sub-harmonic projection $P$ i.e. the closed span of vectors $\{ S_I^*\Omega: |I| < \infty \}$ containing $\Omega$ and $\{v_Iv_J^*: |I|,|J| < \infty \}$ 
need not be even an algebra.        

\vsp
Now we aim to deal with another class of Popescu elements associated with an $\lambda$-invariant state on $\clo_d$. 
In fact this class of Popescu elements will play a significant role for the rest of the text and we will repeatedly use
this proposition in section 3. 

\vsp
\begin{pro} 

Let $(\clh,\pi,\Omega)$ be the GNS representation of a $\lambda$ invariant state 
$\psi$ on $\clo_d$ and $P$ be the support projection of the normal state $\psi_{\Omega}(X)=<\Omega,X\Omega>$ in the 
von-Neumann algebra $\pi(\clo_d)''$. Then the following holds:

\vsp
\NI (a) $P$ is a sub-harmonic projection for the endomorphism $\Lambda(X)=\sum_k S_kXS^*_k$ on $\pi(\clo_d)''$
i.e. $\Lambda(P) \ge P$ satisfying the following:

\vsp 
\NI (i) $\Lambda_n(P) \uparrow I$ as $n \uparrow \infty$;

\vsp 
\NI (ii) $PS^*_kP=S^*_kP,\;\;1 \le k \le d$;

\vsp 
\NI (iii) $\sum_{1 \le k \le d} v_kv_k^*=I$ 

where $S_k=\pi(s_k)$ and $v_k=PS_kP$ for $1 \le k \le d$;

\vsp
\NI (b) For any $I=(i_1,i_2,...,i_k),J=(j_1,j_2,...,j_l)$ with $|I|,|J| < \infty$ we have $\psi(s_Is^*_J) =
<\Omega,v_Iv^*_J\Omega>$ and the vectors $\{ S_If: f \in \clk,\;|I| < \infty \}$ are total in $\clh$;

\vsp
\NI (c) The von-Neumann algebra $\clm=P\pi(\clo_d)''P$, acting on the Hilbert space
$\clk$ i.e. the range of $P$, is generated by $\{v_k,v^*_k:1 \le k \le d \}''$ and the normal state
$\phi(x)=<\Omega,x \Omega>$ is faithful on the von-Neumann algebra $\clm$.

\vsp 
\NI (d) The self-adjoint part of the commutant of $\pi(\clo_d)'$ is norm and order isomorphic to the space
of self-adjoint fixed points of the completely positive map $\tau$. The isomorphism takes $X' \in \pi(\clo_d)'$
onto $PX'P \in \clb_{\tau}(\clk)$, where $\clb_{\tau}(\clk) = \{ x \in \clb(\clk): \tau(x)=x \}$ and 
$\tau(x)=\sum_kv_kxv^*_k,\;x \in \clb(\clk)$. Furthermore $\clm' = \clb_{\tau}(\clk)$.

\vsp
Conversely let $\clm$ be a von-Neumann algebra generated by a family $\{v_k: 1 \le k \le d \}$ of bounded operators
on a Hilbert space $\clk$ so that $\sum_kv_kv_k^*=1$ and the commutant 
\be 
\clm'=\{x \in \clb(\clk): \sum_k v_kxv_k^*=x \}
\ee 
Then the Popescu dilation $(\clh,P,S_k,\;1 \le k \le d)$ described in Theorem 2.1 satisfies the following:

\vsp 
\NI (i) $P \in \{S_k,S^*_k,\;1 \le k \le d \}''$;

\vsp 
\NI (ii) For any faithful normal $\tau-$invariant state $\phi$ on $\clm$ there exists a state $\psi$ on $\clo_d$ defined by
$$\psi(s_Is^*_J)=\phi(v_Iv^*_J),\;|I|,|J| < \infty $$
such that the GNS space associated with $(\clm,\phi)$ is the support projection
for $\psi$ in $\pi(\clo_d)''$ satisfying (a)-(d). 

\vsp
Further for a given $\lambda$-invariant state $\psi$, the family $(\clk,\clm,v_k\;1 \le k \le d,\phi)$ satisfying (a)-(d) 
is determined uniquely up to unitary conjugation. 

\vsp
\NI (e) $\phi$ is a faithful normal $\tau$-invariant state on $\clm$. Furthermore the following statements are equivalent:

\NI (i) $(\clo_d,\lambda,\psi)$ is ergodic;

\NI (ii) $(\clm,\tau,\phi)$ is ergodic;

\NI (iii) $\clm$ is a factor.

\end{pro} 

\vsp
\begin{proof} 
$\Lambda(P)$ is also a projection in $\pi_{\psi}(\clo_d)''$ so that
$\psi_{\Omega}(\Lambda(P))=1$ by invariance property. Thus we have $\Lambda(P) \ge P$ i.e.
$P\Lambda(I-P)P=0$. Hence we have
\be
PS^*_kP=S^*_kP
\ee
Moreover by $\lambda$ invariance property we also note that the faithful normal state
$\phi(x)=<\Omega,x\Omega>$ on the von-Neumann algebra $\clm=P\pi_{\psi}(\clo_d)''P$
is invariant for the reduced Markov map [Mo1] on $\clm$ given by
\be
\tau(x)=P\Lambda(PxP)P
\ee
\vsp
We claim that $\mbox{lim}_{n \uparrow \infty}\Lambda^n(P)= I$. That $\{ \Lambda^n(P): n \ge 1 \}$
is a sequence of increasing projections follows from sub-harmonic property of $P$ and endomorphism
property of $\Lambda$. Let the limiting projection be $Y$. Then $\Lambda(Y)=Y$ and so $Y \in
\{S_k,S^*_k \}'$. Since by our construction, the GNS Hilbert space $\clh_{\pi_{\hat{\omega}}}$ is
generated by $S_IS^*_J\Omega$, $Y$ is a scalar. $Y$ being a non-zero projection, it is
the identity operator in $\clh_{\pi_{\psi}}$.

\vsp
Now it is routine to verify (a) (b) and (c). For the first part of (d) we appeal to Theorem 2.1. For the last part
note that for any invariant element $D$ in $\clb(\clk)$ there exists an element $X'$ in $\pi(\clo_d)'$ so that $PX'P=D$.
Since $P \in \pi(\clo_d)''$ we note that $(1-P)X'P=0$.
Now since $X' \in \{ S_k,S^*_k \}'$, we verify that $Dv^*_k= PXPS^*_kP=PXS^*_kP=PS^*_kXP=PS^*_kPXP=v^*_kD$. Since $D^* \in
\clb_{\tau}(\clk)$ we also have $D^*v^*_k=v^*_kD^*$. Thus $D \in \{v_k,v^*_k: 1 \le k \le d \}'=\clm'$. Since
$P\pi_{\hat{\omega}}(\clo_d)'P=
\clb(\clk)_{\tau}$, we conclude that $\clb(\clk)_{\tau} \subseteq \clm'$. The reverse inclusion is trivial.
This completes the proof for (d).

\vsp
For the converse part: for (i), since by our assumption and the commutant lifting property self-adjoint elements of the commutant
$\{S_k,S^*_k,1 \le k \le d \}'$ is order isometric with the algebra $\clm'$ via the map $X' \raro PX'P$,
$P \in \{S_k,S^*_k,1 \le k \le d \}''$ by Proposition 4.2 in [BJKW]. For (ii) without loss of generality, we assume
that $\phi(x)=<\Omega,x\Omega>$ for all $x \in \clm$ and $\Omega$ is a cyclic and separating vector for $\clm$.
( otherwise we set state $\psi(s_Is^*_J)=\phi(v_Iv_J^*)$ and consider its GNS representation )
We are left to show that $\Omega$ is a cyclic vector for the representation $\pi(s_i) \raro S_i$. To that end 
let $Y \in \pi(\clo_d)'$ be the projection on the subspace generated by the vectors $\{S_IS^*_J\Omega:  
|I|,|J| < \infty \}$. Note that, $P$ being an element in $\pi(\clo_d)''$, $Y$ also commutes with all the elements of 
$P\pi(\clo_d)''P=P \clm P$. Hence $Yx\Omega=x\Omega$ for all $x \in \clm$. Thus $Y \ge P$.  
Since $\Lambda_n(P) \uparrow I$ as $n \uparrow \infty$ by our construction, we conclude that $Y=\Lambda_n(Y) 
\ge \Lambda_n(P) \uparrow I$ as $n \uparrow \infty$. Hence $Y=I$. In other words
$\Omega$ is cyclic for the representation $s_i \raro S_i$. This completes the proof for (ii). 

\vsp
Uniqueness up to unitary isomorphism follows as GNS representation is determined uniquely up to unitary conjugation 
and so its support projection.  

\vsp
The first part of (e) we note that $PS_IS^*_JP=v_Iv_J^*$ for all $|I|,|J| < \infty$ and thus $\clm=P\pi(\clo_d)''P$ is the von-Neumann 
algebra generated by $\{v_k,v_k^*:1 \le k \le d \}$ and thus $\tau(x)=P\Lambda(PxP)P$ for all $x \in \clm$. That $\phi$ is 
$\tau(x)=\sum_k v_kxv_k^*$ invariant follows as $\psi$ is $\lambda$-invariant. We are left to prove equivalence of statements (i)-(iii).  

\vsp
By Theorem 3.6 in [Mo1] the Markov semi-group $(\clm,\tau,\phi)$ is ergodic if and
only if $(\pi(\clo_d)'',\Lambda,\psi_{\Omega})$ is ergodic ( here we need to recall by (a) that
$\Lambda_n(P) \uparrow I$ as $n \uparrow \infty$ ). By a standard result [Ev, also BJKW] $(\clm,\tau,\phi)$
is ergodic if and only if there is no non trivial projection $e$ invariant for $\tau$ i.e. $\cli^{\tau} = \{e \in \clm: e^*=e,e^2=e,\tau(e)=e \}=\{0,1 \}$. 
If $\tau(e)=e$ for some projection $e \in \clm $ then $(1-e)\tau(e)(1-e)=0$ and so $ev^*_k(1-e)=0$. Same is true 
if we replace $e$ by $1-e$ as $\tau(1-e)=\tau(1)-\tau(e)=1-e$ and so $(1-e)v^*_ke=0$. Thus $e$ commutes with $v_k,v_k^*$ for all $1 \le k \le d$. 
Hence $\cli^{\tau} \subseteq \clm \bigcap \clm'$. Inequality in the reverse direction is trivial and thus $\cli^{\tau}$ is trivial if and only if $\clm$ is a factor. 
This shows equivalence of (ii) and (iii) follows by a standard result [La,Fr] in non-commutative ergodic theory. This completes the proof.
\end{proof} 

\vsp
The following two propositions are essentially easy adaptations of results proved in [BJKW, Section 6 and Section 7]. These results are  
crucial in our present framework. 

\begin{pro} 
Let $\psi$ be a $\lambda$ invariant factor state on $\clo_d$ and $(\clh,\pi,\Omega)$
be its GNS representation. Then the following holds:

\NI (a) The closed subgroup $H=\{z \in S^1: \psi \beta_z =\psi \}$ is equal to 

$$\{z \in S^1: \beta_z \mbox{extends to an automorphism of } \pi(\clo_d)'' \} $$ 

\NI (b) Let $\clo_d^{H}$ be the fixed point sub-algebra in $\clo_d$ under the action of gauge group $\{ \beta_z: z \in H \}$. Then  
$\pi(\clo_d^{H})'' = \pi(\mbox{UHF}_d)''$.

\NI (c) If $H$ is a finite cyclic group of $k$ many elements and $\pi(\mbox{UHF}_d)''$ is a factor, 
then $\pi(\clo_d)'' \bigcap \pi(\mbox{UHF}_d)' \equiv \IC^m$ where $1 \le m \le k$.  
\end{pro}

\vsp 
\begin{proof}  
It is simple that $H$ is a closed subgroup. For any fix $z \in H $ we define unitary operator $U_z$ 
extending the map $\pi(x)\Omega \raro \pi(\beta_z(x))\Omega$ and check that the map $X \raro U_zXU^*_z$ extends $\beta_z$ 
to an automorphism of $\pi(\clo_d)''$. For the converse we will use the hypothesis that $\psi$ is a $\lambda$-invariant 
factor state and $\beta_z \lambda= \lambda \beta_z$ to guarantee that $\psi \beta_z (X) = {1 \over n} 
\sum_{1 \le k \le n} \psi \lambda^k \beta_z(X) = {1 \over n} \sum_{1 \le k \le n} \psi \beta_z \lambda^k(X) 
\raro \psi(X)$ as $n \raro \infty$ for any $X \in \pi(\clo_d)''$, where we have used the same symbol $\beta_z$ 
for the extension. Hence $z \in H$. 

\vsp
In the following instead of working with $\clo_d$ we should be working with the inductive limit $C^*$ algebra and their inductive limit states. 
For simplicity of notation we still use $UHF_d,\clo_d$ for its inductive limit of $\clo_d \raro^{\lambda} \clo_d$ 
and $UHF_d \raro^{\lambda} UHF_d$ respectively and so for its inductive limit states.   

\vsp
Now we aim to prove (b). $H$ being a closed subgroup of $S^1$, it is either entire $S^1$ or a finite subgroup 
$\{exp({ 2i \pi l \over k})|l=0,1,...,k-1 \}$ where the integer $k \ge 1$. If $H=S^1$ we have nothing to prove 
for (b). When $H$ is a finite closed subgroup, we identify $[0,1)$ with $S^1$ by the usual map and note that 
if $\beta_t$ is restricted to $t \in [0,{1 \over k})$, then by scaling we check that $\beta_t$ defines a 
representation of $S^1$ in automorphisms of $\clo^H_d$. Now we consider the direct integral representation 
$\pi'$ defined by
$$\pi'= \int^{\oplus}_{[0,{1 \over k})} dt \pi_{|_{\clo^H_d}} \beta_t $$
of $\clo_d^{H}$ on $\clh_{|_{\clo_d^H}}\otimes L^2([0,{1 \over k})\;)$, where $\clh_{|_{\clo_d^H}}$ is the cyclic
space of $\pi(\clo_d^H)$ generated by $\Omega$. That it is indeed direct integral follows as states
$\psi \beta_{t_1}$ and $\psi \beta_{t_2}$ are either same or orthogonal for a factor state $\psi$ (see the proof of Lemma 7.4 in [BJKW] for further details).
Interesting point here to note that the new representation $\pi'$ 
is $(\beta_t)$ co-variant i.e. $\pi' \beta_t=\beta_t \pi' $, hence by simplicity of the $C^*$ algebra $\clo_d$ we 
conclude that $$\pi'(\mbox{UHF}_d)'' = \pi'(\clo_d^{H})''^{\beta_t}$$
 
\vsp
By exploring the hypothesis that $\psi$ is a factor state, 
we also have as in Lemma 6.11 in [BJKW] $I \otimes L^{\infty}([0,{1 \over k})\;) \subset \pi'(\clo_d^H)''$. 
Hence we also have
$$\pi'(\clo_d^H)''= \pi(\clo_d^H)'' \otimes L^{\infty}([0,{1 \over k})\;).$$ 
Since $\beta_t$ is acting as translation on $I \otimes L^{\infty}([0,{1 \over k})\;)$ which being 
an ergodic action, we have 
$$\pi'(\mbox{UHF}_d)'' = \pi(\clo_d^{H})'' \otimes 1$$
Since $\pi'(\mbox{UHF}_d)'' = \pi(\mbox{UHF}_d)'' \otimes 1$, we conclude 
that $\pi(\mbox{UHF}_d)'' = \pi(\clo_d^{H})''$. 

\vsp
A proof for the statement (c) follows from Lemma 7.12 in [BJKW]. The original idea of the proof can be 
traced back to Arveson's work on spectrum of an automorphism of a commutative compact group [Ar1]. 
\end{proof}

\vsp
Let $\omega'$ be an $\lambda$-invariant state on the $\mbox{UHF}_d$ sub-algebra of $\clo_d$. Following [BJKW, section 7], 
we consider the set 
$$K_{\omega'}= \{ \psi: \psi \mbox{ is a state on } \clo_d \mbox{ such that } \psi \lambda =
\psi \mbox{ and } \psi_{|\mbox{UHF}_d} = \omega' \}$$

By taking invariant mean on an extension of $\omega'$ to $\clo_d$, we verify that $K_{\omega'}$ is non empty and 
$K_{\omega'}$ is clearly convex and compact in the weak topology. In case $\omega'$ is an ergodic state ( extremal state )
$K_{\omega'}$ is a face in the $\lambda$ invariant states. Before we proceed to the next section here we recall Lemma 7.4 
of [BJKW] in the following proposition.

\vsp 
\begin{pro} Let $\omega'$ be ergodic. Then $\psi \in K_{\omega'}$ is an extremal point in
$K_{\omega'}$ if and only if $\psi$ is a factor state and moreover any other extremal point in $K_{\omega'}$
has the form $\psi \beta_z$ for some $z \in S^1$.
\end{pro}

\vsp
\begin{proof} 
Though Proposition 7.4 in [BJKW] appeared in a different set up, the same proof goes through for
the present case. We omit the details and refer to the original work for a proof.
\end{proof}

\section{ Dual Popescu system and pure translation invariant states: } 

\vsp
In this section we review the amalgamated Hilbert space developed in [BJKW] and prove a powerful criterion 
for a translation invariant factor state to be pure. Finally we will give proof of Theorem 1.1.  

\vsp
To that end let $\clm$ be a von-Neumann algebra acting on a Hilbert space $\clk$ and $\{v_k,\;1 \le k \le d\}$ 
be a family of bounded operators on $\clk$ so that $\clm=\{v_k,v_k^*,\;1 \le k \le d \}''$ and 
$\sum_kv_kv^*_k=1$. Furthermore let $\Omega$ be a cyclic and separating vector for $\clm$ so that the normal state 
$\phi(x)=<\Omega,x\Omega>$ on $\clm$ is invariant for the Markov map $\tau$ on $\clm$ defined by $\tau(x)=\sum_kv_kxv^*_k$ 
for $x \in \clm$. Let $\omega$ be the translation invariant state on UHF$_d =\otimes_{\IZ}M_d$ defined by
$$\omega(e^{i_1}_{j_1}(l) \otimes e^{i_2}_{j_2}(l+1) \otimes ....\otimes e^{i_n}_{j_n}(l+n-1)) = \phi(v_Iv^*_J)$$
where $e^i_j(l)$ is the elementary matrix at lattice site $l \in \IZ$. 

\vsp
We set $\tilde{v}_k = \overline{ \clj \sigma_{i \over 2}(v^*_k) \clj } \in \clm'$ ( see [BJKW] for details ) 
where $\clj$ and $\sigma=(\sigma_t,\;t \in \IR)$ are Tomita's conjugation operator and modular automorphisms 
associated with $\phi$.

By KMS or modular relation [BR vol-1] we verify that 
$$\sum_k \tilde{v}_k \tilde{v}_k^*=1$$ 
and 
\be
\phi(v_Iv^*_J)= \phi(\tilde{v}_{\tilde{I}}\tilde{v}^*_{\tilde{J}})
\ee
where $\tilde{I}=(i_n,..,i_2,i_1)$ if $I=(i_1,i_2,...,i_n)$. Moreover $\tilde{v}^*_I\Omega = 
\clj \sigma_{i \over 2}(v_{\tilde{I}})^*\clj\Omega= \clj \Delta^{1 \over 2}v_{\tilde{I}}\Omega
=v^*_{\tilde{I}}\Omega$. We also set $\tilde{\clm}$ to be the von-Neumann algebra generated by 
$\{\tilde{v}_k: 1 \le k \le d \}$. Thus $\tilde{\clm} \subseteq \clm'$. A major problem that we 
will have to address is: when equality holds and its relation to Haag duality property (16).  

\vsp
Let $(\clh,P,S_k,\;1 \le k \le d )$ and $(\tilde{\clh},P,\tilde{S}_k,\;1 \le k \le d)$ be the Popescu dilation described 
as in Theorem 2.1 associated with $(\clk,v_k,\;1 \le k \le d)$ and $\clk,\tilde{v}_k,\;1 \le k \le d)$ respectively. 
Following [BJKW] we consider the amalgamated tensor product $\clh \otimes_{\clk} \tilde{\clh}$ of $\clh$ with 
$\tilde{\clh}$ over the joint subspace $\clk$. It is the completion of the quotient of the set $$\IC \bar{I} \otimes 
\IC I \otimes \clk,$$ where $\bar{I},I$ both consist of all finite sequences with elements in $\{1,2, ..,d \}$, 
by the equivalence relation defined by a semi-inner product defined on the set by requiring
$$<\bar{I} \otimes I \otimes f,\bar{I}\bar{J} \otimes IJ \otimes g>=<f,\tilde{v}_{\bar{J}}v_Jg>,$$
$$<\bar{I}\bar{J} \otimes I \otimes f, \bar{I} \otimes IJ \otimes g> =<\tilde{v}_{\bar{J}}f,v_Jg>$$
and all inner product that are not of these form are zero. We also define
two commuting representations $(S_i)$ and $(\tilde{S}_i)$ of $\clo_d$ on
$\clh \otimes_{\clk} \tilde{\clh}$ by the following prescription:
$$S_I\lambda(\bar{J} \otimes J \otimes f)=\lambda(\bar{J} \otimes IJ \otimes f),$$
$$\tilde{S}_{\bar{I}}\lambda(\bar{J} \otimes J \otimes f)=\lambda(\bar{J}\bar{I} \otimes J \otimes f),$$
where $\lambda$ is the quotient map from the index set to the Hilbert space. Note that the subspace generated by
$\lambda(\emptyset \otimes I \otimes \clk)$ can be identified with $\clh$ and earlier $S_I$ can be identified
with the restriction of $S_I$ defined here. Same is valid for $\tilde{S}_{\bar{I}}$. The subspace $\clk$ is
identified here with $\lambda(\emptyset \otimes \emptyset \otimes \clk)$. 
Thus $\clk$ is a cyclic subspace for the representation $$\tilde{s}_j \otimes s_i \raro \tilde{S}_j S_i$$ 
of $\tilde{\clo}_d \otimes \clo_d$ in the amalgamated Hilbert space. Let $P$ be the projection onto $\clk$. Then we
have 
$$S_i^*P=PS_i^*P=v_i^*$$
$$\tilde{S}_i^*P=P\tilde{S}_i^*P=\tilde{v}^*_i$$
for all $1 \le i \le d$. 

\vsp
We start with a simple proposition.

\begin{pro} 
The following holds: 

\NI (a) For any $1 \le i,j \le d$ and $|I|,|J|< \infty$ and $|\bar{I}|,|\bar{J}| < \infty$
$$<\Omega,\tilde{S}_{\bar{I}}\tilde{S}^*_{\bar{J}} S_iS_IS^*_JS^*_j \Omega>=<\Omega, 
\tilde{S}_i \tilde{S}_{\bar{I}}\tilde{S}^*_{\bar{J}}\tilde{S}^*_jS_IS^*_J \Omega>;$$

\NI (b) The vector state $\psi_{\Omega}$ on $$\tilde{\mbox{UHF}}_d \otimes \mbox{UHF}_d \equiv
\otimes_{-\infty}^0 M_d \otimes_1^{\infty}M_d \equiv \otimes_{\IZ} M_d$$ is equal to $\omega$;

\NI (c) $\pi(\tilde{\clo}_d \otimes \clo_d)''= \clb(\tilde{\clh} \otimes_{\clk} \clh)$ if and only if
$\{x \in \clb(\clk): \tau(x)=x,\; \tilde{\tau}(x)=x \}= \{zI: z \in  \IC \}$. 
\end{pro} 

\vsp
\begin{proof} 
By our construction $\tilde{S}^*_i\Omega =\tilde{v}^*_i\Omega=v_i^*\Omega=S_i^*\Omega$. Now (a) and (b) 
follow by repeated application of $\tilde{S}^*_i\Omega=S^*_i\Omega$ and
the commuting property of the two representation $\pi(\clo_d \otimes I)$ and $\pi(I \otimes \tilde{\clo}_d)$. The last
statement (c) follows from a more general fact proved below
that the commutant of $\pi(\clo_d \otimes \tilde{\clo}_d)''$ is
order isomorphic with the set $\{x \in \clb(\clk): \tau(x)=x,\; \tilde{\tau}(x)=x \}= \{zI: z \in  \IC \}$ via the 
map $X \raro PXP$ where $X$ is the weak$^*$ limit of $\{\Lambda^m\tilde{\Lambda}^n(x)$ as $(m,n) \raro 
(\infty,\infty)$. For details let $Y$ be the strong limit of increasing sequence of projections 
$(\Lambda \tilde{\Lambda})^n(P)$ as $n \raro \infty$. Then $Y$ commutes with $S_i\tilde{S}_j, S^*_i\tilde{S}^*_j$ 
for all $1 \le i,j \le d$. As $\Lambda(P)) \ge P$, we also have $\Lambda(Y) \ge Y$. Hence $(1-Y)S_i^*Y=0$. As $Y$ commutes 
with $S_i\tilde{S}_j$ we get $(1-Y)S_i^*S_i\tilde{S}_jY=0$ i.e. $(1-Y)\tilde{S}_jY=0$ for all $1 \le j \le d$. By symmetry 
of the argument we also get $(1-Y)S_iY=0$ for all $1 \le i \le d$. Hence $Y$ commutes with $\pi(\clo_d)''$ and by symmetry of 
the argument $Y$ commutes as well with $\pi(\tilde{\clo}_d)''$. As $Yf=f$ for all $f \in \clk$ and $\clk$ is cyclic for the 
representation $\pi(\tilde{\clo}_d \otimes \clo_d)$ we conclude that $Y=I$ on $\tilde{\clh} \otimes_{\clk} \clh$.

\vsp
Let $x \in \clb(\clk)$ so that $\tau(x)=x$ and $\tilde{\tau}(x) = x $ then as in the proof of Theorem 2.1 
we also check that $(\Lambda \tilde{\Lambda})^k(P) \Lambda^m \tilde{\Lambda}^n(x) (\Lambda \tilde{\Lambda})^k(P)$ is 
independent of $m,n$ as long as $m,n \ge k$. Hence the weak$^*$ limit $\Lambda^m \tilde{\Lambda}^n(x) \raro X$ exists 
as $m,n \raro \infty$. Furthermore the limiting element $X \in \pi(\clo_d \otimes \tilde{\clo}_d)'$ and $PXP=x$. That the 
map $X \raro PXP$ is an order-isomorphic on the set of self adjoint elements follows as in Theorem 2.1. This completes the proof.
\end{proof} 

\vsp
In short Proposition 3.1 also says that $(\tilde{\clh} \otimes_{\clk} \clh,S_i\tilde{S}_j\;1 \le i,j \le d,P)$ 
is the Popescu dilation associated with Popescu elements $(\clk,v_i\tilde{v_j},1 \le i,j \le d \}$. Now we 
will be more specific in our starting Popescu's elements in order to explore the representation $\pi$ of 
$\tilde{\clo}_d \otimes \clo_d$ in the amalgamated Hilbert space $\tilde{\clh} \otimes_{\clk} \clh$.  

\vsp 
Let $\omega$ be a translation invariant factor state on $\clb$ and $\omega'$ be its restriction to $\clb_R$ which we identified with 
$\mbox{UHF}_{d}$ with respect to an orthonormal basis $(e_i)$ of $\IC^d$ (see statement before equation (4). Let $\psi$ be an extremal point in $K_{\omega'}$. 
We consider the Popescu's elements $(\clk,\clm,v_k,1 \le k \le d,\Omega)$ described as in Proposition 2.4 associated with support projection of the state $\psi$ 
in $\pi_{\psi}(\clo_d)''$ and also consider associated dual Popescu's elements $(\clk,\tilde{\clm},\tilde{v_k},\;1 \le k \le d)$ where $\tilde{\clm}$ is the von-Neumann algebra generated 
by $\{ \tilde{v}_k:\; 1 \le k \le d \}$. Thus in general $\tilde{\clm} \subseteq \clm'$ and an interesting question: when do we have $\clm'=\tilde{\clm}$? Going back to our starting 
example of unique KMS state for the automorphisms $\beta_t(s_i)= ts_i,\;t \in S^1$, we check that $v^*_k=S^*_k$,  
$\clj \tilde{v}^*_k \clj = {1 \over d} S_k$ and thus equality holds i.e. $\tilde{\clm}=\clm'$. But 
the corner vector space $\tilde{\clm}_c=P\pi(\tilde{\clo}_d)''P$ generated by the elements 
$\{\tilde{v}_I\tilde{v}^*_J: |I|,|J| < \infty \}$ fails to be an algebra. Thus two questions sounds reasonable 
here.  

\vsp 
\NI (a) Does the equality $\clm'=\tilde{\clm}$ hold in general for an extremal element $\psi \in K_{\omega'}$ and 
a factor state $\omega$?  

\vsp 
\NI (b) When can we expect $\tilde{\clm}_c$ to be a $*$-algebra and so equal to $\tilde{\clm}$?  

\vsp
The dual condition on support projection and equality $\tilde{\clm}=\clm'$ are rather deep and will lead us to a far reaching consequence on the state $\omega$. 
In the paper [BJKW] these two conditions are implicitly assumed to give a criterion for a translation invariant factor state to be pure. Apart from this refined interest, 
we will address the converse problem that turns out to be crucial for our main results. In the following we prove a crucial step towards that goal fixing the basic structure 
which will be repeatedly used in the computation using Cuntz relations. 

\vsp
\begin{pro} 
Let $\omega$ be a translation invariant factor state on $\clb$ and $\psi$ be an extremal point in $K_{\omega'}$. 
We consider the amalgamated representation $\pi$ of $\tilde{\clo}_d \otimes
\clo_d$ in $\tilde{\clh} \otimes_{\clk} \clh$ where the Popescu's elements $(\clk,\clm,v_k,\;1 \le k \le d)$ are 
taken as in Proposition 2.4. Then the following statements hold:

\vsp
\NI (a) $\pi(\tilde{\clo}_d \otimes \clo_d)''= \clb(\tilde{\clh} \otimes_{\clk} \clh)$. Furthermore $\pi(\clo_d)''$ 
and $\pi(\tilde{\clo}_d)''$ are factors and the following sets are equal:

\NI (i) $H=\{ z \in S^1: \psi \beta_z = \psi \}$;

\NI (ii) $H_{\pi}  = \{ z : \beta_z \mbox{ extends to an automorphisms of } \pi(\clo_d)'' \}$;

\NI (iii) $\tilde{H}_{\pi}= \{ z : \beta_z \mbox{ extends to an automorphisms of } \pi(\tilde{\clo}_d)'' \}$.
Moreover $\pi(\tilde{\mbox{UHF}}_d \otimes I)''$ and $\pi(I \otimes \mbox{UHF}_d)''$ are factors.

\vsp 
\NI (b) $z \raro U_z$ is the unitary representation of $H$ in the Hilbert space $\tilde{\clh} \otimes_{\clk} \clh$ 
defined by
$U_z(\pi(\tilde{s}_j \otimes s_i )\Omega=\pi(z\tilde{s}_j \otimes zs_i)\Omega$
 
\vsp  
\NI (c) The commutant of $\pi(\tilde{\mbox{UHF}}_d \otimes \mbox{UHF}_d )''$ is invariant by the canonical endomorphisms
$\Lambda(X)=\sum_i S_iXS_i^*$ and $\tilde{\Lambda}(X) = \sum_i \tilde{S}_iX \tilde{S}^*_i$. Same is true 
for each $i$ that the surjective map $X \raro S^*_iXS_i$ keeps the commutant of $\pi(\tilde{\mbox{UHF}}_d \otimes 
\mbox{UHF}_d)''$ invariant. Same holds for the map $X \raro \tilde{S}^*_iX\tilde{S}_i$.

\vsp 
\NI (d) The centre of $\pi(\tilde{\mbox{UHF}}_d \otimes \mbox{UHF}_d )''$ is invariant by the canonical endomorphisms
$\Lambda(X)=\sum_i S_iXS_i^*$ and $\tilde{\Lambda}(X) = \sum_i \tilde{S}_iX \tilde{S}^*_i$. Moreover for each $i$, the
surjective map $X \raro S^*_iXS_i$ keeps the centre of $\pi(\tilde{\mbox{UHF}}_d \otimes 
\mbox{UHF}_d)''$ invariant. Same holds for the map $X \raro \tilde{S}^*_iX\tilde{S}_i$.
\end{pro} 

\vsp
\begin{proof} 
$P$ being the support projection by Proposition 2.4 we have $\{x \in \clb(\clk): 
\sum_k v_k x v_k^*= x\} = \clm'$. That $(\clm',\tilde{\tau},\phi)$ is ergodic follows from a general 
result [Mo1] ( see also [BJKW] for a different proof ) as $(\clm,\tau,\phi)$ is ergodic for a factor state
$\psi$ being an extremal element in $K_{\omega'}$ (Proposition 2.6). Hence $\{x \in \clb(\clk): \tau(x)=\tilde{\tau}(x)=x \} = \IC$. Hence by 
Proposition 3.1, we conclude that $\pi(\tilde{\clo}_d \otimes \clo_d)''= \clb(\tilde{\clh} \otimes_{\clk} \clh)$. That both $\pi(\clo_d)''$ 
and $\pi(\tilde{\clo}_d)''$ are factors follows trivially as $\pi(\tilde{\clo}_d \otimes \clo_d)''=\clb(\tilde{\clh} \otimes_{\clk} \clh)$ 
and $\pi(\clo_d)'' \subseteq \pi(\tilde{\clo}_d)'$. 

\vsp
By the discussion above we first recall that $\Omega$ is a cyclic vector for the representation of 
$\pi(\tilde{\clo}_d \otimes \clo_d)$. Let $G = \{ z=(z_1,z_2) \in S^1 \times S^1: \beta_z \mbox{
extends to an automorphism }$ $\mbox{on}\; \pi(\tilde{\clo}_d \otimes \clo_d)'' \}$ be the closed 
subgroup where
$$\beta_{(z_1,z_2)}(\tilde{s}_j \otimes s_i)=z_1 \tilde{s}_j \otimes z_2s_i .$$
By repeated application of the fact that $\pi(\clo_d)''$ commutes with $\pi(\tilde{\clo}_d)''$ and 
$S_i^*\Omega=\tilde{S}^*_i\Omega$ as in Proposition 3.1 (a) we verify that $\psi \beta_{(z,z)}=\psi$ 
on $\clo_d \otimes \tilde{\clo}_d$ if $z \in H$. For $z \in H $ we set unitary operator $U_z \pi(x \otimes y) \Omega 
= \pi(\beta_z(x) \otimes \beta_z(y))\Omega$ for all $x \in \tilde{\clo}_d$ and $y \in \clo_d$. 
Thus we have $U_z\pi(s_i)U_z^*=z\pi(s_i)$ and also $U_z\pi(\tilde{s}_i)U_z^*=z\tilde{s}_i.$  By taking 
its restriction to $\pi(\clo_d)''$ and $\pi(\tilde{\clo}_d)''$ respectively we check that
$H \subseteq \tilde{H}_{\pi} $ and $H  \subseteq H_{\pi}$.

\vsp
For the converse let $z \in H_{\pi}$ and we use the same symbol $\beta_z$ for the extension to an automorphism of
$\pi(\clo_d)''$. By taking the inverse map we check easily that $\bar{z} \in H_{\pi}$ and in fact $H_{\pi}$ is
a subgroup of $S^1$. Since $\lambda$ commutes with $\beta_z$ on $\clo_d$, the canonical endomorphism $\Lambda$ defined
by $\Lambda(X) = \sum_k S_kXS_k^*$ also commutes with the extension of $\beta_z$ on $\pi(\clo_d)''$. Note that the map
$\pi(x)_{{|}_{\clh}} \raro \pi(\beta_z(x))_{{|}_{\clh}} $ for $x \in \clo_d$ is a well defined $*$-homomorphism.
Since same is true for $\bar{z}$ and $\beta_z\beta_{\bar{z}}=I$, the map is an isomorphism. Hence $\beta_z$ extends
uniquely to an automorphism of $\pi(\clo_d)''_{{|}_{\clh}}$ commuting with the restriction of the canonical endomorphism
on $\pi(\clo_d)''_{|\clh}$. Since $\pi(\clo_d)''_{{|}_{\clh}}$ is a factor, we conclude as in Proposition 2.5 (a) that 
$z \in H$. Thus $H_{\pi} \subseteq H$. As $\pi(\tilde{\clo}_d)''$ is also a factor, we also have $\tilde{H}_{\pi} \subseteq H$. 
Hence we have $H=H_{\pi}=\tilde{H}_{\pi}$ and $\{(z,z): z \in H \} \subseteq G \subseteq H \times H$. 

\vsp
For the second part of (a) we will adopt the argument used for Proposition 2.5. To that end we first note that $\Omega$
being a cyclic vector for the representation $\tilde{\clo}_d \otimes \clo_d$ in the Hilbert space $\tilde{\clh} \otimes_{\clk} \clh $, 
by Lemma 7.11 in [BJKW] (note that the proof only needs the cyclic property ) the representation of UHF$_d$ on $\tilde{\clh} \otimes_{\clk} \clh $ 
is quasi-equivalent to its sub-representation on the cyclic space generated by $\Omega$. On the other hand, by our hypothesis that $\omega$ is a 
factor state, Power's theorem [Pow] ensures that the state $\omega'$ (i.e. the restriction of $\omega$ to $\clb_R$ which is identified here with $\mbox{UHF}_d$ ) 
is also a factor state on $\mbox{UHF}_d$. Hence quasi-equivalence ensures that $\pi(I \otimes \mbox{UHF}_d)''$ is a factor. We also note that
the argument used in Lemma 7.11 in [BJKW] is symmetric i.e. same argument is also valid for $\tilde{\mbox{UHF}}_d$.
Thus $\pi(\tilde{\mbox{UHF}}_d \otimes I)''$ is also a factor. This completes the proof of (a). We have proved (b) 
while giving proof of (a).

\vsp
For $X \in \clb(\tilde{\clh} \otimes_{\clk} \clh)$, as $\Lambda(X)$ commutes with $\pi(\lambda(\tilde{\mbox{UHF}}_d \otimes \mbox{UHF}_d))''$ and
$\{ S_iS^*_j: 1 \le i,j \le d \}$ we verify by Cuntz's relation that $\Lambda(X)$ is also an 
element in the commutant of $\pi(\lambda(\tilde{\mbox{UHF}}_d \otimes \mbox{UHF}_d))''$ once 
$X$ is so. It is also obvious that $\Lambda(X)$ is an element in $\pi(\tilde{\mbox{UHF}}_d \otimes \mbox{UHF}_d)''$ if 
$X$ is so. Thus $\Lambda(X)$ is an element in the commutant/centre of $\pi(\lambda(\tilde{\mbox{UHF}}_d \otimes \mbox{UHF}_d)''$ once $X$ is so. For the last statement consider the map $X \raro S_i^*XS_i$
on $\pi(\tilde{\mbox{UHF}}_d \otimes \mbox{UHF}_d)''$ which is clearly onto by Cuntz relation (3). Hence 
we need to show that $S^*_iXS_i$ is an element in the commutant whenever $X$ is so. To that end note that $S_i^*XS_i S_i^*YS_i=S_i^*S_iS_i^*XYS_i=S_i^*YXS_i=S^*_iYS_iS_i^*XS_i$ since $X$ commutes with $S_iS_i^*$. 
Thus onto property of the map ensures that $S_i^*XS_i$ is an element in the commutant/centre of $\pi(\tilde{\mbox{UHF}}_d \otimes \mbox{UHF}_d)''$ once $X$ is so. This completes the proof of (c) and (d). 
\end{proof} 

\vsp
One interesting problem here how to describe the von-Neumann algebra $\cli$ which consists of invariant elements of the gauge action $\{ \beta_z: z \in H \}$ in $\clb(\tilde{\clh} \otimes_{\clk} \clh)$. A general result 
due to E. Stormer [So] says that the algebra of invariant elements is a von-Neumann algebra of type-I with 
centre completely atomic. Here the situation is much simple because we know explicitly 
that $\cli=\{U_z: z \in H \}'$ and we write 
spectral decomposition as 
\be 
U_z=\sum_{ k \in \hat{H} } z^k F_k
\ee 
for $z \in H$, $\hat{H}$ is the dual group of $H$, either $\hat{H}=\{z:z^n=1 \}$ or $\IZ$. Thus the centre of $\cli$ is equal to $\{F_k: k \in \hat{H} \}$. 

\vsp
As a first step we describe the center $\clz$ of $\pi(\tilde{\mbox{UHF}}_d \otimes \mbox{UHF}_d)''$ 
by exploring Cuntz relation, that it is also non-atomic even for a factor state $\omega$. In fact we
will show that the centre $\clz$ is a sub-algebra of the centre of $\cli$. In the following proposition 
we give an explicit description. 

\vsp
\begin{pro}  
Let $\omega,\psi$ be as in Proposition 3.2 with Popescu system $(\clk,\clm,$ $v_k,\Omega)$ be 
taken as in Proposition 2.4 i.e. on support projection. Then the centre of 
$\pi(\tilde{\mbox{UHF}}_d \otimes \mbox{UHF}_d)''$ is completely atomic and the element 
$E_0=[\pi(\tilde{\mbox{UHF}}_d \otimes \mbox{UHF}_d)' \vee \pi(\tilde{\mbox{UHF}}_d \otimes 
\mbox{UHF}_d)''\Omega]$ is a minimal projection in the centre of $\pi(\tilde{\mbox{UHF}}_d \otimes \mbox{UHF}_d)''$ and the centre is invariant for both $\Lambda$ and $\tilde{\Lambda}$. Furthermore the 
following holds:

\vsp 
\NI (a) The centre of $\pi(\tilde{\mbox{UHF}}_d \otimes \mbox{UHF}_d)''$ has the following two disjoint
possibilities:

\vsp 
\NI (i) There exists a positive integer $m \ge 1$ such that the centre is generated by the family of minimal
orthogonal projections $\{ \Lambda_k(E_0): 0 \le k \le m-1 \}$ where $m \ge 1$ is the least positive integer so that
$\Lambda^m(E_0)=E_0$. In such a case $\{z: z^m = 1 \} \subseteq H$;

\vsp 
\NI (ii) The family of minimal nonzero orthogonal projections $\{ E_k: k \in \IZ \}$ where $E_k= \Lambda^k(E_0)$
for $k \ge 0$ and $E_k= S^*_IE_0S_I$ for $k < 0$ where $|I|=-k$ and independent of multi-index $I$ generates
the centre and $H=S^1$; 

\vsp
\NI (b) $\Lambda(E)=\tilde{\Lambda}(E) $ for any $E$ in the centre of $\pi(\tilde{\mbox{UHF}}_d \otimes \mbox{UHF}_d)''$

\vsp 
\NI (c) If $\Lambda(E_0)=E_0$ then $E_0=1$.

\end{pro} 

\vsp
\begin{proof} 
Let $E' \in \pi(\tilde{\mbox{UHF}}_d \otimes \mbox{UHF}_d)'$ be the projection onto the subspace generated
by the vectors $\{ S_IS^*_J\tilde{S}_{I'}S^*_{J'}\Omega,\; |I|=|J|,|I'|=|J'| < \infty \}$ and $\pi_{\Omega}$ be the restriction
of the representation $\pi$ of $\tilde{\mbox{UHF}}_d \otimes \mbox{UHF}_d$ to the cyclic subspace $\clh_{\Omega}$ generated by
$\Omega$. Identifying $\clb$ with $\tilde{\mbox{UHF}}_d \otimes \mbox{UHF}_d$ we check that $\pi_{\omega}$ is unitary
equivalent to $\pi_{\Omega}$. Thus $\pi_{\Omega}$ is a factor representation. 

\vsp
For any projection $E$ in the centre of $\pi(\tilde{\mbox{UHF}}_d \otimes \mbox{UHF}_d)''$, via the unitary equivalence, we note 
that $EE'=E'EE'$ is an element in the centre of $\pi_{\Omega}(\tilde{\mbox{UHF}}_d \otimes \mbox{UHF}_d)''$. $\omega$ being a 
factor state we conclude that $EE'$ is a scalar multiple of $E'$ and so we have
\be
EE'=\omega(E)E'
\ee
Thus we also have $EYE'=\omega(E)YE'$ for all $Y \in \pi(\tilde{\mbox{UHF}}_d \otimes \mbox{UHF}_d)'$ and so 
\be
EE_0=\omega(E)E_0
\ee
 
\vsp
Since $EE'$ is a projection and $E' \ne 0$, we have $\omega(E)=\omega(E)^2$. Thus $\omega(E)=1$ or $0$. So, for such an element $E$,
the following is true:

\NI (i) If $E \le E_0$ then either $E=0$ or $E=E_0$ i.e. $E_0$ is a minimal projection in the centre of $\pi(\tilde{\mbox{UHF}}_d \otimes \mbox{UHF}_d)''$

\NI (ii) $\omega(E)=1$ if and only if $E \ge E_0$ 

\NI (iii) $\omega(E)=0$ if and only if $EE_0=0$.

\vsp
As $\Lambda(E_0)$ is a projection in the centre of $\pi(\tilde{\mbox{UHF}}_d \otimes \mbox{UHF}_d)''$ by our last proposition i.e. Proposition 
3.2 (c), we have either $\omega(\Lambda(E_0))=1$ or $0$. Since $\Lambda(E_0) \neq 0$ by the injective property of the endomorphism, 
we have either $\Lambda(E_0) \ge E_0$ or $\Lambda(E_0)E_0=0$. In case $\Lambda(E_0) \ge E_0$ we have $S^*_iE_0S_i \le S^*_i\Lambda(E_0)S_i=E_0$ 
for all $1 \le i \le d$. If so, $S^*_iE_0S_i$ being a non-zero projection in the centre of $\pi(\mbox{UHF}_d \otimes \tilde{\mbox{UHF}}_d)''$ (Proposition 3.2 (c) ), 
by (i) we have $E_0=\Lambda(E_0)$. Thus we have either $\Lambda(E_0)=E_0$ or $\Lambda(E_0)E_0=0$. 

\vsp
If $\Lambda(E_0)E_0=0$, we have $\Lambda(E_0) \le I-E_0$ and by Cuntz's relation we check that
$E_0 \le I - S^*_iE_0S_i$ and $S^*_jS^*_iE_0S_iS_j \le I- S^*_jE_0S_j$ for all $1 \le i,j \le d$.
So we also have $E_0S^*_jS^*_iE_0S_iS_jE_0 \le E_0- E_0S^*_jE_0S_jE_0=E_0$. Thus we have either
$E_0S^*_jS^*_iE_0S_iS_jE_0= 0$ or $E_0S^*_jS^*_iE_0S_iS_jE_0=E_0$ as $S^*_jS^*_iE_0S_iS_j$ is an element
in the centre by Proposition 3.2 (c). So either we have $\Lambda^2(E_0)E_0=0$ or $\Lambda^2(E_0) \le E_0$. 
$\Lambda$ being an injective map we either have $\Lambda^2(E_0)E_0=0$ or $\Lambda^2(E_0)= E_0.$

\vsp
More generally we check that if $\Lambda(E_0)E_0=0, \Lambda^2(E_0)E_0=0,..\Lambda^k(E_0)E_0=0$ for some $k \ge 1$
then either $\Lambda^{k+1}(E_0)E_0=0$ or $\Lambda^{k+1}(E_0)=E_0$. To verify that first we check that in such a
case $E_0 \le I- S^*_IE_0S_I$ for all $|I|=n$ and then following the same steps as before to check that $S^*_iS^*_IE_0S_IS_i
\le I-S^*E_0S_i$ for all $i$. Thus we have $E_0S^*_iS^*_IE_0S_IS_iE_0 \le E_0$ and arguing as before we complete the proof
of the claim that either $\Lambda^{k+1}(E_0)E_0=0$ or $\Lambda^{k+1}(E_0)=E_0.$ 

\vsp
We summarize now by saying that $E_0,\Lambda(E_0),..,\Lambda^{m-1}(E_0)$ are mutually orthogonal projections with $m \ge 1$ possibly 
be infinite, if not then $\Lambda^m(E_0)=E_0$. 

\vsp
Let $\pi_k,\;k \ge 0$ be the representation $\pi$ of $\tilde{\mbox{UHF}}_d \otimes \mbox{UHF}_d$ restricted to the
subspace $\Lambda^k(E_0)$. The representation $\pi_0$ of $\tilde{\mbox{UHF}}_d \otimes \mbox{UHF}_d$ is 
isomorphic to the representation $\pi$ of $\tilde{\mbox{UHF}}_d \otimes \mbox{UHF}_d$ restricted to $E'$ and thus
quasi-equivalent. For a general discussion on quasi-equivalence we refer to section 2.4.4 in [BR vol-1]. Since $\omega$ is a factor state, 
$\pi_0$ is a factor representation. We claim now that each $\pi_k$ is a factor representation. We fix any $k \ge 1$ and let $X$ be an element 
in the centre of $\pi_k(\mbox{UHF}_d \otimes \tilde{\mbox{UHF}}_d)$. Then for any $|I|=k$, $S^*_IE_kS_I=E_0$ and so $S^*_IXS_I$ is an element 
in the centre of $\pi_0(\mbox{UHF}_d \otimes \tilde{\mbox{UHF}}_d)$ by Proposition 3.2 (d). Further $S_I^*XS_I= S^*_IXS_IS^*_JS_J = S^*_JXS_J$ 
for all $|J|=|I|=k$. $\pi_0$ being a factor representation, we have $S^*_IXS_I= c E_0$ for some scalar $c$ independent of the multi-index we choose $|I|=k$. 
Hence $c \Lambda_k(E_0) = \sum _{|J|=k } S_JS^*_IXS_IS_J^*=\sum _{|J|=k } S_JS^*_IS_IS_J^*X= X$ as $X$ is an element in the centre of $\pi(\tilde{\mbox{UHF}}_d \otimes \mbox{UHF}_d)$. 
Thus for each $k \ge 1$, $\pi_k$ is a factor representation as $\pi_0$ is so.  

\vsp
We also note that $\Lambda(E_0)\tilde{\Lambda}(E_0) \neq 0$. Otherwise we have $<S_i\Omega,\tilde{S}_j\Omega>=0$
for all $i,j$ and so $<\Omega,\tilde{S}_jS^*_i\Omega>=0$ for all $i,j$ as $\pi(\clo_d)''$ commutes with 
$\pi(\tilde{\clo}_d)''$. However $\tilde{S}^*_i\Omega=S^*_i\Omega$ and $\sum_i \tilde{S}_i\tilde{S}^*_i=1$ 
which leads a contradiction. Hence $\Lambda(E_0)\tilde{\Lambda}(E_0) \neq 0$. As $\pi$ restricted to 
$\Lambda(E_0)$ is a factor state and both $\Lambda(E_0)$ and $\tilde{\Lambda}(E_0)$ are elements in the 
centre of $\pi(\tilde{\mbox{UHF}}_d \otimes \mbox{UHF})''$, by Proposition 3.2 (d), 
we conclude that $\Lambda(E_0)=\tilde{\Lambda}(E_0).$ Using the commuting property of the endomorphisms $\Lambda$ 
and $\tilde{\Lambda}$, we verify by a simple induction method that $\Lambda^k(E_0) = \tilde{\Lambda}^k(E_0)$ 
for all $k \ge 1$. Thus the sequence of orthogonal projections $E_0, \tilde{\Lambda}(E_0),...$, are also periodic with same 
period or aperiodic according as the sequence of orthogonal projections $E_0,\Lambda(E_0),...$ is. 

\vsp
If $\Lambda^m(E_0)=E_0$ for some $m \ge 1$ then we check that $\sum_{0 \le k \le m-1} \Lambda^k(E_0)$ is a 
$\Lambda$ and as well $\tilde{\Lambda}$-invariant projection and thus equal to $1$ by the cyclic property of 
$\Omega$ for $\pi(\clo_d \otimes \tilde{\clo}_d)''$. In such a case we set $V_z=\sum_{0 \le k \le m-1} z^kE_k$ for $z \in S^1$ 
for which $z^m=1$ and check that $\Lambda(V_z)= \sum_{0 \le k \le m-1}z^k\Lambda(E_k)=\sum_{0 \le k \le m-1} z^k E_{k+1}= \bar{z} V_z$ 
where $E_m=E_0$ and so by the Cuntz relations we have $V^*_zS_iV_z= \bar{z}S_i$ for all $1 \le i \le d$. Following the same 
steps we also have $\tilde{\Lambda}(V_z)=\bar{z}V_z$ and so $V^*_z\tilde{S}_iV^*_z=\bar{z}\tilde{S}_i$ for $ 1 \le i \le d$. 
Thus $V_z=U_z$ for all $z \in H_0 = \{z: z^m=1 \} \subseteq H$. 

\vsp
Now we consider the case where $E_0,\Lambda(E_0),..\Lambda^k(E_0),..$ is a sequence of aperiodic orthogonal
projections. We extend the family of projections $\{E_k:\;k \in \IZ \}$ to all integers by 
 
$$E_k=\Lambda^k(E_0)\;\; \mbox{for all}\;\; k \ge 1$$
and
$$E_k = S^*_IE_0S_I\;\; \mbox{for all}\;\;k \le 1,\;\;\mbox{where}\;\; |I|=-k$$
We claim that the definition of $\{ E_k;\;k \le -1 \}$ does depends only on length of the multi-index 
$I$ that we choose. We may choose any other $J$ so that $|J|=|I|$ and check the following identity: 
$$S^*_IE_0S_I = S^*_IE_0S_IS^*_JS_J=S^*_IS_IS^*_JE_0S_J=S^*_JE_0S_J$$ 
where $E_0$, being an element in the centre of $\pi(\tilde{\mbox{UHF}}_d \otimes \mbox{UHF}_d)''$, 
commutes with $S_IS_J^*$ as $|I|=|J|$. Further $\Lambda^k(E_0)=\tilde{\Lambda}^k(E_0)$ ensures that 
$S_I\tilde{S}_J^*$ commutes with $E_0$ for all $|I|=|J|=k$ and $k \ge 1$. Hence we also have 
$$E_{-k}=S^*_IE_0S_I\tilde{S}^*_J\tilde{S}_J=\tilde{S}^*_JE_0\tilde{S}_J$$ 
for all $|J|=|I|=k$ and $k \ge 1$. Now we claim that 
$$\Lambda(E_k) = \tilde{\Lambda}(E_k) = E_{k+1}$$ 
for all $k \in \IZ$. For $k \ge 0$ we have nothing to prove. For $k \le -1$ we check that the following steps
$$\Lambda(S^*_IE_0S_I)$$
$$=\sum_j S_jS^*_iS^*_{I'}E_0S_{I'}S_iS_j^*$$
$$=\sum_jS^*_{I'}E_0S_{I'} S_jS^*_iS_iS_j^* = S^*_{I'}E_0S_{I'}$$
where we have $I=(I',i)$ and used the fact that elements $S_jS^*_i$ commutes with $\{ E_k: k \in \IZ \}$, which are elements in the centre of 
$\pi(\mbox{UHF}_d \otimes \tilde{\mbox{UHF}}_d)''$. For a proof that 
$\tilde{\Lambda}(E_k)=E_{k+1}$ we may follow the same steps as 
$E_k=\tilde{S}^*_IX\tilde{S}_I$ where $|I|=-k$ and $k \le -1$. 

\vsp
We also claim that $\{ E_k: k \in \IZ \}$ is an orthogonal family of non-zero projections. To that end we choose any two elements
say $E_k,E_m,\;k \ne m$ and use endomorphism $\Lambda^n$ for $n$ large enough so that both $n+k \ge 0 , n+m  \ge 0$
to conclude that $\Lambda^n(E_kE_m) = E_{k+n}E_{k+m}=0$ as $k+n \ne k+m$. $\Lambda$ being an injective map we get
the required orthogonal property. Thus $\sum_{ k \in \IZ} E_k$ being an invariant projection for both $\Lambda$ and
$\tilde{\Lambda}$ we get by cyclic property of $\Omega$ that $\sum_{ k \in \IZ} E_k=I$. Let $\pi_k,\; k \le -1$ be the representation $\pi$ 
of $\tilde{\mbox{UHF}}_d \otimes \mbox{UHF}_d$ restricted to the subspace $E_k$. Going along the same line as above, we verify that for 
each $k \le -1$, $\pi_k$ is a factor representation of $\tilde{\mbox{UHF}}_d \otimes \mbox{UHF}_d$. We also set $V_z= \sum_{-\infty < k 
< \infty} z^k E_k$ for all $z \in S^1$ and check that $\Lambda(V_z)=\bar{z}V_z$ and also $\tilde{\Lambda}(V_z)=\bar{z}V_z$. Hence $S^1=H$ 
as $H$ is a closed subset of $S^1$. This completes the proof of (a). Proof of (b) and (c) are now simple consequence of the proof of (a). 
\end{proof} 

\vsp
It is clear that $\cli$ contains $\cli_0:=^{\mbox{def}} \pi(\tilde{\mbox{UHF}}_d \otimes \mbox{UHF}_d)'' \vee \{U_z: z \in H \}''$. 
By the last proposition the centre of $\cli$, which is equal to $\{U_z: z \in H \}''$, contains the centre of $\pi(\tilde{\mbox{UHF}}_d \otimes \mbox{UHF}_d)''$ 
and thus by taking the commutant we also have $\cli \subseteq \pi(\tilde{\mbox{UHF}}_d \otimes \mbox{UHF}_d)''\vee \pi(\tilde{\mbox{UHF}}_d \otimes 
\mbox{UHF}_d)'$. In the last proposition we have described explicitly the factor decomposition of the representation 
$\pi$ of $\pi(\mbox{UHF}_d \otimes \tilde{\mbox{UHF}}_d)''$. One central issue is when such an factor decomposition 
is also an extremal decomposition. A clear answer at this stage seems to be somewhat hard. However the following 
proposition makes an attempt for our purpose. To that end we set few more notations and elementary properties. 

\vsp
For each $k \in \hat{H}$, let $\pi'_k$ be the representation $\pi$ of $\tilde{\mbox{UHF}}_d \otimes \mbox{UHF}_d$ restricted to $F_k$. We claim that each $\pi'_k$ 
is pure if $\pi'_0$ is pure. Fix any $k \in \hat{H}$ and let $X$ be an element in the commutant of $\pi'_k(\tilde{\mbox{UHF}}_d \otimes \mbox{UHF}_d)''$ then 
$S^*_IF_kS_I = S^*_I\Lambda^k(F_0)S_I = F_0$ as $S^*_JS_I$ commutes with $F_0$ for $|I|=|J|=k$ and further for any $|I|=|J|$, $S_I^*XS_IS_J^*S_J=S_I^*S_IS_J^*XS_J=S^*_JXS_J$ 
as $X$ commutes with $S_IS_J^*$ with $|I|=|J|$. Thus by Proposition 3.2 (c) $S^*_IXS_I$ is an element in commutant of $\pi'_0(\tilde{\mbox{UHF}}_d \otimes \mbox{UHF}_d)''$ 
for any $|I|=k$ and thus $S_I^*XS_I= c F_0$ for some scalar $c$ independent of $|I|=k$ as $\pi'_0$ is pure. We use 
the commuting property of $X$ with $\pi(\tilde{\mbox{UHF}}_d \otimes \mbox{UHF}_d)''$ to conclude that 
$X = c \Lambda^k(E_0)$ for some scalar $c$. If $k \le -1$ we employ the same method but with the endomorphism 
$\Lambda^{-k}$ so that $\Lambda^{-k}(X)$ is an element in the commutant of 
$\pi'_0(\tilde{\mbox{UHF}}_d \otimes \mbox{UHF}_d)''$. Thus $\sum_{I:|I|=-k }S_IXS^*_I = c I$ and by 
the injective property of the endomorphism we get $X$ is a scalar. Thus we conclude that each $\pi'_k$ is 
pure if $\pi'_0$ is pure. 

\vsp 
Next we claim that for each fixed $k \in \hat{H_0}$, the representation $\pi_k$ of $\tilde{\mbox{UHF}}_d \otimes \mbox{UHF}_d$ defined in the proof of Proposition 3.3 is quasi-equivalent to representation $\pi'_k$ ( here we recall 
$\hat{H_0} \subseteq \hat{H}$ as $H_0 \subseteq H$ ). That $\pi'_0$ is quasi-equivalent to 
$\pi_0$ follows as they are isomorphic by construction. In the proof of Proposition 3.3 we defined representation $\pi_k$ of $\tilde{\mbox{UHF}}_d \otimes \mbox{UHF}_d)'$ for all $k \in \hat{H}_0$ associated with minimal projections 
$\{E_k: k \in \hat{H}_0 \}$. More generally for any $k \in \hat{H}$, we denote $\pi_k$ for the restriction of $\pi$ to the minimal central projections 
$E_k$ on the subspace span by $\{\pi(\tilde{\mbox{UHF}}_d \otimes \mbox{UHF}_d)'f: \forall \;\; f\;\; \in \clh \otimes_{\clk} \tilde{\clh}, \;\;F_kf=f \}$. 
So each $E_k$ is a minimal central element containing $F_k$. However two such elements i.e. 
$E_k$ and $E_j$ are either equal or mutually orthogonal being minimal. Thus $\{E_k: k \in \hat{H} \} = \{E_k: k \in \hat{H}_0 \}$
and quasi-equivalence follows as $\pi_k$ is isomorphic with $\pi'_k$ for all $k \in \hat{H}_0$.     

\vsp 
We now set 
$$F'_0=[\pi(\tilde{\mbox{UHF}}_d \otimes \mbox{UHF}_d)''\Omega]$$ 
For a vector $f$ if $F'_0f=f$ then $U_zf=f$ for all $z \in H$ and thus $F'_0 \le F_0$. We prove in following text that equality holds if $\omega$ is pure.
 
\vsp 
First we consider the case when $H=\{z: z^n=1\}$. Projections $\Lambda(F'_0)$ and $\tilde{\Lambda}(F'_0)$ are elements in $\pi(\tilde{\mbox{UHF}}_d \otimes \mbox{UHF}_d)'$ by Proposition 3.3. The representation $\pi(\tilde{\mbox{UHF}}_d \otimes \mbox{UHF}_d)''$ restricted to both the projections $\Lambda(F'_0),\tilde{\Lambda}(F'_0)$ are pure as well. A pull back by the map $X \raro S^*_iXS_i$ with any $1 \le i \le d$ will do the job for the projection $\Lambda(F'_0)$. Thus 
$\Lambda(F'_0)\tilde{\Lambda}(F'_0)\Lambda(F'_0) = c \Lambda(F'_0)$ for 
some scalar. By pulling back with the action $X \raro S^*_i X S_i$ we get
$F'_0 S^*_i\tilde{\Lambda}(F'_0)S_i F'_0 = cF'_0$ and so 
$$c=<\Omega,S^*_i\tilde{\Lambda}(F'_0)S_i\Omega>$$
$$=\sum_k <\Omega,S_k S^*_iF'_0S_iS^*_k\Omega>$$
as $\tilde{S}^*_k\Omega=S^*_k\Omega$ and further $F'_0$ commutes 
with $\pi(\mbox{UHF}_d)$ and thus 
$$c=\sum_k <\Omega, S_k S_k^*\Omega>=1$$
This shows that $\tilde{\Lambda}(F'_0) \ge \Lambda(F'_0)$. Interchanging the role
of $\Lambda$ and $\tilde{\Lambda}$ we conclude that $\Lambda(F'_0)=\tilde{\Lambda}(F'_0)$. Now it essentially follows along the same line  
$\Lambda^k(F'_0)=\tilde{\Lambda}^k(F'_0)$ for all $k \ge 1$. 
By Proposition 2.5 we also note that 
$$\Lambda^n(F'_0)=F'_0=\tilde{\Lambda}^n(F_0)$$ 
as $H=\{z : z^n =1 \}$. Thus $F'=\sum_{0 \le k \le n-1}\Lambda(F'_0)$ is a $\Lambda$ and as well $\tilde{\Lambda}$ invariant projection. Since 
$F'\Omega=\Omega$ we conclude by the cyclic property of $\Omega$ for $\pi(\clo_d \otimes \tilde{\clo}_d)''$ that $F'=1$. Since $\Lambda^k(F'_0) 
\le F_k$ and $\sum_k F_k=1 $, we conclude that $\Lambda^k(F_0)=F_k$. In such a case we may check that  
$$F_k = [\pi(\tilde{\mbox{UHF}}_d \otimes \mbox{UHF}_d)''S^*_I\Omega: |I|=n-k]$$ 
for $1 \le k \le n-1$. 

\vsp
Similarly in case $H=S^1$ and $\omega$ is pure we also have $F_0=F'_0$ and 
for $k \ge 1$ 
$$F_k=[\pi(\tilde{\mbox{UHF}}_d \otimes \mbox{UHF}_d)''S_I\Omega: |I|= k]$$
$$F_{-k}=[\pi(\tilde{\mbox{UHF}}_d \otimes \mbox{UHF}_d)''S^*_I\Omega: |I|=k]$$
Thus we have got an explicit description of the complete atomic centre of $\cli$ 
when $\omega$ is a pure state.

\begin{pro}
Let $\omega,\psi$ and Popescu system $(\clk,\clm,v_k,\Omega)$ be as in Proposition 3.3. Then 

\NI (a) $\{\beta_z: z \in H \}$ invariant elements in $\pi(\tilde{\mbox{UHF}}_d \otimes \clo_d)''$ ( as well as in 
$\pi(\tilde{\clo}_d \otimes \mbox{UHF}_d)''$ ) are equal to $\pi(\mbox{UHF}_d \otimes \tilde{\mbox{UHF}}_d)''$. 

\NI (b) $\cli=\cli_0$ if and only if $\omega$ is pure. 

Further the following statements are equivalent: 

\NI (c) $\cli = \pi(\tilde{\mbox{UHF}}_d \otimes \mbox{UHF}_d)''$; 

\NI (d) $\pi(\tilde{\mbox{UHF}}_d \otimes \clo_d)''= \clb(\tilde{\clh} \otimes_{\clk} \clh)$;

\NI (e) $\pi(\tilde{\clo}_d \otimes \mbox{UHF}_d )''= \clb(\tilde{\clh} \otimes_{\clk} \clh)$;

\vsp
In such a case (i.e. if any of (c),(d) and (e) is true ) the following statements are also true: 

\NI (f) $\pi(\tilde{\mbox{UHF}}_d \otimes \mbox{UHF}_d)''$ is a type-I von-Neumann algebra with centre equal 
to $\{ U_z: z \in H \}''$ where $U_z$ is defined in Proposition 3.2. 

\NI (g) $\omega$ is a pure state on $\clb$. 

\vsp
Conversely if $\omega$ is a pure state then $\pi(\tilde{\mbox{UHF}}_d \otimes \mbox{UHF}_d)''$ is a type-I von-Neumann 
algebra with centre equal to $\{U_z: z \in H_0 \}''$ where $H_0$ is a subgroup of $H$.     
\end{pro} 

\vsp
\begin{proof} 
Along the same line of the proof of Proposition 2.5 (b) we get $\{\beta_z: z \in H \}$ invariant elements 
in $\pi(\clo_d \otimes \tilde{\mbox{UHF}}_d)''$ is $\pi(\tilde{\mbox{UHF}}_d \otimes \mbox{UHF}_d)''$ where factor property 
of $\pi(\clo_d)''$ is crucial as in proof of Proposition 2.5 (b). Same holds for $\pi(\mbox{UHF}_d \otimes \tilde{\clo}_d)''$ as 
$\pi(\tilde{\clo}_d)''$ is a factor. Here we comment that factor property of $\pi(\tilde{\clo}_d)''$ can be
ensured whenever $\psi$ is an extremal element in $K_{\omega'}$ (See Proposition 3.2 (a) ). 

\vsp
For (b) we will first prove $\cli_0=\cli$ if $\omega$ is pure. As by definition $\cli_0 \subseteq \cli$, it is enough if we show 
$\cli'_0 \subseteq \cli'$. Let $X \in \cli'_0$ i.e. $X$ commutes with $\{ U_z:z \in H\}''$ and $\pi(\mbox{UHF}_d \otimes 
\tilde{\mbox{UHF}}_d)''$. For each $k \in \hat{H}$, $F_kXF_k$ is an element in the commutant of $F_k\pi(\mbox{UHF}_d \otimes 
\tilde{\mbox{UHF}}_d)''F_k$. $\omega$ being pure each representation $\pi$ restricted to $F_k$ is irreducible and thus
$F_kXF_k = c_k F_k$ for some scalars $c_k$. Hence $X=\sum_k c_k F_k \in \cli'=\{U_z: z \in H \}''$. 

\vsp
For the converse we need to show that the restriction of $\pi(\mbox{UHF}_d \otimes 
\tilde{\mbox{UHF}}_d)''$ to $F'_0$ is pure. Let $X$ be an element on the subspace $F'_0$ and in the commutant 
of $F'_0\pi(\mbox{UHF}_d \otimes \tilde{\mbox{UHF}}_d)''F'_0$, ( which in our earlier notation $E'$ in Proposition 3.3 ). 
Then $X$ commutes with each $F_k$ for $k \in \hat{H}$ and $\pi(\mbox{UHF}_d \otimes 
\tilde{\mbox{UHF}}_d)''F_k$ as $F'_0 \le F_0$ and $F_k$ are orthogonal to 
$F'_0$ for $k \neq 0$. So $X$ commutes with $\{U_z: z \in H \}''$ and 
$\pi(\mbox{UHF}_d \otimes \tilde{\mbox{UHF}}_d)''$ i.e. $X \in \cli'_0$. By our assumption $\cli_0=\cli$, we have now
$X \in \cli'$ which is equal to $\{F_k: k \in \hat{H} \}''$ and so $X=cF_0$ 
for some scalar $c_0$. This shows that $F'_0=F_0$ and $\omega$ is pure.    

\vsp
(c) implies (d): $\{U_z: z \in H \}$ is a commuting family of unitaries such that $\beta_z(X)=U_zXU_z^*$ and 
thus by (c) $\{U_z: z \in H \}'' \subseteq \pi(\tilde{\mbox{UHF}}_d \otimes \mbox{UHF}_d)''$. Let $X$ be an 
element in the commutant of $\pi(\tilde{\mbox{UHF}}_d \otimes \clo_d)''$. Then 
$X$ commutes also with $\{U_z: z \in H \}''$ and thus $X \in \pi(\tilde{\mbox{UHF}}_d \otimes \mbox{UHF}_d)''$ by (c).
Hence $X$ is an element in the centre of $\pi(\tilde{\mbox{UHF}}_d \otimes \mbox{UHF}_d)''$ and so 
$X=\sum_k c_kE_k$ where $E_k$ are the minimal projections in the centre of $\pi(\tilde{\mbox{UHF}}_d \otimes \mbox{UHF}_d)''$ 
given in Proposition 3.3. However $X$ also 
commutes with $\pi(\clo_d)''$ by our assumption (c) and $\Lambda(E_k)=E_{k+1}$ for $k \in \hat{H}$.  
So $c_k=c_{k+1}$ and $X$ is a scalar multiple of unit operator. Hence (d) follows from (c). Along the same line we prove 
(c) implies (e). For a proof for (d) implies (c) and (e) implies (c), we simply apply (a).   

\vsp
Now we will prove (f) and (g). That $\pi(\tilde{\mbox{UHF}}_d \otimes \mbox{UHF}_d)''$ is 
a type-I von-Neumann algebra ( with completely atomic centre ) follows by a theorem of [So] once we use (c). In the proof 
of Proposition 3.3 we have proved that the centre of $\pi(\tilde{\mbox{UHF}}_d \otimes \mbox{UHF}_d)''$ is 
$\{U_z: z \in H_0 \}''$ where $H_0 \subseteq H$. For equality in the present situation we simply use (c), as $\beta_w(U_z)=U_z$
for all $w,z \in H$, to conclude that $U_z$ is in the centre of $\pi(\tilde{\mbox{UHF}}_d \otimes \mbox{UHF}_d)''$. 

\vsp
If (c) holds then $\cli_0=\cli$ and thus (g) follows by (b). Here we will give another proof using the same idea to prove (f). 
Let $X$ be an element in the commutant of $\pi_0(\tilde{\mbox{UHF}}_d \otimes \mbox{UHF}_d)''$, where $\pi_0$ is the factor 
representation on the minimal central projection $E_0$ defined in Proposition 3.3. Then $X$ commutes with 
$\{U_z: z \in H \}''$ and so by (c) $X$ in an element in $\pi(\tilde{\mbox{UHF}}_d \otimes \mbox{UHF}_d)''$. 
So $X$ is in the centre of $\pi_0(\tilde{\mbox{UHF}}_d \otimes \mbox{UHF}_d)''$. $\pi_0$ being a factor representation
$X$ is a scalar multiple of $E_0$. Thus $\pi_0$ is an irreducible representation and so $\omega$ is pure. 

\vsp
By Proposition 3.1 we recall that $\pi'_0$ is unitary equivalent to the GNS representation of $(\clb,\omega)$. 
Thus $\pi'_0$ is irreducible if and only if $\omega$ is pure. So for a pure state $\omega$, 
for each $k \in \hat{H}_0$, $\pi_k$ being quasi-equivalent to $\pi'_k$, $\pi_k$ is a type-I factor 
representation of $\pi(\mbox{UHF}_d \otimes \tilde{\mbox{UHF}}_d)''$. This completes the proof. 
\end{proof} 
  
\vsp
The following theorem is the central step that will be used repeatedly.

\vsp
\begin{pro} 
Let $\omega$ be an extremal translation invariant state on $\clb$ and $\psi$ be an extremal 
element $\psi$ in $K_{\omega}$. We consider the Popescu elements $(\clk,v_k:1 \le k \le d, \clm,\Omega)$
as in Proposition 2.4 for the dual Popescu elements and associated amalgamated representation $\pi$ of
$\clo_d \otimes \tilde{\clo}_d$ as described in Proposition 3.1. Let $\cle$ and $\tilde{\cle}$ be the support 
projections of the state $\psi$ in $\pi(\clo_d)''$ and $\pi(\tilde{\clo}_d)''$ respectively. Let $\clf$ and $\tilde{\clf}$ be the projections 
$[\pi(\clo_d)''\Omega]$ and $[\pi(\tilde{\clo}_d)\Omega]$ respectively. Then the following holds:

\vsp 
\NI (a) $\pi(\tilde{\clo}_d \otimes \clo_d)''= \clb( \tilde{\clh} \otimes_{\clk} \clh)$;

\vsp 
\NI (b) $\pi(\tilde{\clo}_d)''=\pi(\clo_d)'$ if and only if $\pi(\tilde{\clo}_d)''\cle=\pi(\clo_d)'\cle$; 

\vsp 
\NI (c) $Q=\cle \tilde{\cle}$ is the support projection of the state $\psi$ in $\pi(\clo_d)''\tilde{\cle}$
and also in $\pi(\tilde{\clo}_d)''\cle$. Further $P=\cle \clf \le Q$;

\vsp 
\NI (d) If $\cle \clf = \tilde{\cle}\tilde{\clf}$ then $\cle=\tilde{\clf},\;\tilde{\cle}=\clf,\;P=Q$;

\vsp 
\NI (e) If $P=Q$ then the following statements are true:

\vsp 
\NI (i) $\clm'= \tilde{\clm}$ where $\tilde{\clm} = \{ P\tilde{S}_iP: 1 \le i \le d
\}''$;

\vsp 
\NI (ii) $\pi(\clo_d)' = \pi(\tilde{\clo}_d)''$;

\vsp 
\NI (f) If $P=[\tilde{\clm}\Omega]$ then $\clm'=\tilde{\clm}$; 

\vsp 
\NI (g) $\omega$ is pure on $\clb$ if and only if there exists a sequence of elements $x_n \in \clm$ such that for each $m \ge 0$ 
$x_{n+m}\tau_n(x) \raro \phi(x)1$ as $n \raro \infty$ in strong operator topology, equivalently $\phi(\tau_n(x)x_{n+m}^*x_{n+m}\tau_n(y)) \raro \phi(x)\phi(y)$ as $n \raro \infty$ for 
all $x,y \in \clm$ where $\clm = \{ v_i=PS_iP: 1 \le i \le d \}''$ and $\tau(x)=\sum_{1 \le k \le d}v_kxv_k^*,\;\;x \in \clm$; Same holds true if we replace $\clm_0$ for 
$\clm$ where $\clm_0=\{x \in \clm: \beta_z(x)=x; z \in H \}$.  

\end{pro} 

\vsp
\begin{proof} 
(a) is a restatement of Proposition 3.2 (a). $\cle$ ( $\tilde{\cle}$ ) being the support projection of
the state $\psi$ in $\pi(\clo_d)''$ ( $\pi(\tilde{\clo}_d)''$ ) and $\psi = \psi \Lambda$ we have $\Lambda(\cle) \ge \cle$ 
and further we have $\cle=[\pi(\clo_d)'\Omega] \ge [\pi(\tilde{\clo}_d)''\Omega]$ and hence increasing projections in $\pi(\clo_d)''$, 
$\Lambda^n(\cle) \uparrow I$ as $n \raro \infty$ because $\Omega$ is cyclic for $\pi(\clo_d \otimes \tilde{\clo}_d)''$ in $\clh \otimes_{\clk} \tilde{\clh}$ 
and limiting projection being in the commutant of $\pi(\clo_d)''$ as well i.e. since $\Lambda(Y)=Y$ implies $Y \in \pi(\clo_d)''$ by Cuntz relations.   

\vsp 
We set von-Neumann algebras $\cln_1=\pi(\clo_d)'\cle$ and $\cln_2=\pi(\tilde{\clo}_d)''\cle$. By our construction we have $\pi(\tilde{\clo}_d)'' 
\subseteq \pi(\clo_d)' $ and so $\cln_2 \subseteq \cln_1$. Since $\Lambda^n(\cle) \uparrow I$ as $n \raro \infty$ in strong operator topology, 
two operators in $\pi(\clo_d)'$ are same if their actions are same on $\cle$. So (b) is true. 

\vsp 
For (c) we note that $Q=\cle\tilde{\cle} \in \cln_2 \subseteq \cln_1$ and claim that $Q$ is the support projection of the state $\psi$ in $\cln_2$. 
To that end let $x\cle \ge 0$ for some $x \in \pi(\tilde{\clo}_d)''$ so that $\psi(QxQ)=0$. As $\Lambda^k(x\cle) \ge 0$ for all $k \ge 1$ and $\Lambda^k(\cle) \raro I$ 
we conclude that $x \ge 0$. As $\cle \Omega=\Omega$ and thus $\psi(\tilde{\cle}x\tilde{\cle})=\psi(QxQ)=0$, we conclude
$\tilde{\cle}x\tilde{\cle}=0$, $\tilde{\cle}$ being the support projection for $\pi(\tilde{\clo}_d)''$. Hence $QxQ=0$.
As $\psi(Q)=1$, we complete the proof of the claim that $Q$ is the support of $\psi$ in $\cln_2$. Similarly $Q$ is also
the support projection of the state $\psi$ in $\pi(\clo_d)''\tilde{\cle}$. This completes the proof of (c).

\vsp 
Thus if $\cle \clf=\tilde{\cle}\tilde{\clf}$, we get $\Lambda^n(\cle)\clf=\tilde{\cle} \Lambda^n(\tilde{\clf})$ and 
$\cle \tilde{\Lambda}(\clf) = \tilde{\Lambda}(\tilde{\cle})\tilde{\clf}$ and thus taking limit we get 
$\clf=\tilde{\cle}$ and $\cle=\tilde{\clf}$. It is obvious now that $P=\cle \clf = \cle \tilde{\cle}=Q$. This completes 
the proof of (d). 

\vsp
As $\cle \in \pi(\clo_d)''$ and $\tilde{\cle} \in \pi(\tilde{\clo}_d)''$ we check that von-Neumann algebras 
$\clm^1 = Q \pi(\clo_d)'' Q$ and $\tilde{\clm}^1 = Q \pi(\tilde{\clo}_d) Q$ acting on $Q$ satisfy 
$\tilde{\clm}^1 \subseteq \clm^{1'}$. Now we explore that $\pi(\tilde{\clo}_d \otimes \clo_d)''= 
\clb(\clh \otimes_{\clk} \tilde{\clh})$ and note that in such a case $Q\pi( \tilde{\clo}_d \otimes \clo_d )''Q$ 
is the set of all bounded operators on the Hilbert subspace $Q$. As $\cle \in \pi(\clo_d)''$ and 
$\tilde{\cle} \in \pi(\tilde{\clo}_d)''$ we check that together $\clm^1 = Q \pi(\clo_d)'' Q$ and 
$\tilde{\clm}^1 = Q \pi(\tilde{\clo}_d) Q$ generate all bounded operators on $Q$. Thus both $\clm^1$ 
and $\tilde{\clm}^1$ are factors. The canonical states $\psi$ on $\clm^1$ and $\tilde{\clm}^1$ are 
faithful and normal. We set $l_k=QS_kQ$ and $\tilde{l}_k=Q\tilde{S}_kQ,\; 1 \le k \le d$ and recall 
that $v_k=PS_kP$ and $\tilde{v}_k=P\tilde{S}_kP,\; 1 \le k \le d$. We note that $Pl_kP=v_k$ and 
$P\tilde{l}_kP = \tilde{v}_k$ where we recall that, by our construction, $P$ is the support projection of 
the state $\psi$ in $\pi(\clo_d)''_{|}[\pi(\clo_d)\Omega]$. $Q$ being the support projection of 
$\pi(\clo_d)\tilde{\cle}$, by Theorem 2.4 applied to Cuntz elements $\{S_i\tilde{E}: 1 \le i \le d \}$, 
$\tilde{\cle} \pi(\clo_d)'\tilde{\cle}$ is order isomorphic to $\clm^{1'}$ via the map $X \raro QXQ$. As the projection 
$\clf=[\pi(\clo_d)''\Omega] \in \pi(\clo_d)'$, we check that the element $Q\clf\tilde{\cle}Q \in \clm^{1'}$. However 
$Q\clf \tilde{\cle}Q=\cle\tilde{\cle}\clf\tilde{\cle}\cle=QPQ=P$ and thus $P \in \clm^{1'}$. We also check that 
$\clm^1\Omega=\clm^1P\Omega=P\clm^1\Omega=\clm \Omega$ and thus $P=[\clm^1\Omega]$. We set 
$\tilde{\clm}$ for the von-Neumann algebra generated by $\{\tilde{v}_k:\; 1 \le k \le d \}$. 

\vsp 
So far our analysis did not use the hypothesis on statement (e) i.e. $P=Q$. For $P=Q$, we have 
$\clm^1=\clm$ and $\tilde{\clm}^1=\tilde{\clm}$. By order isomorphic property 
we get (i) is equivalent to $\tilde{\cle} \pi(\clo_d)' \tilde{\cle} = \tilde{\cle} \pi(\tilde{\clo}_d)''\tilde{\cle}$ 
and taking commutant again we get $\pi(\clo_d)''\tilde{\cle} = \pi(\tilde{\clo}_d)'\tilde{\cle}$. Now we invoke the first 
part of the argument changing the role or using the endomorphism $\tilde{\Lambda}$ we conclude that 
$\pi(\clo_d)''=\pi(\tilde{\clo}_d)'$. This completes the proof of (e) provided we find a independent proof for (i) 
which is not so evident and this crucial point was not noticed in the proof given for Lemma 7.8 in [BJKW].  

\vsp 
Now we will analyze the representation $\pi$ of $\tilde{\clo}_d \otimes \clo_d$ which is pure to prove (i). To that end we note 
since $P=Q$ by our assumption, $\Omega$ is a common cyclic and separating vector for $\tilde{\clm}$ and $\clm'$.   
Thus we can get an endomorphism $\alpha: \clm' \raro \tilde{\clm}$ defined by 
$$\alpha(y)=\tilde{\clj}\clj y \clj \tilde{\clj}$$ 
where $\tilde{\clj}$ is the Tomita's conjugate operator associated with the cyclic and separating vector $\Omega$ for $\tilde{\clm}$  
( i.e. for $y \in \clm'$, we have $\clj y \clj \in \clm \subseteq \tilde{\clm}'$ since $\tilde{\clm} \subseteq \clm'$. 
Since $\tilde{\clj}\tilde{\clm}'\tilde{\clj}=\tilde{\clm}$, we have $\alpha(y) \in \tilde{\clm}$ ). We note that the general theory 
does not guarantee [AcC] that the endomorphism be Takesaki's canonical conditional expectation associated with $\phi$. If so then the modular 
automorphism group $(\sigma_t)$ of $\clm'$ also preserves $\tilde{\clm}$. Thus $\sigma_z(x) \in \tilde{\clm}$ for $-1 \le Im(z) \le 1$ 
if $x$ is an analytic element in $\tilde{\clm}$. Thus we would have got $\clj v_k(\delta) \clj = 
\sigma_{i \over 2}(\tilde{v}_k(\delta)) \in \tilde{\clm}$ where $x(\delta)$ is average of 
$\sigma_t(x)$ with respect to Gaussian measure with variance $\delta > 0$. That $\tilde{v}_k(\delta)$ 
is an analytic element follows from the general Tomita-Takesaki theory [BR1]. Since $v_k(\delta) \raro v_k$ in strong 
operator topology as $\delta \raro 0$ and $\clm = \{v_k: 1 \le k \le d \}''$ together with 
$\clj \clm \clj = \clm'$ we arrived at $\tilde{\clm}=\clm'$. In the following we avoid this tempting 
route and aim to explore the general representation theory of $C^*$-algebras [BR1,chapter 2].   

\vsp 
We claim that $\clm'=\tilde{\clm}$. Suppose not i.e. $\tilde{\clm} \subset \clm'$. Then $\alpha(\tilde{\clm})$ is a proper von-Neumann subalgebras 
of $\alpha(\clm') \subseteq \tilde{\clm} $ being an into map and hence $\alpha(\tilde{\clm})$ is a proper von-Neumann 
sub-algebra of $\tilde{\clm}$. Now consider the Popescu elements $(\clk,\alpha(\tilde{v}_i),\Omega)$ and its dilation 
as in Theorem 2.1. Then by the commutant lifting theorem applied to pairs 
$(\tilde{v}_i),\alpha(\tilde{v}_i)$ we find a unitary operator $U$ on $\tilde{\clh}$ so that 
$U\pi(\tilde{\clo}_d)''U^*$ is strictly contained in $\pi(\tilde{\clo}_d)''$ ( Without loss of generality we can 
take the dilated Hilbert space for $(\clk,\alpha(\tilde{v}_i),\Omega)$ to be same as $\tilde{\clh}$ as there exists an 
isomorphism preserving $\clk$, see the remark that follows after Theorem 2.1 ). We extend $U$ to a unitary operator 
on $\tilde{\clh} \otimes_{\clk} \clh$ and denote $\pi_u(x)=U \pi(x) U^*$ for $x \in \tilde{\clo}_d \otimes \clo_d$ which 
is unitary equivalent to the pure representation $\pi$ and $\pi_u(\tilde{\clo}_d)''$ is strictly contained in 
$\pi(\tilde{\clo}_d)''$. Now $\pi_u$ is also an amalgamated representation over the subspace $\clk$ with 
$P_u=Q_u$. Thus we can repeat now same procedure with $\pi_u$ and so on. Note that the process won't terminate in finite 
time. Our aim is to find a contradiction from this using formal set theory.  

\vsp 
To that end we use temporary notation $\pi_0$ for $\pi$ defined in last paragraph and $\pi$ will be used for a generic representation. 
Let $\clp$ be the collection of representation $(\pi,H_{\pi},\Omega)$ quasi-equivalent to $\pi_0:\tilde{\clo}_d \otimes \clo_d  
\raro \clb(\tilde{\clh} \otimes_{\clk} \clh )$ with a shift invariant vector state $\omega(x)=<\Omega,\pi(x)\Omega>$ i.e. 
$\omega(\pi(\theta(x))=\omega(\pi(x))$. So there exists cardinal numbers $n_{\pi},n_0(\pi)$ so that $n_{\pi}H_{\pi}$ is unitary 
equivalent to $n_0(\pi)\pi_0$. Thus given an element $(\pi,H_{\pi},\Omega)$ we can associate two cardinal numbers $n_{\pi}$ and $n_0(\pi)$ 
and without loss of generality we assume that $H_{\pi} \subseteq n_0(\pi)H_0$ and $n_{\pi}H_{\pi}=n_0(\pi)H_0.$  $\pi_0$ being a pure representation, 
any element $\pi \in \clp$ is a type-I factor representation of $\tilde{\clo}_d \otimes \clo_d $. The interesting point here that $\dsp{\oplus_{\pi \in \clp }}\pi$ 
is also an element in $\clp$ with associated cardinal numbers $\sum_{\pi} n_{\pi}$ and $\sum_{\pi} n_0(\pi)$. We say $(\pi_1,H_{\pi_1},\Omega^1) \prec (\pi_2,H_{\pi_2},\Omega^2)$ 
if there exists an isometry $U:n_{\pi_1} H_{\pi_1} \raro n_{\pi_2} H_{\pi_2}$ so that 

\vsp 
\NI (C1) For each $1 \le \alpha \le n_{\pi_1} $ we have 
$U \Omega^1_{\alpha} = \Omega^2_{\alpha'}$ for some 
$1 \le \alpha' \le n_{\pi_2}$;

\NI (C2) $n_{\pi_2}\pi_2(x)E'_2 = U n_{\pi_1}\pi_1(x) U^* $ where $\cle'_2 \in n_{\pi_2}\pi_2(\tilde{\clo}_d)'$; 
 
\NI (C3) $U\oplus_{1 \le \alpha \le n_{\pi_1} }\pi^{\alpha}_1(\tilde{\clo}_d)''U^* \subset \oplus_{1 \le 
\alpha \le n_{\pi_2}}\pi^{\alpha}_1(\tilde{\clo}_d)''E'_2$. 

\vsp 
That the partial order is non-reflexive follows as $(\pi,H_{\pi},\Omega) \prec (\pi,H_{\pi},\Omega)$ contradicts (C3) as $I=E_2'$. By our starting assumption that $\clm' \neq \tilde{\clm}$ 
we check that $\pi_0 \prec \pi_u$. Thus going via the isomorphism we also check that for a given element $\pi \in \clp$ there exists an element $\pi' \in \clp$ so that $\pi \prec \pi'$. 
Thus $\clp_0$ is a non empty set and has at least one infinite chain. Partial order property follows easily. If $\pi_1 \prec \pi_2$ and $\pi_2 \prec \pi_3$ then $\pi_1 \prec \pi_3$. 
If $U_{12}$ and $U_{23}$ are isometric operators that satisfies (C1)-(C3) respectively, then $U_{13}=U_{23}U_{12}$ will do the job for $\pi_1$ and $\pi_3$. 
 
\vsp 
However by Hausdorff maximality theorem there exists a non-empty maximal totally ordered subset $\clp_0$ of $\clp$. We claim that 
$\pi_{max}=\oplus_{\pi \in \clp_0} \pi$ on $H_{\pi_{max}}=\oplus_{\pi \in \clp_0} H_{\pi} $ is an upper bound 
in $\clp_0$. That $\pi_{max} \in \clp$ is obvious. Further given an element $(H_1,\pi_1,\Omega_1) \in \clp_0$ there exists 
an element $(H_2,\pi_2,\Omega_2) \in \clp_0$ so that $\pi_1 \prec \pi_2$ by our starting remark as $\pi_0 \prec \pi_u$. 
By extending isometry $U_{12}$ to an isometry from $H_1 \raro n_{\pi_{max}} H_{\pi_{max}}$ trivially we get the 
required isometry that satisfies (C1),(C2) and (C3) where cardinal numbers $n_{\pi_{max}} = \sum_{\pi \in \clp_0} 
n_{\pi} \in \aleph_0$. Thus by maximal property of $\clp_0$ we have $\pi_{max} \in \clp_0$. This brings a contradiction 
as by our construction $(\pi_{max},H_{\pi_{max}},\Omega) \prec (\pi_{max},H_{\pi_{max}},\Omega)$ 
as $\pi_{max} \in \clp_0$ but partial order is strict. This contradicts our starting hypothesis that $\tilde{\clm}$ is 
a proper subset of $\clm'$. This completes the proof for (i) of (e) $\clm'=\tilde{\clm}$ when $P=Q$.  

\vsp 
In the proof of $\clm'=\tilde{\clm}$ in (e), we have used equality $P=Q$ just to ensure that $\Omega$ is also a cyclic for $\tilde{\clm}$ and $P=Q$ is used 
to prove $\pi(\clo_d)' = \pi(\tilde{\clo}_d)''$. So (f) follows by the proof of (e). 

\vsp 
A proof for (g) is given in [Mo3] with $\clm_0$. Here we will also give an alternative proof relating the criteria obtained in Proposition 3.4. 
To that end we claim that $$\bigcap_{n \ge 1}\tilde{\Lambda}^n(\pi(\mbox{UHF}_d)')=\pi(\tilde{\mbox{UHF}}_d)'
\bigcap \pi(\mbox{UHF}_d)'.$$ That $\tilde{\Lambda}^n(\pi(\mbox{UHF}_d)') \subseteq \{ \tilde{S}_I\tilde{S}^*_J: |I|=|J| < \infty \}'$ follows by 
Cuntz relation and thus $\bigcap_{n \ge 1}\tilde{\Lambda}^n(\pi(\mbox{UHF}_d)') \subseteq \pi(\tilde{\mbox{UHF}}_d)'\bigcap 
\pi(\mbox{UHF}_d)'$. For the reverse inclusion let $X \in \pi(\tilde{\mbox{UHF}}_d)'\bigcap \pi(\mbox{UHF}_d)'$. For $n \ge 1$, 
we choose $|I|=n$ and set $Y_n=\tilde{S}^*_IX\tilde{S}_I$. 
We check that it is independent of the index that we have chosen as 
$Y_n=\tilde{S}^*_IX\tilde{S}_I\tilde{S}^*_J\tilde{S}_J=\tilde{S}^*_I\tilde{S}_I\tilde{S}^*_JX \tilde{S}_J=\tilde{S}^*_JX\tilde{S}_J$ 
where in second equality we have used $X \in \pi(\tilde{\mbox{UHF}}_d)'$ and also 
$\Lambda^n(Y_n)=\sum_{|J|=n } \tilde{S}_J\tilde{S}^*_IX\tilde{S}_I\tilde{S}^*_J=X$. This proves the equality in the claim. Going along the same line we also get 
$$\bigcap_{n \ge 1}\tilde{\Lambda}^n(\pi(\clo_d)') = \pi(\tilde{\mbox{UHF}}_d)'\bigcap \pi(\clo_d)'= \pi(\tilde{\mbox{UHF}}_d \otimes \clo_d)'.$$ 

\vsp 
By Proposition 3.4 $\omega$ is pure if and only if the set above is trivial. Thus once more by Proposition 1.1 in [Ar2] and 
Theorem 2.4 in [Mo2], purity is equivalent to asymptotic relation $||\psi \tilde{\tau}^n -\phi|| \raro 0$ as $n \raro \infty$ 
for any normal state $\psi$ on $\clm'$ ( Here we recall by Proposition 2.4 
$P\pi(\clo_d)'P=\clm'$ as $P$ is also the support projection 
in $\pi(\clo_d)''\clf$ and support projection of $\phi_{\Omega}$ in $\pi(\clo_d)'$ 
is $\clf=[\pi(\clo_d)''\Omega]$, where commutant is taken in $\clb(\clk)$ ). 
By a duality argument Theorem 2.4 in [Mo3] we conclude that $\omega$ is pure if and only if there exists 
a sequence of elements $x_n \in \clm$ so that for each $m \ge 0$, $x_{m+n}\tau_n(x) \raro \phi(x)1$ as $n \raro \infty$ for all $x \in \clm \subseteq \clb(\clk)$. 
This completes the proof of (d) with $\clm$. For the proof with $\clm_0$ we need to show if part as only if part follows $\clm_0$ being a subset of $\clm$ and $\tau$ 
takes elements of $\clm_0$ to itself. For if part we refer to Theorem 3.2 in [Mo3].
\end{proof} 

\vsp
We set 
$$(\clm')_0 = \{x \in \clm' : \beta_z(x)=x,\;\; z \in H \}.$$
Similarly we also set $\tilde{\clm}_0$ and $(\tilde{\clm}')_0$ as $(\beta_z:\;z \in H)$ invariant elements
of $\tilde{\clm}$ and $(\tilde{\clm}')$ respectively. We note that as a set $(\tilde{\clm}_0)'$ could be different from 
$(\tilde{\clm}')_0$. We note also that $P\tilde{\clm}^1P \subseteq \tilde{\clm}$ and unless $P$ is an element 
in $\tilde{\clm}^1$, equality is not guaranteed for a factor state $\omega$. The major problem is to show that 
$P$ is indeed an element in $\tilde{\clm}^1$ when $\omega$ is a pure state.   

\vsp
We warn here an attentive reader that in general for a factor state $\omega$, the set $\clf \pi(\tilde{\clo}_d)'' \clf$, which is a 
subset of $\clf \pi(\clo_d)'\clf$, need not be an algebra. However by commutant lifting theorem applied to dilation $v_i \raro S_i\clf$, 
$\pi(\clo_d)'\clf$ is order isomorphic to $\clm'$ as $P=\clf \cle$ is the support projection. Thus the von-Neumann sub-algebra generated 
by the elements $\clf\pi(\tilde{\clo}_d)''\clf$ is order isomorphic to $\tilde{\clm}$. However 
$\tilde{\clm}_0$ may properly include $\tilde{\clm}_{00}=\{ P\pi(\tilde{\mbox{UHF}}_d)P \}''$ 
( as an example take $\psi$ to be the unique KMS state on $\clo_d$ and $\omega$ be the unique trace on $\clb$ 
for which we get $\tilde{\clm}_{00} = \IC$ and $P\pi(\tilde{\clo}_d)''P$ is the linear span of 
$\{\tilde{v}_J^*,\;I,\;\tilde{v}_J:|J| < \infty \}$. 

\vsp 
Existence of a $\phi$ preserving normal conditional expectation $\int_{z \in H} \beta_z dz: \clm \raro \clm_0$ by Proposition 2.5 ensures that modular operator of $\phi$ preserves $\clm_0$ [Ta] and so does on 
$(\clm')_0$. However there is no reason to take it granted for $\tilde{\clm}_0$ to be 
invariant by the modular group of $((\clm')_0,\phi)$. By Takesaki's theorem such a property 
is true if and only if there exists a $\phi$-invariant norm one projection from $(\clm')_0$ 
onto $\tilde{\clm}_0$. In the following we avoid this tempted route.  

\vsp 
At this stage it is not clear how we can ensure existence of a normal conditional expectation from $\clm'$ onto 
$\tilde{\clm}$ directly and so the equality $\clm'=\tilde{\clm}$ when $\omega$ is a pure state. Further interesting point here that the equality $\tilde{\clm}=\clm'$ 
holds when $\omega$ is the unique trace on $\clb$ as $v^*_k=S^*_k$ and $\clj \tilde{v}^*_k \clj = {1 \over d} S_k$ for all $1 \le k \le d$ where $P \neq Q$ and $\pi(\clo_d)' 
\supset \pi(\tilde{\clo}_d)''$. In the last proposition we have also proved if $[\tilde{\clm}\Omega]=P$ then 
$\clm'=\tilde{\clm}$. Thus a natural question that arises here: how the equality $P=Q$ is related to 
purity of $\omega$?  We are now in a position to state the main mathematical result of this section. 

\vsp 
\begin{thm} Let $\omega$ be as in Theorem 3.5. Then the following holds:

\vsp 
\NI (a) $P$ is also the support projection of $\psi$ in $\pi(\tilde{\clo}_d))''_{|}\tilde{\clh}$ if and only if $\omega$ is pure. 

\vsp 
\NI (b) If $\omega$ is pure then the following holds:   

\NI (i) $\clm'= \tilde{\clm}$ where $\tilde{\clm} = \{ P\tilde{S}_iP: 1 \le i \le d
\}''$;

\NI (ii) $\pi(\clo_d)' = \pi(\tilde{\clo}_d)''$.

\NI (iii) $\pi_{\omega}(\clb_R)' = \pi_{\omega}(\clb_L)''$;

\end{thm}

\vsp
\begin{proof} 
First we will prove that $\omega$ is pure if $P$ is also the support projection of the state $\psi$ in $\pi(\tilde{\clo}_d)''\tilde{\clf}$, 
where $\tilde{\clf}=[\pi(\tilde{\clo}_d)''\Omega]$. The support projection of $\psi$ in $\pi(\tilde{\clo}_d)''\tilde{\clf}$ is $\tilde{\cle}\tilde{\clf}$ and thus we also 
have $P = \tilde{\cle}\tilde{\clf}$ by our hypothesis. Since $\Lambda^n(P)=\Lambda^n(\cle)\clf \uparrow \clf$ and now
$\Lambda^n(P)=\tilde{\cle}\Lambda^n(\tilde{\clf}) \uparrow \tilde{\cle}$ as $n \uparrow \infty$, we also have 
$\clf = \tilde{\cle}$. Similarly we also have for each $n$, $\cle\tilde{\Lambda}^n(\clf)=\tilde{\Lambda}^n(\tilde{\cle})\tilde{\clf}$ 
and thus taking limit we also get $\cle=\tilde{\clf}$.  

\vsp 
So we have $P=\cle \clf=\cle \tilde{\cle}=Q$. $\tilde{\clm}=P\pi(\tilde{\clo}_d)''P$ is cyclic in $\clk$ i.e. 
$[\tilde{\clm}\Omega]=[P\pi(\tilde{\clo}_d)''P\Omega]= P\tilde{\clf}=P\cle=P$ as $\tilde{\clf}=\cle$.  

\vsp 
However $\bigcap_{n \raro \infty} \tilde{\Lambda}^n(\pi(\tilde{\mbox{UHF}}_d)) = \pi(\tilde{\mbox{UHF}}_d)'' \bigcap \pi(\tilde{\mbox{UHF}}_d)'$ ( for a proof which is a simple 
application of Cuntz relation, we refer to section 5 of [Mo2]). Further $\psi$ being a factor state in $K_{\omega'}$, by Proposition 3.2 $\pi(\tilde{\mbox{UHF}}_d)''$ is a factor. 
In particular we have $\bigcap_{n \raro \infty} \tilde{\Lambda}^n(\pi(\tilde{\mbox{UHF}}_d))\tilde{\clf} = \IC \tilde{\clf}$. Thus by Proposition 1.1 in [Ar2] we conclude that $||\Psi \circ \tilde{\tau}_n -\phi|| \raro 0$ 
as $n \raro \infty$ for all normal state $\Psi$ on $\tilde{\clm}_0$ where $\tilde{\clm}_0=P \pi(\tilde{\mbox{UHF}}_d)''P$ as $\tilde{\clf}=\cle$ and support projection of $\psi$ in $\pi(\tilde{\mbox{UHF}}_d)''$ 
is $\tilde{\cle}$ and $P=\cle\tilde{\cle}\cle$. 

\vsp 
Note that $\tilde{\clm}_0 \subseteq \clm_0'$ where $\clm_0=P\pi(\mbox{UHF}_d)''P$. Further by Proposition 2.5 
$\clm_0 = \{x \in \clm:\beta_z(x)=x;z \in H \}\;\;\;$ and  $\tilde{\clm}_0=\{x \in \tilde{\clm}:
\beta_z(x)=x;z \in H \}$. Once we set $P_0=[\clm_0\Omega]$ then we also have $P_0=[\tilde{\clm}_0\Omega]$ 
as $[\tilde{\clm}\Omega]=P=[\clm\Omega]$ by expending $u_z=\sum_{k \in \hat{H}} z^kP_k$ where $z \raro u_z=
PU_zP$ is a unitary representation of group $H$.      

\vsp
For $x \in \tilde{\clm}, y \in \tilde{\clm}'$ we have  
$$\phi(\tilde{\tau}(x)y)=\sum_k <\tilde{v}^*_k\Omega,x\tilde{v}^*_ky\Omega> = \sum_k <v^*_k\Omega,x y \tilde{v}^*_k\Omega>$$  
(as $v_k^*\Omega=\tilde{v}^*_k\Omega$ )
$$=\sum_k <\Omega, x v_kyv_k^* \Omega>=\phi(x \tau(y))$$
The dual group of $(\tilde{\clm},\tilde{\tau},\phi)$ is given on the commutant by 
$(\tilde{\clm}',\tau,\phi)$ where $\tau(x)=\sum_k v_k x v_k^*$ for $x \in \tilde{\clm}'$. 
where commutant is taken in $\clb(\clk)$. 
Now moving to $\{\beta_z: z \in H \}$ invariant elements in the duality relation above, we verify that 
adjoint Markov map of $(\tilde{\clm}_0,\tilde{\tau},\phi)$ is given by 
$(\tilde{\clm}'_0,\tau,\phi)$ where $\tilde{\clm}'_0$, the commutant of $\tilde{\clm}_0$ 
is taken in $\clb(\clk_0)$ and $\clk_0$ is the Hilbert subspace $P_0$ with $\Omega$ as
cyclic and separating vector for $\tilde{\clm}_0$ in $\clk_0$. Thus by Theorem 2.4 in [Mo3], 
there exists a sequence of elements $y_n \in \tilde{\clm}'_0$ such that for each $m \ge 1$, 
$y_{m+n}\tau_n(y) \raro \phi(y)1$ as $n \raro \infty$ for all $y \in \tilde{\clm}'_0 \subseteq \clb(\clk_0)$. 
Thus $\omega$ is pure by Proposition 3.5 (e) once we recall $\clm'=\tilde{\clm}$ as $P=Q$ 
( and so $\clm'_0 = \tilde{\clm}_0$ ) by Proposition 3.5 (d). This completes the proof of
 purity property of $\omega$.

\vsp  
In the following we now prove the hardest part of the theorem namely $\tilde{\clf}=\cle$ and $\clf=\tilde{\cle}$ if $\omega$ is pure. Proof uses extensively 
the general theory of quasi-equivalent representation of a $C^*$ algebra and we refer [BR1, Chapter 2.4.4] as a general reference.   

\vsp 
We set unitary operator $V=\sum_k S_k\tilde{S}^*_k$. That $V$ is a unitary operator follows by Cuntz's relations and commuting property 
of $(S_i)$ and $(\tilde{S}_i)$. Further a simple computation shows that $V\pi(x)V^*=\pi(\theta(x))$ for all $x \in \clb = \clb_L \otimes \clb_R$, 
identified with $\tilde{\mbox{UHF}}_d \otimes \mbox{UHF}_d$ and $\theta$ is the right shift. We also have    
\be 
\theta(\cle)=V \cle V^*= \sum_{k,k'}S_k\tilde{S}^*_k \cle S^*_{k'}\tilde{S}_{k'}=\Lambda(\cle) \ge \cle
\ee 
Similarly we also have
\be 
\theta(\tilde{\clf}) \ge \tilde{\clf} 
\ee
So in particular we have $V(I-\cle)V^* \le I-\cle$, i.e. $(I-\cle)V^*\cle=0$.  
Also for any $X \in \pi(\clo_d)'$ we have $V^*\tilde{\clf}X\Omega = 
\tilde{\clf} \sum_k \tilde{S}_k X \tilde{S}_k^*\Omega$ as $S_k^*\Omega 
=\tilde{S}_k^*\Omega$. Thus $(I-\tilde{\clf})V^*\tilde{\clf}=0$ i.e. 
$V\tilde{\clf}V^*=\Lambda(\tilde{\clf}) \ge \tilde{\clf}$. Similarly we also 
have 
\be 
\theta^{-1}(\tilde{\cle})=V^*\tilde{\cle}V \ge \tilde{\cle}\;\mbox{and}\; V^* \clf V \ge \clf
\ee 
We also set two family of increasing projections for all natural numbers $n \in \IZ$ as follows 
\be 
\cle_n=V^n \cle (V^n)^*,\;\;\tilde{\clf}_n=V^n\tilde{\clf}(V^n)^*
\ee

\vsp 
Since $\beta_z(V)=V$ for all $z \in H$, $V \in \pi(\tilde{\mbox{UHF}}_d \otimes \mbox{UHF}_d)''$ by Proposition 3.4 as $\omega$ is pure. $\omega$ being also a factor state, we have $<f,V^ng> \raro <f,\Omega><\Omega,g>$ as $n \raro + \mbox{or} - \infty$ for any $f,g \in \pi(\clb_{loc})\Omega$ by Power's criteria [Pow] given in (1). Since such vectors are dense in the Hilbert space topology and the family $\{V^n:n \ge 1\}$ is uniformly bounded, we get $V^n \raro |\Omega><\Omega|$ in weak operator topology as $n \raro + \mbox{or} - \infty$.   

\vsp 
For the time being we assume that $H$ is trivial. Otherwise the argument that follows here we can use for the representation $\pi_0$ of $\tilde{\mbox{UHF}}_d \otimes \mbox{UHF}_d$ i.e. $\pi$ restricted to $[\pi(\tilde{\mbox{UHF}}_d \otimes \mbox{UHF}_d)\Omega]$. 

\vsp 
We have following distinct cases: 

\vsp 
\NI {\bf Case 1.} $\cle \neq I \;( \tilde{\cle} \neq I)$. Let $\cle_n \raro \cle_{-\infty}$ as $n \raro -\infty$ and thus $V\cle_{-\infty}V^*=\cle_{-\infty}$. We claim that either $\cle_{-\infty} = |\Omega><\Omega|$ or $\cle_{-\infty}$ is a proper infinite dimensional projection i.e. if $\cle_{- \infty}$ is a finite projection then $\cle_{- \infty}=|\Omega><\Omega|$. Suppose not then the finite subspace is shift invariant. In particular there exists a unit  vector $f$ orthogonal to $\Omega$ such that $Vf=zf$ for some $z \in S^1$ and this contradicts weak mixing property i.e. $V^n \raro |\Omega><\Omega|$ in weak operator topology proved above as point spectrum of $V$ has only $1$ with spectral multiplicity $1$. 

\vsp 
If $\cle_{-\infty}$ is infinite dimensional we can get a unitary operator $U_0$ from 
$F_0=\tilde{\clh} \otimes_{\clk} \clh$ onto $\cle_{- \infty}$ and 
via the unitary map we can get a sequence of increasing projections $U_0\cle_nU^*_0$ in $\cle_{- \infty}$ and note that $U_0\cle_nU_0^*=V^nU_0 \cle U_0^*(V^n)^*$. Note that if $\cle_{- \infty}$ is infinite dimension the process will not stop in finite step. Thus we have $F_0 \ominus \Omega = \oplus_{1 \le k \le n_{\cle}} F(k)$ where the index set is either singleton or infinity and each $F(k)$ will give a system of imprimitivity with respect to $V$, where $F(1)
=F_0-\cle_{-\infty}$. Further $\tilde{\mbox{UHF}}_d$ being a simple $C^*$-algebra, each such imprimitivity system 
is of Mackey index $\aleph_0$ [Mo3,section 4]: We fix a nonzero $f \in \cle-\theta^{-1}(\cle) \neq 0$ otherwise $\cle=I$ as $\theta^n(\cle) \uparrow I$ as $n \raro \infty$. $\pi_f:x \raro \theta^{-1}(x)f$ gives a representation of $\tilde{\mbox{UHF}}_d=\clb_L$ and we check that $[\pi_f(\theta^{-1}(\clb_L)''f] \le \cle-\theta^{-1}(\cle)$ 
as $f \perp [\theta^{-1}(\pi(\clb_R)')\Omega] \ge [\theta^{-1}(\pi(\clb_L))\Omega]$. Thus simplicity 
ensures that $\cle - \theta^{-1}(\cle)$ is a projection of dimension $\aleph_0$. Further $n_\cle$ is either $1$
or $\aleph_0$ since $F_0$ is separable.   

\vsp 
Since $\tilde{\clf}$ is also a proper projection, same argument is valid for $\tilde{\clf}$ with $\tilde{\clf}_{-\infty}=
\mbox{lim}_{n \raro -\infty} \theta^n(\tilde{\clf})$ i.e. we can write $F_0 \ominus \Omega = \oplus_{1 \le k \le n_{\tilde{\clf}}} \tilde{H}(k)$, where each $\tilde{H}(k)$ give rises to a system of imprimitivity with respect to $V$ where each system of imprimitivity is of Mackey index $\aleph_0$ where $\tilde{F}(1) = F_0-\tilde{\clf}_{-\infty}$ 
and $n_{\tilde{\clf}}$ is either $1$ or $\aleph_0$.        

\vsp 
In the following we use temporary notation $H$ for Hilbert subspace $F_0$. For a cardinal number $n$, we 
amplify a representation $\pi:\clb \raro \clb(H)$ of the $C^*$ algebra $\clb$ to $n$ fold direct sum 
$n\pi=\oplus_{1 \le k \le n } \pi_k $ acting on $nH=\oplus_{1 \le k \le n } H_k$  defining by 
$$n\pi(x) (\oplus \zeta_k) = \oplus (\pi(x)\zeta_k)$$ 
where $\pi_k = \pi$ is the representation of $\clb=\mbox{UHF}_d \otimes \tilde{\mbox{UHF}}_d$ 
on $H_k=H$ where $H =[\pi(\tilde{\mbox{UHF}}_d \otimes \mbox{UHF}_d)\Omega]$. We also 
extend $\bar{\tilde{F}} = \oplus \tilde{F}_{\alpha}$, $\bar{E} = \oplus E_{\alpha}$ and $\bar{V}=\oplus_{1 \le k \le n} V_k$ respectively. We also set notation $\Omega_k=\oplus_{1 \le k \le n}\delta^k_j\Omega$. 

\vsp 
Thus by Mackey's theorem, there exists a cardinal number $n \in \aleph_0$ and a unitary operator $U: nH \raro nH$ so that 
$\bar{V}=U\bar{V}U^*$ and $\bar{\cle}=U\bar{\tilde{\clf}}U^*$. We set a representation 
$\pi^U:\clb \raro \clb(nH)$ by $\pi^U(x)= Un\pi(x)U^*$
and rewrite the above identity as 
$$\oplus_{1 \le k \le n}[\pi_k(\mbox{UHF}_d)'\Omega_k] = \oplus_{1 \le k \le n}
[\pi^U_k(\tilde{\mbox{UHF}}_d)''\Omega_k]$$
where $\pi^U_k(x)=U\pi_k(x)U^*$. Note that by our construction we 
can ensure $U\Omega_k=\Omega_k$ for all $1 \le k \le n$ as the operator 
intertwining between two imprimitivity systems are acting on the orthogonal 
subspace of the projection generated by vectors $\{\Omega_k:1 \le k \le n\}$. 

\vsp 
We claim $\cle=\tilde{\clf}$. Suppose not i.e. $\tilde{\clf} < \cle$. In such a case we 
have
$$\oplus_{1 \le k \le n}[\pi_k(\mbox{UHF}_d)'\Omega_k] < \oplus_{1 \le k \le n}
[\pi^u_k(\mbox{UHF}_d)'\Omega_k]$$
Alternatively 
$$\oplus_{1 \le k \le n}[\pi_k(\tilde{\mbox{UHF}}_d)''\Omega_k] < \oplus_{1 \le k \le n}
[\pi^u_k(\tilde{\mbox{UHF}}_d)''\Omega_k]$$

\vsp 
Thus in principle we can repeat our construction now with $\pi^U$ and so we 
get a strict partial ordered set of quasi-equivalent representation of $\clb$. 
In the following we now aim to employ formal set theory to bring a contradiction 
on our starting assumption that $\tilde{\clf} < \cle$. 

\vsp 
To that end we need to deal with more then one representation of $\clb$. For the rest of the proof we reset notation $\pi_0$ for $\pi$ used for the pure representation of 
$\clb$ in $H_0 =[\pi_0(\clb)\Omega_0]$ where $\Omega_0$ is the cyclic vector, the reset notation for $\Omega$. Let $\clp$ be the collection of representation $(\pi,H_{\pi},\Omega)$ 
quasi-equivalent to $\pi_0:\clb \raro \clb(H_0)$ with a shift invariant vector state $\omega(x)=<\Omega,\pi(x)\Omega>$ i.e. $\omega(\pi(\theta(x))=\omega(\pi(x))$. So there exists 
minimal cardinal numbers $n_{\pi},n_0(\pi)$ so that $n_{\pi}H_{\pi}$ is unitary equivalent to $n_0(\pi)\pi_0$. Thus for such an element $(\pi,H_{\pi},\Omega_{\pi})$ we can associate 
two cardinal numbers $n_{\pi}$ and $n_0(\pi)$ and without loss of generality we assume that $H_{\pi} \subseteq n_0(\pi)H_0$ and $n_{\pi}H_{\pi}=n_0(\pi)H_0.$  $\pi_0$ being a pure 
representation, any element $\pi \in \clp$ is a type-I factor representation of $\clb$. The interesting point here that $\dsp{\oplus_{\pi \in \clp }}\pi$ is also an element in $\clp$ 
with associated cardinal numbers $\sum_{\pi} n_{\pi}$ and $\sum_{\pi} n_0(\pi)$. We say $(\pi_1,H_{\pi_1},\Omega^1) \prec (\pi_2,H_{\pi_2},\Omega^2)$ if there exists an 
isometry $U:n_{\pi_1} H_{\pi_1} \raro n_{\pi_2} H_{\pi_2}$ so that 

\vsp 
\NI (C1) For each $1 \le \alpha \le n_{\pi_1} $ we have 
$U \Omega^1_{\alpha} = \Omega^2_{\alpha'}$ for some 
$1 \le \alpha' \le n_{\pi_2}$;        

\NI (C2) $n_{\pi_2}\pi_2(x)E'_2 = U n_{\pi_1}\pi_1(x) U^* $ where $\cle'_2 \in n_{\pi_2}\pi_2(\clb)'$; 
 
\NI (C3) $\oplus_{1 \le \alpha \le n_{\pi_1} }[\pi^{\alpha}_1(\mbox{UHF}_d)'\Omega^1_{\alpha}] < \oplus_{1 \le 
\alpha \le n_{\pi_2}}[\pi^{\alpha}_2(\mbox{UHF}_d)'\Omega^2_{\alpha}]E'_2$.

\vsp 
In the inequality we explicitly used that both Hilbert spaces are subspaces of $nH_0$ for some possibly larger cardinal 
number $n$. That the partial order is non-reflexive follows as $(\pi,H,\Omega) \prec (\pi,H,\Omega)$ contradicts 
(C3) as $I=E_2'$. Partial order property follows easily. If $\pi_1 \prec \pi_2$ and $\pi_2 \prec \pi_3$ then $\pi_1 \prec \pi_3$. 
If $U_{12}$ and $U_{23}$ are isometric operators that satisfies (C1)-(C3) respectively, then $U_{13}=U_{23}U_{12}$ will 
do the job for $\pi_1$ and $\pi_3$. Thus $\pi^U \in \clp$ and by our starting assumption that $\tilde{\clf} \neq \cle$ we also 
check that $\pi_0 \prec \pi^U$. Thus going via the isomorphism we also check that for a given element $\pi \in \clp$ there exists 
an element $\pi' \in \clp$ so that $\pi \prec \pi'$. Thus $\clp_0$ is a non empty set and has at least one infinite chain containing $\pi_0$.
 
\vsp 
However by Hausdorff maximality theorem there exists a non-empty maximal totally ordered subset $\clp_0$ of $\clp$ 
containing $\pi_0$. We claim that $\pi_{max} = \oplus_{\pi \in \clp_0} \pi$ on $H_{\pi_{max}}=\oplus_{\pi \in \clp_0} H_{\pi} $ is an upper bound in $\clp_0$. That $\pi_{max} \in \clp$ is obvious. Further given an element $(H_{\pi_1},\pi_1,\Omega_1) \in \clp_0$ there exists an element $(H_{\pi_2},\pi_2,\Omega_2) \in \clp_0$ so that $\pi_1 \prec \pi_2$ by our starting remark as $\pi_0 \prec \pi^U$. By extending isometry $U_{12}$ to an isometry from $H_{\pi_1} \raro n_{\pi_{max}} H_{\pi_{max}}$ trivially we get the required isometry that satisfies (C1),(C2) and (C3) where cardinal numbers $n_{\pi_{max}} = \sum_{\pi \in \clp_0} n_{\pi} \in \aleph$. Thus by maximal property of $\clp_0$ we have $\pi_{max} \in \clp_0$. This brings a contradiction as by our construction $(\pi_{max},H_{\pi_{max}},\Omega) \prec (\pi_{max},H_{\pi_{max}},\Omega)$ as $\pi_{max} \in \clp_0$ but partial order is strict. This contradicts our starting hypothesis that $\tilde{\clf} < \cle$. This completes the proof that $\tilde{\clf}=\cle$ when $\cle \neq I$. By symmetry of the argument we also get $\clf=\tilde{\cle}$ when $\tilde{\cle} < 1$.  

\vsp 
\NI {\bf Case 2:} $\cle=I \;(\tilde{\cle}=I)$. We need to show $\tilde{\clf}=I\;(\clf=I)$ respectively. Suppose not and assume that both $\tilde{\clf}$ is a proper non-zero projection. 

\vsp 
We set projection $G$ on the closed linear span of elements in the subspaces $[\theta^{-n}(\clf)\pi(\tilde{\mbox{UHF}}_d)''\Omega]$ for all $n \ge 0$. We recall that $\theta(X)=VXV^*$ 
where $V=\sum_k S_k\tilde{S}_k^*$ and $\theta^{-1}(X)=\tilde{\Lambda}(X)$ 
for $X \in \pi(\tilde{\mbox{UHF}}_d)''$. Thus we have 
$$V^*\theta^{-n}(\clf)\pi(\tilde{\mbox{UHF}}_d)''\Omega$$
$$=\theta^{-n-1}(\clf)V^*\pi(\tilde{\mbox{UHF}}_d)''V\Omega$$
$$=\theta^{-n-1}(\clf) \tilde{\Lambda}(\pi(\tilde{\mbox{UHF}}_d)''\Omega.$$
Thus $(1-G)V^*G=0$ i.e. $\theta(G) \ge G$. It is also clear that $\tilde{\clf} \le G$ 
as the defining sequence of subspaces of $G$ goes to precisely $\tilde{\clf}$ as $n \raro \infty$ ( recall 
that $\theta^{-n}(\clf)=\tilde{\Lambda}_n(\clf) \uparrow I$ strongly as $n \uparrow \infty$ ). Once more we have 
$\theta^n(G) \ge \theta^n(\tilde{\clf}) = \Lambda^n(\tilde{\clf}) \uparrow I$ as $n \uparrow 
\infty$.  

\vsp 
If $G$ is a proper projection we can follow the steps as in the case 1 to find a unitary operator 
$U:nH \raro nH$ with $U\Omega_k=\Omega_k$ and $U\bar{V}U^*=\bar{V}$ so that $U\bar{G}U^*=\bar{\tilde{\clf}}$. 
We consider the subset $\clp_G$ of elements in $\clp$ for which $\cle_{\pi}=1$ and 
$\{\theta^{-n}(\clf_{\pi}):n \ge 0 \}$ commutes with $\tilde{\clf}_{\pi}$ 
and modify the strict partial ordering by modifying (C3) as  

\vsp 
\NI (C3') $\oplus_{1 \le \alpha \le n_{\pi_1} }[\pi^{\alpha}_1(\tilde{\mbox{UHF}}_d)''\Omega^1_{\alpha}] <  \oplus_{1 \le 
\alpha \le n_{\pi_2}}[\pi^{\alpha}_2(\tilde{\mbox{UHF}}_d)''\Omega^2_{\alpha}]E'_2$

\NI So we also get $\pi^U \in \clp_G$ and $\pi_0 \prec \pi^U$ and going along the same line we conclude that 
$G=\tilde{\clf}$. Thus we conclude that $G$ is either equal to $1$ or $G=\tilde{\clf}$. 

\vsp 
\NI {\bf Sub-case 1 of case 2:} If $G=I$ then $\clf G = \clf$ and so $[\clf\pi(\tilde{\mbox{UHF}}_d)''\Omega]=\clf$ as $\theta^{-n}(\clf) \ge \clf$. 
Thus $\tilde{\clf} \ge \clf$. So $\tilde{\clf} \ge \tilde{\Lambda}^n(\clf)$ for all $n \ge 1$ and taking limit we get $\tilde{\clf} \ge I$ i.e. $\tilde{\clf}=I$. 
This contradicts our starting assumption that $\tilde{\clf}$ is a proper projection.  

\vsp 
\NI {\bf Sub-case 2 of case 2:} Now we consider the case $G=\tilde{\clf} < I$. In such a case we have 
$(1- \tilde{\clf})\theta^{-n}(\clf)\tilde{\clf}=0$ and so $\theta^{-n}(\clf)$ commutes with $\tilde{\clf}$ for all $n \ge 0$. 

\vsp 
First we rule out the simplest possibility in the present situation for $\clf \tilde{\clf} = |\Omega><\Omega|$. If so then $\omega$ is
a Bernoulli state and a proof follows once we compute the following using the property $\Lambda(\tilde{\clf}) \ge \tilde{\clf}$ i.e.
$\tilde{\clf}\pi(s_i)^*\tilde{\clf} = \pi(s_i)^*\tilde{\clf}$ 
and commuting property of $\clf$ with $\pi(s_i)^*$ to get some scalars $\bar{\lambda_i}$ such
that 
$$\bar{\lambda}_i= \clf\tilde{\clf}\pi(s_i)^*\tilde{\clf}\clf = \pi(s_i)^*\tilde{\clf}\clf$$
and so 
$$\omega(s_Is_J^*)=<\Omega,\pi(s_I) \pi(s_J)^*\Omega>$$
$$=<\Omega, \clf \tilde{\clf} \pi(s_I)\pi(s_J)^*\tilde{\clf}\clf \Omega>$$
$$=\lambda_I\lambda_J^*$$
The inductive limit state $\omega$ on $\clb$, give rises to a pure state once restricted to $\clb_R$ and thus we have 
$\cle=\tilde{\clf}$ by Haag duality when $\pi(\clb_R)''$ is a type-I factor. This contradicts our starting assumption that $\cle=I$.   

\vsp 
Now we set projection $\clf'$ defined by 
\be 
\clf'=\clf-\clf\tilde{\clf}+|\Omega><\Omega|
\ee 
and check by commuting property of $\tilde{\clf}$ with $\clf$ that 
$$ 
\clf'\theta^{-1}(\clf')\clf' = (I-\tilde{\clf})\clf\theta^{-1}(\clf)(I-\theta^{-1}(\tilde{\clf}))(I-\tilde{\clf})\clf 
+ |\Omega><\Omega|$$ 
$$= \clf(I-\tilde{\clf}) + |\Omega><\Omega|=\clf'$$ 
where we have used $\theta^{-1}(\clf) = \tilde{\Lambda}(\clf) \ge \clf,\; \theta(\clf) \le \clf$ and  
$\theta(\tilde{\clf}) = \Lambda(\tilde{\clf}) \ge \tilde{\clf}$. 
Thus we get 
\be 
\theta^{-1}(\clf') \ge \clf'
\ee                                                                                                                                                                                                                                                                                                                                                                                                                                                                                                                     

\vsp 
We also rule out the possibility that $\clf \tilde{\clf} = \clf$. If so then $\clf \le \tilde{\clf}$ and $\tilde{\Lambda}^n(\clf) \le 
\tilde{\Lambda}^n(\tilde{\clf}) = \tilde{\clf}$. Taking limit we get $\tilde{\clf}=I$ as $\tilde{\Lambda}^n(\clf) \uparrow I$ as $n \raro \infty$. 
This brings a contradiction to our hypothesis.  
  
\vsp 
So we have in particular $\clf' <  \clf \le I$ and $\clf' - |\Omega><\Omega| \neq 0$. Now we will rule out the possibility of $\clf=I$ under our hypothesis $\tilde{\clf} < \cle=I$. 
Suppose so i.e. $\clf=I$, then $\tilde{\cle}=I$ since $\tilde{\cle} \ge \clf$. Then $Q=\cle \tilde{\cle} = \cle = \cle \clf=P$. Thus by Proposition 3.5 (e) we get $\pi(\clo_d)'= 
\pi(\tilde{\clo}_d)''$ and so we have in particular $\tilde{\clf}=[\pi(\tilde{\clo}_d)''\Omega]=[\pi(\clo_d)'\Omega]=\cle$. This contradicts our starting assumption once more that $\tilde{\clf} < \cle=I$.    

\vsp 
We also have $\theta^{-1}(\clf') \ge \clf'$ and $\theta^{-1}(\clf) \ge \clf$. Thus we can follow the steps of Case-1 with elements $\clf',\clf,\theta^{-1}$ replacing the role of $\tilde{\clf},\cle,\theta$ 
to get a unitary operator $U:nH \raro nH$ so that $U\bar{V}=\bar{V}U$ and $U\bar{\clf}U^*=\bar{\clf'}$ for a 
cardinal number $n$.    

\vsp 
Now we consider a further subset $\clp_{G'}$ of $\clp_G$ consist of quasi-equivalent representations $\pi$ to $\pi_0$ 
of $\clb$ where $\pi$ admits the additional property: $\clf_{\pi} < I, \cle_{\pi}=I$ and $\{\theta^{-n}(\tilde{\clf_{\pi}}): n \ge 0 \}$ 
commutes with $\clf_{\pi}$ with the strict partial ordering $(\clh_{\pi_1},\pi_1,\Omega_{\pi_1}) \prec (\clh_{\pi_2},\pi_2,\Omega_{\pi_2})$ 
given by modifying condition (C3') as

\vsp 
\NI (C3'') $\oplus_{1 \le \alpha \le n_{\pi_1} }[\pi^{\alpha}_1(\mbox{UHF}_d)''\Omega^1_{\alpha}] >  \oplus_{1 \le 
\alpha \le n_{\pi_2}}[\pi^{\alpha}_2(\mbox{UHF}_d)''\Omega^2_{\alpha}]E'_2$

\vsp 
Since $\pi^U$ also satisfies the conditions that of $\pi_0 \in \clp_{G'}$ by covariance relation 
of $U$ with respect to shifts once more we get $\pi^U \in \clp_{G'}$ and $\pi_0 \prec \pi^U$. Thus we can 
repeat the process and so $\clp_{G'}$ has at least one infinite chain of totally ordered containing $\pi_0$.   
Once more by Hausdorff maximality principle we bring a contradiction to our starting assumption 
that $\clf' < \clf$. In other words this brings a contradiction to our starting hypothesis 
that $\tilde{\clf}$ is a proper projection i.e. $\clf < \cle=I$. Thus we arrive at $\tilde{\clf}=\cle$ 
when $\cle=I$. 

\vsp 
By symmetry of argument used here it also follows that $\clf=\tilde{\cle}$ when $\tilde{\cle}=I$. This completes the proof of 
$\tilde{\clf}=\cle,\;\; \clf =\tilde{\cle}$ for the case when $H$ is the trivial closed subgroup of $S^1$. 

\vsp 
Now we will remove the assumption that $H$ is trivial using Proposition 3.4. Let $(H,\pi_0,\Omega)$ be the GNS space of the state $\omega$ on $\clb$. 
Let $e_0$ and $\tilde{e}_0$ be the support projections of $\omega$ in $\pi_0(\clb_R)''$ and $\pi_0(\clb_L)''$ respectively. Similarly we also 
set projections $\clf_0=[\pi_0(\cla_R)''\Omega]$ and $\tilde{\clf}_0=[\pi_0(\clb_L)''\Omega]$. We also set projections 
$q_0=[\pi_0(\mbox{UHF})'\Omega] [\pi_0(\tilde{\mbox{UHF}}_d)'\Omega]$ and $p_0=[\pi_0(\mbox{UHF})'\Omega] [\pi_0(\mbox{UHF}_d)''\Omega]$.

\vsp 
$\omega$ being pure, by Theorem 3.4 any $\{\beta_z:z \in H \}$ invariant element of $\clb(\tilde{\clh} \otimes_{\clk} \clh)$ is an element in 
$\pi(\tilde{\mbox{UHF}}_d \otimes \mbox{UHF}_d)''$. $\cle,\tilde{\cle},\clf, \tilde{\clf}$ are $\{\beta_z:z \in H \}$ elements. Thus once 
we identify cyclic space $[\pi_0(\clb)\Omega]$ with $F_0=[\pi(\tilde{\mbox{UHF}}_d \otimes \mbox{UHF}_d)''\Omega]$, we get obvious relations 
\be 
\cle F_0=e_0,\tilde{\cle}F_0=\tilde{e}_0,\clf F_0=f_0, \tilde{\clf}F_0=\tilde{f}_0
\ee
and 
\be 
P F_0 = p_0 \;\;\mbox{and} \;\; QF_0=q_0
\ee
as 
$\cle=[\pi(\mbox{UHF})'\Omega]$, $\tilde{\cle} = [\pi(\tilde{\mbox{UHF}}_d)'\Omega]$, $Q=\cle\tilde{\cle}$ and $P=\cle \clf$. Further $V$ is also 
$\{ \beta_z:z \in H \}$ invariant and $V \pi(x)V^*=\pi(\theta(x))$ for all $x \in \clb$ which we have identified with $\tilde{\mbox{UHF}}_d \otimes \mbox{UHF}_d$. 

\vsp 
By applying the first part of the argument with representation $\pi_0$, for pure $\omega$, we have 
$e_0=\tilde{f}_0$ and $\tilde{e}_0=f_0$ and $p_0=q_0$.

\vsp 
Now we write the equality $p_0=q_0$ as $\cle \clf F_0=\cle\tilde{\cle}F_0$ and apply $\Lambda$ on both side to 
conclude that $\Lambda(\cle)\clf F_1=\Lambda(\cle)\tilde{\cle} F_1$ and multiplying by $\cle$ from left we get 
$\cle\clf F_1=\cle\tilde{\cle}F_1$ as $\Lambda(\cle)\cle=\cle$ and thus we get $PF_1=QF_1$. By repeated 
application of $\Lambda$, we get $PF_m=QF_m$. 

\vsp 
If $H=\{z:z^n=1 \}$ then we get $P=\sum_k PF_k=\sum_k QF_k=Q$. This completes the proof for $P=Q$. Similarly 
$\tilde{\clf} = \sum_k \tilde{\clf}F_k = \sum_k \cle F_k = \cle$ and also $\clf=\tilde{\cle}$. 

\vsp 
If $H=S^1$ then $\hat{H}=\IZ$ and for $m \ge 0$, we have $PF_m=QF_m$. For $m < 0$ we take $k=-m$ and  
check that $\Lambda^k(QF_m-PF_m)=\Lambda^k(Q)F_{m+k}-\Lambda^k(P)F_{m+k} = 
\Lambda^k(\cle)\tilde{\cle} F_0 - \Lambda^k(\cle)\clf F_0= \Lambda^k(\cle)(\tilde{e}_0 - f_0)=0$
Since $\Lambda$ is an injective map, we get $QF_m=PF_m$ for all $m < 0$.

\vsp
Now we are left to prove those three statements given in (b). $\omega$ being pure we have $P=Q$ and thus 
by Proposition 3.5 we have $\tilde{\clm}=\clm'$ and $\pi(\clo_d)'=\pi(\tilde{\clo}_d)''$. We are left to show
$\pi_{\omega}(\clb_R)'=\pi_{\omega}(\clb_L)''$. For that we recall $F_0$ and check few obvious relation 
$F_0\pi(\clo_d)'F_0=F_0\pi(\tilde{\clo}_d)''F_0$ and $\pi_{\omega}(\clb_R)' \subseteq F_0\pi(\clo_d)'F_0$.
Since $F_0\pi(\tilde{\mbox{UHF}}_d)''F_0$ is equal to $(\beta_z:z \in H)$ invariant elements in 
$F_0\pi(\tilde{\clo}_d)''F_0$ and elements in $\pi_{\omega}(\clb_R)'$ are $(\beta_z: z \in H)$ invariant 
we conclude that $\pi_{\omega}(\clb_R)' \subseteq \pi_{\omega}(\clb_L)''$. Inclusion in other direction 
is obvious and thus Haag duality property (iii) holds. 
\end{proof}

\vsp 
\begin{proof} 
(of Theorem 1.1:) (a) implies that $q_0=p_0$. (d) also says that $p_0=q_0$. Thus in either case, following last part of the proof of Theorem 3.6 we get $P=Q$. The statement (e) also implies $P=Q$. 
That shows now that (a),(d) as well as (e) implies (f) by the if part of Theorem 2.6 (a). That (b) implies (a) is trivial as $\Omega$ is separating for $\clm_1,\tilde{\clm}_1$ by faithful 
property. That (c) implies (f) is trivial as $\pi_{\omega}(\clb_R)''$ is a factor. Thus we have showed so far any of the statement (a),(b),(c),(d),(e) implies (f). For the converse we appeal to
the only if part of Theorem 3.6. 
\end{proof}

\bigskip
{\centerline {\bf REFERENCES}}

\begin{itemize} 
\bigskip
\item{[Ac]} Accardi, L. : A non-commutative Markov property, (in Russian), Functional.  anal. i Prilozen 9, 1-8 (1975).

\item{[AcC]} Accardi, Luigi; Cecchini, Carlo: Conditional expectations in von Neumann algebras and a theorem of Takesaki.
J. Funct. Anal. 45 (1982), no. 2, 245–273. 

\item{[AM]} Accardi, L., Mohari, A.: Time reflected
Markov processes. Infin. Dimens. Anal. Quantum Probab. Relat. Top., vol-2
,no-3, 397-425 (1999).

\item{[AKLT]} Affleck, L.,Kennedy, T., Lieb, E.H., Tasaki, H.: Valence Bond States in Isotropic Quantum Antiferromagnets, Commun. Math. Phys. 115, 477-528 (1988). 

\item{[AHP]} Akiho, Nobuyuki; Hiai, Fumio; Petz, Dénes Equilibrium states and their entropy densities in gauge-invariant C∗-systems. Rev. Math. Phys. 17 (2005), 
no. 4, 365-389.

\item{[Ara1]} Araki, H.:  Gibbs states of a one dimensional quantum lattice. Comm. Math. Phys. 14 120-157 (1969). 

\item{[Ara2]} Araki, H.: On uniqueness of KMS-states of one-dimensional quantum lattice systems, Comm. Maths. Phys. 44, 1-7 (1975).

\item{[AMa]} Araki, H., Matsui, T.: Ground states of the XY model, Commun. Math. Phys. 101, 213-245 (1985).

\item{[Ar1]} Arveson, W.: On groups of automorphisms of operator algebras, J. Func. Anal. 15, 217-243 (1974).

\item{[Ar2]} Arveson, W.: Pure $E_0$-semigroups and absorbing states. 
Comm. Math. Phys 187 , no.1, 19-43, (1997)

\item{[BBN]} Baumgartner, Bernhard; Benatti, Fabio; Narnhofer, Heide: Translation invariant states on twisted algebras on a lattice. J. Phys. A 43 (2010), no. 11, 115301, 13 pp. 

\item{[BR]} Bratteli, Ola,: Robinson, D.W. : Operator algebras
and quantum statistical mechanics, I,II, Springer 1981.

\item{[BJP]} Bratteli, Ola,: Jorgensen, Palle E.T. and Price, G.L.: 
Endomorphism of $\clb(\clh)$, Quantisation, nonlinear partial differential 
equations, Operator algebras, ( Cambridge, MA, 1994), 93-138, Proc. Sympos.
Pure Math 59, Amer. Math. Soc. Providence, RT 1996.

\item{[BJKW]} Bratteli, Ola,: Jorgensen, Palle E.T., Kishimoto, Akitaka and
Werner Reinhard F.: Pure states on $\clo_d$, J.Operator Theory 43 (2000),
no-1, 97-143.    

\item{[BJ]} Bratteli, Ola: Jorgensen, Palle E.T. Endomorphism of $\clb(\clh)$, II, 
Finitely correlated states on $\clo_N$, J. Functional Analysis 145, 323-373 (1997). 

\item{[Cun]} Cuntz, J.: Simple $C\sp*$-algebras generated by isometries. Comm. Math. Phys. 57, 
no. 2, 173--185 (1977).  

\item{[DHR]} Doplicher,S., Haag, R. and Roberts, J.: Local observables and particle statistics I, II, Comm. Math. Phys. 23 (1971) 119-230 and 35, 49–85 (1974)

\item{[Ev]} Evans, D.E.: Irreducible quantum dynamical
semigroups, Commun. Math. Phys. 54, 293-297 (1977).

\item{[Ex]} Exel, Ruy : A new look at the crossed-product of a $C^*$-algebra by an endomorphism. (English summary)
Ergodic Theory Dynam. Systems 23 (2003), no. 6, 1733–1750.

\item{[FNW1]} Fannes, M., Nachtergaele, B., Werner, R.: Finitely correlated states on quantum spin chains,
Commun. Math. Phys. 144, 443-490(1992).

\item{[FNW2]} Fannes, M., Nachtergaele, B., Werner, R.: Finitely correlated pure states, J. Funct. Anal. 120, 511-
534 (1994). 

\item{[FNW3]} Fannes, M., Nachtergaele, B. Werner, R.: Abundance of translation invariant states on quantum spin chains, Lett. Math. Phys. 25 no.3, 249-258 (1992).

\item{[Fr]} Frigerio, A.: Stationary states of quantum dynamical semigroups, Comm. Math. Phys. 63 (1978) 269-276.

\item{[Hag]} Haag, R.: Local quantum physics, Fields, Particles, Algebras, Springer 1992. 

\item{[HMHP]} Hiai, Fumio; Mosonyi, Milán; Ohno, Hiromichi; Petz, Dénes, Free energy density for mean field perturbation of states of a one-dimensional spin chain. Rev. Math. Phys. 20 (2008), no. 3, 335-365. 

\item{[Ki]} Kishimoto, A.: On uniqueness of KMS-states of one-dimensional quantum lattice systems, Comm. Maths. Phys. 47, 167-170 (1976).

\item{[La]} Christopher, Lance, E.: Ergodic theorems for convex sets and operator algebras, Invent. Math. 37, no. 3, 201-214 (1976). 

\item{[Mac]} Mackey, George W.: Imprimitivity for representations of locally compact groups. I. Proc. Nat. Acad. Sci. U. S. A. 35, 
537-545 (1949). 

\item{[Ma1]} Matsui, A.: Ground states of fermions on lattices, Comm. Math. Phys. 182, no.3 723-751 (1996). 

\item{[Ma2]} Matsui, T.: A characterization of pure finitely correlated states. 
Infin. Dimens. Anal. Quantum Probab. Relat. Top. 1, no. 4, 647--661 (1998).

\item{[Ma3]} Matsui, T.: The split property and the symmetry breaking of the quantum spin chain, Comm. 
Maths. Phys vol-218, 293-416 (2001) 

\item{[Ma4]} Matsui, Taku, On the absence of non-periodic ground states for the antiferromagnetic XXZ model. Comm. Math. Phys. 253 (2005), no. 3, 585-609.

\item{[Na]} Nachtergae, B. Quantum Spin Systems after DLS1978, `` Spectral Theory and Mathematical Physics: A Festschrift in Honor of Barry Simon's 60th Birthday '' Fritz Gesztesy et al. (Eds), 
Proceedings of Symposia in Pure Mathematics, Vol 76, part 1, pp 47--68, AMS, 2007.

\item{[Mo1]} Mohari, A.: Markov shift in non-commutative probability, J. Funct. Anal. vol- 199 
, no-1, 190-210 (2003) Elsevier Sciences. 

\item{[Mo2]} Mohari, A.: Pure inductive limit state and Kolmogorov's property, J. Funct. Anal. vol 253, no-2, 584-604 (2007)
Elsevier Sciences.

\item{[Mo3]} Mohari, A: Jones index of a completely positive map, Acta Applicandae Mathematicae. Vol 108, Number 3, 665-677 
(2009).  

\item{[Mo4]} Mohari, A.: Pure inductive limit state and Kolmogorov's property-II,  http://arxiv.org/abs/1101.5961 To appear Journal of Operator Theory. 

\item{[Mo5]} Mohari, A.: Translation invariant pure states in quantum spin chain and its split property, http://arxiv.org/abs/0904.2104. 

\item{[Mo6]} Mohari, A.: A complete weak invariance for Kolmogorov states on $\clb = \dsp{\otimes_{k \in \IZ}}\!M^{(k)}_d(\IC)$, http://arxiv.org/abs/

\item {[OP]} Ohya, M., Petz, D.: Quantum entropy and its use, Text and monograph in physics, Springer 1995. 

\item{[Or1]} Ornstein, D. S.: Bernoulli shifts with the same entropy are isomorphic, Advances in Math. 4 1970 337-352 (1970).

\item{[Or2]} Ornstein, D. S.: A K-automorphism with no square root and Pinsker's conjecture, 
Advances in Math. 10, 89-102. (1973).

\item{[Pa]} Parry, W.: Topics in Ergodic Theory, Cambridge University Press, 1981. 

\item{[Po]} Popescu, G.: Isometric dilations for infinite sequences of non-commutating operators, Trans. Amer. Math.
Soc. 316 no-2, 523-536 (1989)

\item{[Pow1]} Powers, Robert T.: Representations of uniformly hyper-finite algebras and their associated von Neumann. rings, Annals of Math. 86 (1967), 138-171.

\item{[Ru]} Ruelle, D. : Statistical Mechanics, Benjamin, New York-Amsterdam (1969) . 

\item{[Sa]} Sakai, S. : Operator algebras in dynamical systems. The theory of unbounded derivations in $C^*$-algebras. 
Encyclopedia of Mathematics and its Applications, 41. Cambridge University Press, Cambridge, 1991. 

\item{[So]} Stormer E.: On projection maps of von Neumann algebras, Math. Scand. 30, 46-50 (1972). 

\item{[Ta]} Takesaki, M.: Conditional Expectations in von Neumann Algebras, J. Funct. Anal., 9, pp. 306-321 (1972)

\end{itemize}

\end{document}